\documentclass[a4paper,12pt]{article}
\usepackage[utf8]{inputenc}
\usepackage[english]{babel}

\usepackage{amscd}
\usepackage{amssymb}
\usepackage{amsthm}
\usepackage{amsmath}
\usepackage{amstext}
\usepackage{amsxtra}
\usepackage{amsfonts}
\usepackage{color}
\usepackage{csquotes}
\usepackage{enumerate}
\usepackage{graphicx}
\usepackage{graphics}
\usepackage{epstopdf}
\usepackage[pdftex]{hyperref}
\usepackage{indentfirst}
\usepackage{lmodern}
\usepackage{listings}
\usepackage{multicol}
\usepackage{mathrsfs}
\usepackage{multirow}
\usepackage[mathscr]{euscript}
\usepackage{mathtools}
\usepackage{pdflscape}
\usepackage{pgf,tikz}
\usepackage{pdfpages}
\usepackage{pb-diagram}
\usepackage{pgothic}
\usepackage{rotating}
\usepackage[scr]{rsfso}
\usepackage{setspace}
\usepackage[shortlabels]{enumitem}
\usepackage{stmaryrd}
\usepackage{textcomp}
\usepackage[normalem]{ulem}
\usepackage{upgreek}
\usepackage{verbatim}
\usepackage[all]{xy}

\newcommand{\R}{\mathbb R}
\newcommand{\C}{\mathbb C}

\newcommand{\N}{\mathbb N}

\newcommand{\Z}{\mathbb Z}
\newcommand{\T}{\mathbb T}
\newcommand{\F}{\mathbb{F}}
\newcommand{\G}{\mathcal{G}}

\newcommand{\p}{\mathfrak{p}}

\usepackage{geometry}
\newtheorem{teorema}{Theorem}[section]
\newtheorem{lema}[teorema]{Lemma}
\newtheorem{corolario}[teorema]{Corollary}
\newtheorem{defi}[teorema]{Definition}
\newtheorem{prop}[teorema]{Proposition}
\newtheorem{exemplo}[teorema]{Example}
\newtheorem{exemplos}[teorema]{Examples}
\newtheorem{obs}[teorema]{Remark}
\theoremstyle{definition}

\author{Felipe}

\begin{document}

\author{{Felipe A. Tasca\footnote{This study was financed in part by Coordenação de Aperfeiçoamento de Pessoal de Nível Superior - Brasil (CAPES) - Finance Code 001.}}, and Daniel Gon\c{c}alves\footnote{Partially supported by Conselho Nacional de Desenvolvimento Cient\'ifico e Tecnol\'ogico (CNPq) grant numbers 304487/2017-1 and 406122/2018-0 and Capes-PrInt grant number 88881.310538/2018-01 - Brazil.}}

\title{KMS states and continuous orbit equivalence for ultragraph shift spaces with sinks}

\maketitle

\begin{abstract}
We extend ultragraph shift spaces and the realization of ultragraph C*-algebras as partial crossed products to include ultragraphs with sinks (under a mild condition, called (RFUM2), which allow us to dismiss the use of filters) and we describe the associated transformation groupoid. Using these characterizations we study continuous orbit equivalence of ultragraph shift spaces (via groupoids) and KMS and ground states (via partial crossed products). 
\end{abstract}

\vspace{0.5pc}



\section{Introduction}

Mathematicians have long sought a definition for shift spaces over countable alphabets. The more direct (and most common) approach is to consider a countable alphabet $A$ with the discrete topology and the shift space $A^\N$ with the product topology, what yields a space that is not locally compact. This strongly differs from the finite alphabet case, where shift spaces are always compact. Motivated (among other things) by this difference, a few definitions of shift spaces over infinite alphabets has arisen over the years, see for example \cite{Fie, OTW}, (also in \cite{OTW} there is an excellent overview of the topic). The definition proposed in \cite{OTW} has connections with C*-algebras, and dynamical properties of this shift are studied in \cite{Mstep, GSS, GSS1} for example. Very recently, building from the ideas in \cite{OTW, webster}, the authors of \cite{danidani} introduce the notion of an ultragraph shift space, which is a generalization (to the infinite alphabet) of shifts of finite type (SFT are among the most important shifts in symbolic dynamics, see for example \cite{LindMarcus}). These shifts have interesting dynamics, as their chaotic behavior coincide with the chaotic behavior of shifts of finite type over finite alphabets (see \cite{danibruno, danibruno2}), and Curtis–Hedlund–Lyndon type results can be proved (see \cite{danimarc}). Furthermore, ultragraph shift spaces are showed in \cite{danidani} to have a strong connection with C*-algebras: if two ultragraphs have associated shift spaces that are conjugate, via a conjugacy that preserves length, then the associated ultragraph C*-algebras are isomorphic. 

One should note that the ultragraph shift spaces defined in \cite{danidani} have two restrictions: They are defined only for ultragraphs without sinks and the invariance of the associated C*-algebras is showed only for ultragraphs that satisfy Condition~(RFUM). While one could, a priori, think that lifting the restriction on having no sinks is a straightforward task this is not the case. In fact, when adding sinks to ultragraphs, or to graphs, many unforeseen difficulties arise. For example, the results describing topological full groups as invariant for continuous orbit equivalence of shift spaces (or groupoids) associated to graphs (or ultragraphs) are only valid for graphs (or ultragraphs) without sinks (see \cite{Ort_Nyl} and \cite{danielgillesdanie}), and there is no clear way to extend these results to include graphs and ultragraphs with sinks. Furthermore, when studying KMS states associated to graphs the existence of sinks induce new KMS states (see \cite{GM}). 
Therefore it is interesting to extend the definition of ultragraph shift spaces given in \cite{danidani} to include ultragraphs with sinks. This is one of our main goals in this paper.

Once we have extended the definition of ultragraph shift spaces to include sinks, we provide further evidence that these shift spaces are well connected with C*-algebras, namely we describe continuous orbit equivalence of ultragraph shift spaces in terms of isomorphism of the associated groupoids and in terms of isomorphism of the associated C*-algebras. Furthermore, using partial crossed product theory, we study the dynamics associated with these shifts: we describe the KMS and ground states associated to ultragraph C*-algebras, extending results in \cite{danielgilles} to include ultragraphs with sinks.

Let us explain the structure of the paper. After this introduction we include a section of preliminaries, where we briefly recall concepts regarding ultragraphs and set up notation. In Section~3 we describe the topological space that, with the shift map, will form the shift space associated with an ultragraph. This space is called the boundary ultrapath space and, under a mild condition (called Condition~(RFUM2)), that allows us to dismiss the use of filters, we prove that this space has a countable basis of compact open sets. We remark that for ultragraphs without sinks Condition~(RFUM2) is the same as Condition~(RFUM) in \cite{danidani}, and that all graphs (including the ones with sinks) satisfy Condition~(RFUM2) (therefore our results are a complete generalization of the results in \cite{CarLar, danidani}). Proceeding, in Section~4 we realize the C*-algebras associated to ultragraphs that satisfy Condition~(RFUM2) as partial crossed products. We use this realization in Section~5 to, via the transformation groupoid, describe ultragraph C*-algebras as groupoid C*-algebras. In Section~6 we use the characterization of ultragraph C*-algebras as groupoid C*-algebras to describe stabiliser-preserving, continuous orbit equivalence of ultragraph shift spaces in terms of isomorphism of the associated C*-algebras, in terms of the associated groupoids, and in terms of continuous orbit equivalence that preserve eventually periodic isolated points. We also characterize eventual conjugacy in this section. Finally, using the description of ultragraph C*-algebras as partial crossed products done in Section~4, we describe KMS states associated to C*-algebras of ultragraphs that satisfy Condition~(RFUM2). In particular we compute the KMS states associated to an ultragraph (with sinks) such that the ultragraph Leavitt path algebra is not isomorphic to any (graph) Leavitt path algebra nor to any algebraic Exel-Laca algebra.

\section{Preliminaries}
\label{cap_ultragrafo}

We start recalling concepts regarding ultragraphs, as defined in \cite{tom2003} (see also \cite{KMST,mu&mu}). 

\begin{defi}\label{def of ultragraph}
An \emph{ultragraph} is a quadruple $\mathcal{G}=(G^0, \mathcal{G}^1, r,s)$ consisting of two countable sets $G^0, \mathcal{G}^1$, a map $s:\mathcal{G}^1 \to G^0$, and a map $r:\mathcal{G}^1 \to P(G^0)\setminus \{\emptyset\}$, where $P(G^0)$ is the power set of $G^0$.
\end{defi}

A key object when studying ultragraphs are generalized vertices, which we define below. 

\begin{defi}\label{prop_vert_gen}
Let $\mathcal{G}$ be an ultragraph. Define $\mathcal{G}^0$ to be the smallest subset of $P(G^0)$ that contains $\{v\}$ for all $v\in G^0$, contains $r(e)$ for all $e\in \mathcal{G}^1$, and is closed under finite unions and nonempty finite intersections. Elements of $\G^0$ are called \emph{generalized vertices}.
\end{defi}

An useful description of the generalized vertices is given below.

\begin{lema}\label{description}\cite[Lemma~2.12]{tom2003}
If $\mathcal{G}$ is an ultragraph, then 
\begin{equation*}
\G^0 =
\left\{
\begin{array}{cc}
\displaystyle\bigcap_{e\in X_1} r(e)\cup \ldots \cup \bigcap_{e \in X_n} r(e) \cup F : & X_1,\ldots, X_n\text{ are finite subsets of } \mathcal{G}^1  \\
 & \text{ and }F\text{ is a finite subset of }G^0  \\
\end{array}
\right\}.
\end{equation*}
Furthermore, $F$ may be chosen to be disjoint from $\bigcap_{e \in X_1} r(e) \cup \ldots \cup \bigcap_{e \in X_n} r(e)$.
\end{lema}

\begin{defi}\label{defi_alg_ult}
Let $\mathcal{G}$ be an ultragraph. The \emph{ultragraph algebra} $C^*(\mathcal{G})$ is the universal C*-algebra generated by a family of partial isometries with orthogonal ranges $\{s_e:e\in \mathcal{G}^1\}$ and a family of projections $\{p_A:A\in \mathcal{G}^0\}$ satisfying
\begin{enumerate}
\item\label{p_Ap_B=p_{A cap B}}  $p_\emptyset=0,  p_A p_B=p_{A\cap B},  p_{A\cup B}=p_A+p_B-p_{A\cap B}$, for all $A,B\in \mathcal{G}^0$;
\item $s_es_e^*\leq p_{s(e)}$ for all $e\in \mathcal{G}^1$; and
\item\label{CK-condition} $p_v=\sum\limits_{s(e)=v}s_es_e^*$ whenever $0<\vert s^{-1}(v)\vert< \infty$.
\end{enumerate}
\end{defi}

Next we set up notation that will be used throughout the paper. This agrees with notation introduced in \cite{danidani} and \cite{mu&mu}. Let $\G$ be an ultragraph. A \emph{finite path} in $\G$ of is either an element of $\G^{0}$ or a sequence of edges $e_{1}\ldots e_{k}$ in $\G^{1}$ where $s\left(e_{i+1}\right)\in r\left(e_{i}\right)$ for $1\leq i<k$. If we write $\alpha=e_{1}\ldots e_{k}$, the length $\left|\alpha\right|$ of $\alpha$ is $k$. The length $|A|$ of a path $A\in\G^{0}$ is zero. We define $r\left(\alpha\right)=r\left(e_{k}\right)$ and $s\left(\alpha\right)=s\left(e_{1}\right)$. For $A\in\G^0$,
we set $r\left(A\right)=A=s\left(A\right)$. The set of finite paths in $\G$ is denoted by $\G^{\ast}$. 
An \emph{infinite path} in $\G$ is an infinite sequence of edges $\gamma=e_1e_2\ldots$ in $\prod \G^1$, where $s\left(e_{i+1}\right)\in r\left(e_i\right)$ for all $i\geq 1$. Keep of course $s(\beta)=s(e_1)$. The length $\left|\gamma\right|$ of $\gamma\in\p^\infty$ is defined to be $\infty$. The set of infinite paths in $\G$ is denoted by $\p^\infty$. A vertex $v$ in $G^0$ is called a \emph{sink} if $\left|s^{-1}\left(v\right)\right|=0$, it is called an \emph{infinite emitter} if $\left|s^{-1}\left(v\right)\right|=\infty$, and it is called a \emph{source} if $v\notin r(e)$ for all $e\in\G^1$. 
If a vertex is a sink or an infinite emitter, then it is called a \emph{singular} vertex. Otherwise, i.e, if $0<|s^{-1}(v)|<\infty$, then $v$ is called a \emph{regular} vertex. All of these nomenclatures were inherited from graph theory. 

For $n\geq1,$ we define $$\p^n:=\{(\alpha,A):\alpha\in\G^*,|\alpha |=n,A\in\G^0,A\subseteq r(\alpha)\}.$$ 
We specify that $\left(\alpha,A\right)=(\beta,B)$ if, and only if,
$\alpha=\beta$ and $A=B$. Letting $\p^0:=\G^0$ we define the \label{def_esp_ult} \emph{ultrapath space} associated with the ultragraph $\G$ to be $\p:=\displaystyle\bigsqcup_{n\geq 0}\p^n$.  The elements of $\p$ are called \emph{ultrapaths}. We embed the set of finite paths $\G^*$ in $\p$ by sending $\alpha$ to $(\alpha,r(\alpha))$. We define the length $\left\vert\left(\alpha,A\right)\right\vert$ of a pair $\left(\alpha,A\right)$ to be $\left\vert\alpha\right\vert$. Each $A\in\G^{0}$ is regarded as an ultrapath of length zero and can be identified with the pair $(A,A)$. Hence, we extend the range map $r$ and the source map $s$ to $\p$ by declaring that $r\left(\left(\alpha,A\right)\right)=A$, $s\left(\left(\alpha,A\right)\right)=s\left(\alpha\right)$ and $r(A,A)=r(A)=A=s(A)=s(A,A)$.

We concatenate elements in $\p$ in the following way: If $x=(\alpha,A)$ and $y=(\beta,B)$, with $|x|\geq 1, |y|\geq 1$, then $x\cdot y$ is defined if, and only if, $s(\beta)\in A$, in which case, $x\cdot y:=(\alpha\beta,B)$. 
Also we specify that:
	$$x\cdot y=
	\left\{
	\begin{array}{ll}
	x\cap y, & \text{if $x,y\in \G^0$ and if $x\cap y\neq\emptyset$}\\
	y, & \text{if $x\in\G^0$, $|y|\geq 1$, and if $s(y)\in x$}\\
	x_y, & \text{if $y\in\G^0$, $|x|\geq 1$, and if $r(x)\cap y\neq \emptyset$},
	\end{array}
	\right.$$  
where, if $x=\left(\alpha,A\right)$, $\left|\alpha\right|\geq1$ and if $y\in\G^{0}$, the expression $x_{y}$ is defined to be $\left(\alpha,A\cap y\right)$. Given $x,y\in\p$, we say that $x$ has $y$ as an initial segment, or that $y$ is an initial segment of $x$, if $x=y\cdot x^{\prime}$, for some $x^{\prime}\in\p$, with $s\left(x^{\prime}\right) \cap r\left(y\right)\neq\emptyset$. 
We can also concatenate ultrapaths in $\p$ with paths in $\p^\infty$. If $s(\beta)\in r(\alpha,A)=A$, we define $(\alpha,A).\beta=\alpha\beta\in \p^\infty$, where $(\alpha,A)\in \p$ and $\beta \in \p^{\infty}$. Furthermore, if $\alpha=A$ then $(A,A).\beta=\beta\in\p^\infty$. Of course if $s(\beta)\notin A$, $(\alpha,A).\beta$ is not defined. We say $x\in \p\sqcup \p^\infty$ has $y\in\p$ as a \emph{initial segment}\label{seg_ini} if $x=y.x'$, for some $x'\in\p$, with $s(x')\cap r(y)\neq \emptyset$.

\begin{obs} To simplify notation we omit the dot in the definition of concatenation, so that $x\cdot y$ will be denoted by $xy$.
\end{obs}


\section{The boundary ultrapath space}
	
	In this section we will define a topological space (boundary ultrapath space) associated to an ultragraph, which generalizes the boundary path space of a graph studied in \cite{BCW}, \cite{CarLar} and \cite{webster} for example. Our construction builds on the construction done in \cite{danidani} in order to include ultragraphs with sinks. We start with a couple of definitions.
	
\begin{defi} 
	Let $\G$ be an ultragraph. For each $A\in \G^0$, we define $\varepsilon(A):=\{e\in \G^1:s(e)\in A\}$. We say that $A\in \G^0$ is an infinite emitter if $|\varepsilon(A)|=\infty$. Otherwise we say that $A$ is a finite emitter. 
\end{defi}

\begin{defi}
\label{def_G^0}
	Let $\G$ be an ultragraph and $A\in \G^0$. We say that $A$ is an \emph{minimal infinite emitter} if $A$ is an infinite emitter, contains no proper subsets (in $\G^0$) that are infinite emitters, and contains no proper subsets (in $\G^0$) which are finite emitters and have
	infinite cardinality. Equivalently, $A$ is a minimal infinite emitter if it is an infinite emitter (i.e $|\varepsilon (A)|=\infty $) and $ \nexists B\subsetneq A\cap \G^0$ such that $|\varepsilon (B)|=\infty$ and $\nexists B \subsetneq A\cap\G^0$ such that $|B|=\infty$ and $|\varepsilon (B)|<\infty$. We denote by $A_\infty$ the set of all minimal infinite emitters in $\G^0$. Therefore,
	$$ A_\infty=\{ A\in\G^0:|\varepsilon(A)|=\infty; \nexists B\subsetneq A\cap \G^0 \text{ / } |\varepsilon(B)|=\infty \text{ or }|\varepsilon(B)|<\infty\text{ with } |B|=\infty \}.$$ 
\end{defi}
	
\begin{obs}
Notice that if the ultragraph has no sinks, the definition above coincides with the definition of minimal infinite emitters given in \cite{danidani}. 
\end{obs}

\begin{exemplo}
Below are illustrate two different kinds of minimal infinite emitters:
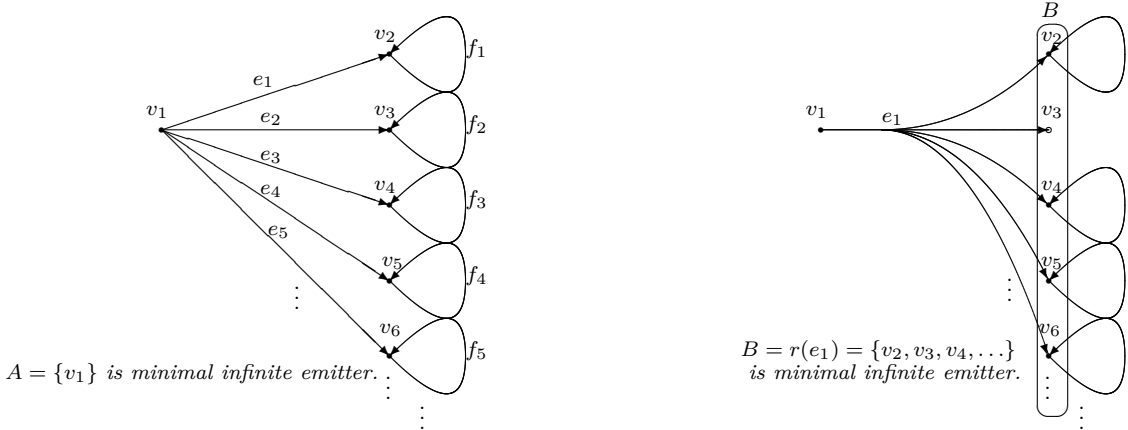
\begin{center}
\begin{multicols}{2}
\setlength{\unitlength}{1mm}
\begin{picture}(50,45)
		\put(0,32){\scriptsize$v_1$}
		\put(2,30){\circle*{0.7}}
		\put(30,42){\scriptsize$v_2$}
		\put(32,40){\circle*{0.7}}
		\put(30,32){\scriptsize$v_3$}
		\put(32,30){\circle*{0.7}}
		\put(30,22){\scriptsize$v_4$}
		\put(32,20){\circle*{0.7}}
		\put(31,12){\scriptsize$v_5$}
		\put(32,10){\circle*{0.7}}
		\put(30.7,3){\scriptsize$v_6$}
		\put(32,0){\circle*{0.7}}
		
		\put(2,30){\vector(3,1){30}}
		\put(2,30){\vector(1,0){30}}
		\put(2,30){\vector(3,-1){30}}
		\put(2,30){\vector(3,-2){30}}
		\put(2,30){\vector(1,-1){30}}
		
		\put(14,36){\scriptsize$e_1$}
		\put(15,31){\scriptsize$e_2$}
		\put(15,26){\scriptsize$e_3$}
		\put(15,21.5){\scriptsize$e_4$}
		\put(16,16){\scriptsize$e_5$}
		
		\put(42,40){\scriptsize$f_1$}
		\put(42,30){\scriptsize$f_2$}
		\put(42,20){\scriptsize$f_3$}
		\put(42,10){\scriptsize$f_4$}
		\put(42,0){\scriptsize$f_5$}
		
		\begin{rotate}{90}
		\put(6,-20){\scriptsize$\ldots$}
		\end{rotate}

		\qbezier(32,40)(42,50)(42,40)
		\qbezier(32,40)(42,30)(42,40)
		\qbezier(32,30)(42,40)(42,30)
		\qbezier(32,30)(42,20)(42,30)
		\qbezier(32,20)(42,30)(42,20)
		\qbezier(32,20)(42,10)(42,20)
		\qbezier(32,10)(42,20)(42,10)
		\qbezier(32,10)(42,0)(42,10)
		\qbezier(32,0)(42,10)(42,0)
		\qbezier(32,0)(42,-10)(42,0)
		
		\put(32,40){\vector(-1,-1){0}}
		\put(32,30){\vector(-1,-1){0}}
		\put(32,20){\vector(-1,-1){0}}
		\put(32,10){\vector(-1,-1){0}}
		\put(32,0){\vector(-1,-1){0}}
		
		\begin{rotate}{90}
		\put(-6,-32){\scriptsize$\ldots$}
		\end{rotate}
		
		\begin{rotate}{90}
		\put(-10,-35){\scriptsize$\ldots$}
		\end{rotate}

		\put(-20,-3){\scriptsize $A=\{v_1\}$ is minimal infinite emitter.}
\end{picture}
	
		\setlength{\unitlength}{1mm}
\begin{picture}(50,45)
		\put(0,32){\scriptsize$v_1$}
		\put(2,30){\circle*{0.7}}
		\put(31,41.8){\scriptsize$v_2$}
		\put(32,40){\circle*{0.7}}
		\put(31,32){\scriptsize$v_3$}
		\put(32,30){\circle{0.7}}
		\put(31,22){\scriptsize$v_4$}
		\put(32,20){\circle*{0.7}}
		\put(31,12){\scriptsize$v_5$}
		\put(32,10){\circle*{0.7}}
		\put(30.7,3){\scriptsize$v_6$}
		\put(32,0){\circle*{0.7}}
		\put(10,31){\scriptsize$e_1$}
		\put(31,45){\scriptsize$B$}

		\linethickness{0.2mm}
		\put(2,30){\line(1,0){10}}
		\thinlines
		\qbezier(10,30)(22,30)(32,40)
		\qbezier(10,30)(22,30)(32,30)
		\qbezier(10,30)(22,30)(32,20)
		\qbezier(10,30)(22,30)(32,10)
		\qbezier(10,30)(22,30)(32,0)
		\put(31.8,39.8){\vector(1,1){0}}
		\put(32,30){\vector(1,0){0}}
		\put(31.8,20.2){\vector(1,-1){0}}
		\put(31.8,10.2){\vector(1,-2){0}}
		\put(31.8,0.2){\vector(1,-3){0}}
		\begin{rotate}{90}
		\put(7,-27){\scriptsize$\ldots$}
		\end{rotate}
		
		\qbezier(32,40)(42,50)(42,40)
		\qbezier(32,40)(42,30)(42,40)
		\qbezier(32,20)(42,30)(42,20)
		\qbezier(32,20)(42,10)(42,20)
		\qbezier(32,10)(42,20)(42,10)
		\qbezier(32,10)(42,0)(42,10)
		\qbezier(32,0)(42,10)(42,0)
		\qbezier(32,0)(42,-10)(42,0)
		
		\put(32,40){\vector(-1,-1){0}}
		\put(32,20){\vector(-1,-1){0}}
		\put(32,10){\vector(-1,-1){0}}
		\put(32,0){\vector(-1,-1){0}}
		
		\begin{rotate}{90}
		\put(-6,-32){\scriptsize$\ldots$}
		\end{rotate}
		
		\begin{rotate}{90}
		\put(-10,-35){\scriptsize$\ldots$}
		\end{rotate}
	
		\put(31,18){\oval(4,52)}
		
		\put(-10,0){\scriptsize$B=r(e_1)=\{v_2,v_3,v_4,\ldots\}$}
		\put(-9,-3){\scriptsize is minimal infinite emitter.}
		
		\end{picture}
\end{multicols}
\end{center}
\end{exemplo}
	\vspace{20pt}
\begin{obs}
Note that in the second ultragraph of the example above, the range $r(e_1)$ contains a sink. In fact, if there was a finite number of sinks in $r(e_1)$ we would still have a minimal infinite emitter. Furthermore, notice that for $r(e_1)$ to be a minimal infinite emitter it can not contain a proper subset in $\G^0$ which has infinite cardinality and is a finite emitter. 
\end{obs}

The following proposition will be useful throughout our work.

\begin{prop}\label{obs3}
	If $A$ is a minimal infinite emitter and $B$ is an infinite emitter in $\G^0$, then either $A\subseteq B$ or the intersection $A\cap B$ is at most finite.
\end{prop}
\begin{proof}
Indeed, if $A\cap B=\emptyset$ or $|A\cap B|<\infty$, then their intersection is finite and the result follows. If $|A\cap B|=\infty$ then we have two cases: If $|\varepsilon(A\cap B)|=\infty$, we must have  $A\cap B=A$, otherwise we would have $A\cap B\subsetneq A$ and $|\varepsilon(A\cap B)|=\infty$, which contradicts the minimality of $A$. Hence $A\subseteq B$. On the other hand, if $|\varepsilon(A\cap B)|<\infty$, then again we must have  $A\cap B=A$, since otherwise $A\cap B$ is a proper subset of $A$ that is a finite emitter of infinite cardinality.
 \end{proof}	
	
 The following lemma characterizes minimal infinite emitters. In fact, it is a generalization of Lemma~3.3 in \cite{danidani}, and since its proof is analogous we omit it. 
	
\begin{lema}\cite[Lemma~3.3]{danidani}
\label{Lema_1}
	Let $x=(\alpha,A)\in\p$ be such that $A$ is a minimal infinite emitter. If the cardinality of $A$ is finite, then it is equal at to one. If the cardinality of $A$ is infinite, then $A=\displaystyle\bigcap_{e\in Y}r(e)$ for some finite set $Y\subseteq \G^1$.
\end{lema}

To define our topological space we need to extend the notion of minimal infinite emitters to generalized vertices that are finite emitters. For this, denote the set of sinks by $G^0_{sink}$, or alternatively, let $G^0_{sink}:=G^0_s:=\{v\in G^0:s^{-1}(v)=\emptyset\}\subseteq G^0$. Define $\G^0_{sink}:=\G^0_s:=\displaystyle\bigsqcup_{v_i\in G^0_s} \{\{v_i\}\}\subseteq \G^0$, i.e, $\G^0_s$ is the collection of all singletons of $\G^0$ whose element is a sink. 
		
	\begin{defi}
	Let $\G$ be an ultragraph and $A\in \G^0$. We say that $A$ is a \emph{minimal sink} if $|A|=\infty$, $|\varepsilon(A)|<\infty$ and $A$ has no subsets (in $\G^0$) with infinite cardinality. We call $A_s$ the set of all minimal sinks in $\G^0$. Therefore, $$A_s=\{A\in\G^0 : |A|=\infty;\varepsilon(A)<\infty; \nexists B\subsetneq A\cap\G^0\text{ such that }|B|=\infty\}.$$  
	\end{defi}
	
As with minimal infinite emitters, minimal sinks have some interesting properties.

\begin{prop}\label{prop1.16}
If $A_1$, $A_2\in A_s$ then either $|A_1\cap A_2|<\infty$ or $A_1=A_2$.
\end{prop}
\begin{proof}
Indeed, if $A_1\cap A_2=\emptyset$, the result follows. If $A_1\cap A_2\neq\emptyset$, then either $|A_1\cap A_2|<\infty$ (and the result follows) or $|A_1\cap A_2|=\infty$. In this case, if there exists a vertex $v\in A_1\backslash A_2$, then we have that $A_1\cap A_2\in \G^0$ is such that $A_1\cap A_2\subsetneq A_1$ and $|A_1\cap A_2|=\infty$, which contradicts the minimality of $A_1$. Therefore, $A_1=A_2$.
\end{proof}
	
The final step to construct the boundary path space is to consider the following sets:
	$$X_{min}:=\{ (\alpha,A)\in\p:|\alpha|\geq 1, A\in A_\infty\}\cup
	\{ (A,A)\in \p^0:A\in A_\infty\};$$
	$$ X_{sin}:=\{ (\alpha,A)\in\p:|\alpha|\geq 1, A\in \G^0_s\sqcup A_s\}\cup \{ (A,A)\in \p^0:A\in \G_s^0\sqcup A_s\}.$$
	
\begin{defi}
		\label{defi_esp_ult_fro}
		The \emph{boundary ultrapath space} $X$ is the set
		$$X:=\p^\infty \sqcup X_{min}\sqcup X_{sin}.$$
		We also define $$X_{fin}:=X_{min}\sqcup X_{sin}.$$ 
\end{defi}
\begin{obs} Clearly $X=\p^\infty\sqcup X_{fin}$. Also, for ultragraphs without sinks $ X_{sin} = \emptyset $, and therefore $ X_{fin} = X_{min} $. Hence in this case the definition above coincides with the definition of boundary path space in \cite{danidani}.
\end{obs}

\begin{obs}
Notice that if $(\alpha,A)\in X_{fin}$, then either $|A|=1$ or $|A|=\infty$. 
\end{obs}
	
\begin{exemplo} Consider the ultragraph below. Examples of elements in the boundary ultrapath space associated are the ultrapaths  $(\beta_1\beta_2\beta_3,\{v_1\})$, $(\beta_3,\{v_2\})$,  $\beta_1\beta_2\beta_3\beta_2\beta_3\ldots$, $ (\beta_2\beta_3,\{v_5\})$, $ \beta_3\gamma_1\gamma_2\gamma_3\ldots$, $(\{v_1\},\{v_1\})$, $(\{v_2\},\{v_2\})$, $(\{v_5\},\{v_5\})$,  $\beta_2\beta_3\beta_2\beta_3\gamma_1\gamma_2\gamma_3\ldots$,  and $(\delta_1\delta_2\delta_3,r(\delta_3))$. In fact any ultrapath whose range is equal to $\{v_1\}$, $\{v_5\}$, or $r(\delta_3)$ belong to $X_{min}$. If the ultrapath has range equal to $\{v_2\}$ or $r(\beta_5)$, it belongs to $X_{sin}$. Finally, a path that goes through $\gamma_3\ldots$ must be infinite and so it belongs to $\p^\infty$.
\end{exemplo}
\vspace*{10pt}
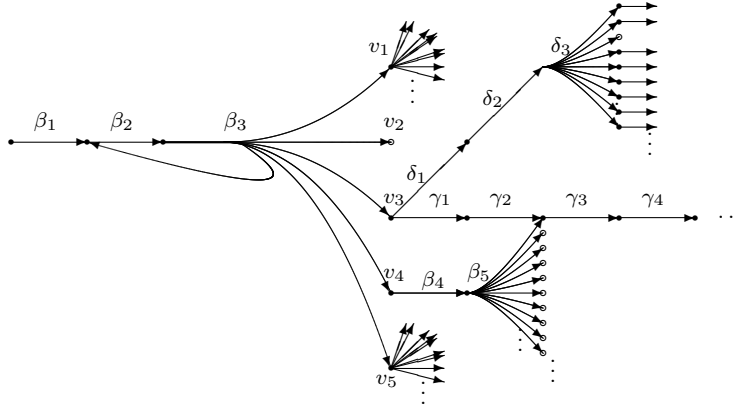
\begin{center}
\setlength{\unitlength}{1mm}
\begin{picture}(60,45)
		\put(2,30){\circle*{0.7}}
		\put(29,42){\scriptsize$v_1$}
		\put(32,40){\circle*{0.7}}
		\put(31,32){\scriptsize$v_2$}
		\put(32,30){\circle{0.8}}
		\put(31,22){\scriptsize$v_3$}
		\put(32,20){\circle*{0.7}}
		\put(31,12){\scriptsize$v_4$}
		\put(32,10){\circle*{0.7}}
		\put(30,-2){\scriptsize$v_5$}
		\put(32,0){\circle*{0.7}}
		\put(-8,30){\circle*{0.7}}
		\put(-18,30){\circle*{0.7}}

		\put(-15,32){\scriptsize$\beta_1$}
		\put(-5,32){\scriptsize$\beta_2$}
		\put(10,32){\scriptsize$\beta_3$}		
		\put(36,11){\scriptsize$\beta_4$}
		\put(42,12){\scriptsize$\beta_5$}
		
		\linethickness{0.3mm}
		\put(2,30){\line(1,0){10}}
		\thinlines
		\put(-8,30){\vector(1,0){10}}
		\put(-18,30){\vector(1,0){10}}
		\qbezier(10,30)(22,30)(32,40)
		\qbezier(10,30)(22,30)(32,30)
		\qbezier(10,30)(22,30)(32,20)
		\qbezier(10,30)(22,30)(32,10)
		\qbezier(10,30)(22,30)(32,0)
		\qbezier(12,30)(27,20)(-8,30)
		\qbezier(42,10)(45,10)(52,12)
		\qbezier(42,10)(45,10)(52,14)
		\qbezier(42,10)(45,10)(52,16)
		\qbezier(42,10)(45,10)(52,18)				
		\qbezier(42,10)(45,10)(52,20)		
		\qbezier(42,10)(45,10)(52,8)	
		\qbezier(42,10)(45,10)(52,6)	
		\qbezier(42,10)(45,10)(52,4)				
		\qbezier(42,10)(45,10)(52,2)		
		\qbezier(52,40)(55,40)(62,42)
		\qbezier(52,40)(55,40)(62,44)
		\qbezier(52,40)(55,40)(62,46)
		\qbezier(52,40)(55,40)(62,48)				
		\qbezier(52,40)(55,40)(62,40)		
		\qbezier(52,40)(55,40)(62,38)	
		\qbezier(52,40)(55,40)(62,36)	
		\qbezier(52,40)(55,40)(62,34)				
		\qbezier(52,40)(55,40)(62,32)		
		
		\put(62,42){\vector(4,1){0}}
		\put(62,44){\vector(2,1){0}}
		\put(62,46){\vector(4,3){0}}
		\put(62,48){\vector(1,1){0}}
		\put(62,40){\vector(1,0){0}}
		\put(62,38){\vector(4,-1){0}}
		\put(62,36){\vector(2,-1){0}}
		\put(62,34){\vector(4,-3){0}}
		\put(62,32){\vector(1,-1){0}}
		
		\put(62,42){\vector(1,0){5}}
		\put(62,46){\vector(1,0){5}}
		\put(62,48){\vector(1,0){5}}
		\put(62,40){\vector(1,0){5}}
		\put(62,38){\vector(1,0){5}}
		\put(62,36){\vector(1,0){5}}
		\put(62,34){\vector(1,0){5}}
		\put(62,32){\vector(1,0){5}}

		\put(31.8,39.8){\vector(1,1){0}}
		\put(32,30){\vector(1,0){0}}
		\put(31.8,20.2){\vector(1,-1){0}}
		\put(31.8,10.2){\vector(1,-2){0}}
		\put(31.8,0.2){\vector(1,-3){0}}
		\put(-7.8,30){\vector(-3,1){0}}
		
		\put(32,40){\vector(1,3){2}}
		\put(32,40){\vector(1,2){3}}
		\put(32,40){\vector(1,1){5}}
		\put(32,40){\vector(4,3){6}}
		\put(32,40){\vector(3,1){7}}
		\put(32,40){\vector(3,2){6}}
		\put(32,40){\vector(4,1){7}}
		\put(32,40){\vector(1,0){7}}
		\put(32,40){\vector(4,-1){7}}
		
		\put(32,20){\vector(1,0){10}}
		\put(42,20){\vector(1,0){10}}
		\put(52,20){\vector(1,0){10}}
		\put(62,20){\vector(1,0){10}}
		\put(75,20){\scriptsize$\ldots$}

		\put(42,10){\vector(1,0){10}}
		\put(52,12){\vector(4,1){0}}
		\put(52,14){\vector(2,1){0}}
		\put(52,16){\vector(4,3){0}}
		\put(52,18){\vector(1,1){0}}
		\put(52,20){\vector(2,3){0}}
		\put(52,8){\vector(4,-1){0}}
		\put(52,6){\vector(2,-1){0}}
		\put(52,4){\vector(4,-3){0}}
		\put(52,2){\vector(1,-1){0}}

		\put(37,22){\scriptsize$\gamma_1$}
		\put(45,22){\scriptsize$\gamma_2$}
		\put(55,22){\scriptsize$\gamma_3$}
		\put(65,22){\scriptsize$\gamma_4$}				
		
		\put(34,25){\scriptsize$\delta_1$}
		\put(44,35){\scriptsize$\delta_2$}	
		\put(53,42){\scriptsize$\delta_3$}
		\put(42,30){\circle*{0.7}}
		\put(62,40){\circle*{0.8}}
		\put(62,42){\circle*{0.8}}
		\put(62,44){\circle{0.8}}
		\put(62,46){\circle*{0.8}}
		\put(62,48){\circle*{0.8}}
		\put(62,38){\circle*{0.8}}
		\put(62,36){\circle*{0.8}}
		\put(62,34){\circle*{0.8}}
		\put(62,32){\circle*{0.8}}

		\put(42,30){\vector(1,1){10}}
		\put(32,20){\vector(1,1){10}}
		\put(32,10){\vector(1,0){10}}
		\put(42,20){\circle*{0.7}}
		\put(52,20){\circle*{0.7}}
		\put(62,20){\circle*{0.7}}
		\put(72,20){\circle*{0.7}}
		\put(42,10){\circle*{0.8}}
		\put(52,10){\circle{0.8}}
		\put(52,12){\circle{0.8}}
		\put(52,14){\circle{0.8}}
		\put(52,16){\circle{0.8}}
		\put(52,18){\circle{0.8}}
		\put(52,8){\circle{0.8}}
		\put(52,6){\circle{0.8}}
		\put(52,4){\circle{0.8}}
		\put(52,2){\circle{0.8}}
		\put(32,10){\vector(1,0){10}}
		
		\put(32,0){\vector(1,3){2}}
		\put(32,0){\vector(1,2){3}}
		\put(32,0){\vector(1,1){5}}	
		\put(32,0){\vector(4,3){6}}
		\put(32,0){\vector(3,1){7}}
		\put(32,0){\vector(3,2){6}}
		\put(32,0){\vector(4,1){7}}
		\put(32,0){\vector(1,0){7}}
		\put(32,0){\vector(4,-1){7}}
		
		\begin{rotate}{90}
		\put(35,-35){\scriptsize$\ldots$}
		\end{rotate}
		
		\begin{rotate}{90}
		\put(-5,-35){\scriptsize$\ldots$}
		\end{rotate}
		
		\begin{rotate}{90}
		\put(-2,-50.7){\scriptsize$\ldots$}
		\end{rotate}
		
		\begin{rotate}{90}
		\put(2,-45){\scriptsize$\ldots$}
		\end{rotate}
		
		\begin{rotate}{90}
		\put(28,-60.7){\scriptsize$\ldots$}
		\end{rotate}
		
		\begin{rotate}{90}
		\put(32,-55){\scriptsize$\ldots$}
		\end{rotate}
		
\end{picture}
\end{center}
	
\vspace{20pt}
	
To define a basis for a topology we need the following sets, called \emph{cylinders}. 
	
For each $(\beta,B)\in \p$, let 
$$D_{(\beta,B)}:=\left\{y\in X: y=\beta\gamma ; s(\gamma)\in B\right\},$$
and for each $(\beta,B)\in X_{fin}$, $F\subseteq \varepsilon(B)$ finite, and $S\subseteq B\cap G_{s}^0$ finite, define 
\begin{align*}
D_{(\beta,B),F,S}:=\left\{(\beta,B)\right\}\sqcup&\left\{ y\in X:y=\beta\gamma',|\gamma'|\geq 1, \gamma_1'\in\varepsilon (B)\setminus F \right\}\\ 
\sqcup &\left\{y\in X_{sin}:y=(\beta,\{v\}):v\in B\backslash S\right\}.
\end{align*}	
\begin{obs}\label{obs1.22}
	Notice that if $(\beta,B)\in \p^0\cap X_{fin}$ then $(\beta,B)=(B,B)$, and hence
	$$D_{(\beta,B),F,S}=D_{(B,B),F,S}=\left\{(B,B)\right\}\sqcup \left\{\gamma\in X: \gamma_1\in\varepsilon (B)\setminus F\right\}\sqcup \left\{(\{v\},\{v\})\in X_{sin}:v\in B\backslash S\right\}.$$
Thus, if $B=\{v\}$ is such that $v$ is a sink, then $D_{(B,B),F,S}=D_{(B,B)}=\big\{(\{v\},\{v\})\big\}$, i.e, the cylinder is a singleton. Similarly, $D_{(\beta,\{v\}),F,S}=\big\{(\beta,\{v\})\big\}$.

If either $F=\emptyset$ or $S=\emptyset$ we will write  $D_{(\beta,B),F}$ or $D_{(\beta,B),S}$ to refer to the cylinder $D_{(B,B),F,S}$ and this should not cause any confusion, as $F\subseteq \G^1$ and $S\subseteq G^0$.
\end{obs}

\begin{obs}
For ease of writing and reading, when referring to the finite subset $S\subseteq B\cap 
G_s^0$ that characterizes the cylinder sets we will only mention that $S\subset B$ and we will assume that $S$ is also a subset of $G_s^0$.  
\end{obs}
	
We can now show that a certain collection of cylinders form a basis for a topology on the boundary ultrapath space $X$.
	
\begin{prop}
	\label{prop_base_top}
	The collection of cylinders $$\big\{\{ D_{(\beta,B)}:(\beta,B)\in\p, |\beta|\geq 1\}\cup\{D_{(\beta,B),F,S}:(\beta,B)\in X_{fin},F\subseteq\varepsilon(B),S\subseteq B,|F|,|S|<\infty\}\big\}$$ is a countable basis for a topology in $X$. Furthermore, if $\gamma=e_1e_2\ldots\in \p^\infty$ then a neighborhood basis for $\gamma$ is given by $$\{D_{(e_1\ldots e_n,r(e_n))}:n\in\N\},$$
	and if $x=(\alpha,A)\in X_{fin}$ is such that $A$ is a minimal infinite emitter or a minimal sink, then a neighborhood basis for $x$ is given by $$\{D_{(\alpha,A),F,S}:F\subseteq \varepsilon(A),S\subseteq A,|F|,|S|<\infty\}.$$
	\end{prop}
\begin{proof}
	Clearly the collection of cylinders described is countable and $X$ can be written as a union of the cylinders in the collection. 
We have to show that if $x\in X$ is such that $x\in A\cap B$ with $A$ and $B$ cylinders, then there is a cylinder $C$ such that $x\in C\subseteq A\cap B$. 

First, if $\gamma\in\p^\infty$ is such that $\gamma\in A\cap B$ with $A, B$ cylinders, then it is clearly possible to find an initial segment of $\gamma $, say $\gamma_1$, with $|\gamma_1|$ larger than the ultrapath sizes that define $A$ and $B$. Then $D_{(\gamma_1,r(\gamma_1))}$ contains $\gamma$ and is contained in $A\cap B$.
	
Let $x=(\alpha,A)\in X_{min}$. Suppose that $x\in D_{(\beta_1,B_1),F_1,S_1}\cap D_{(\beta_2,B_2),F_2,S_2}$ with $|\beta_2|>|\beta_1|$ (notice that $|\beta_1|$ can be equal to zero). Then $(\beta_1,B_1)$ is an initial segment of $(\beta_2,B_2)$ and hence $x\in D_{(\beta_2,B_2),F_2,S_2}\subseteq D_{(\beta_1,B_1),F_1,S_1}$. If $|\beta_1|=|\beta_2|=0$, since $B_1$ and $B_2$ are minimal, we have $\beta_1=B_1$ and $\beta_2=B_2$. Hence, if $|x|>0$, just note that $x\in D_x\subseteq  D_{(\beta_1,B_1),F_1,S_1}\cap D_{(\beta_2,B_2),F_2,S_2}$. On the other hand if $|x|=0$ then $x=(A,A)$, so we cannot have $B_1$ nor $B_2$ as minimal sinks (or even as a singleton with only a sink), since minimal sinks do not emit infinite edges and $A$ emits. Thus, $B_1$ and $B_2$ must be minimal infinite emitters. Therefore, by Proposition~\ref{obs3}, we have that either $|B_1\cap B_2|<\infty$ or $B_1=B_2$. In both cases, since $A \subseteq B_1\cap B_2$ and $A$, $B_1$ and $B_2$ are minimal infinite emitters, we obtain that $x\in D_{(A,A),F_1\cup F_2,S_1\cup S_2}\subseteq D_{(\beta_1,B_1),F_1,S_1}\cap D_{(\beta_2,B_2),F_2,S_2}$.
	
Suppose that $x=(\alpha,A)\in D_{(\beta,B_1),F_1,S_1}\cap D_{(\beta,B_2),F_2,S_2}$, with $|\beta|\geq 1$. If $|x|>|\beta|\geq 1$ just note that $x\in D_{(\alpha, A)}\subseteq D_{(\beta,B_1),F_1,S_1}\cap D_{(\beta,B_2),F_2,S_2}$. On the other hand, if $|x|=|\beta|$ then $x=(\beta, A)$, where $A$ is an infinite emitter. Thus, $B_1$ and $B_2$ cannot be minimal sinks (nor singletons with only a sink). Hence $B_1$ and $B_2$ are minimal infinite emitters. By Proposition~\ref{obs3}, we have that either $|B_1\cap B_2|<\infty$ or $B_1=B_2$. As in the previous case, we get $x\in D_{(\alpha,A),F_1\cup F_2,S_1\cup S_2}\subseteq D_{(\beta,B_1),F_1,S_1}\cap D_{(\beta,B_2),F_2,S_2}$.

	If we suppose that $x\in D_{(\beta,B)}\cap D_{(\beta_1,B_1),F_1,S_1}$, then we have 3 cases: if $|\beta|>|\beta_1|$ then $D_{(\beta,B)}\subseteq D_{(\beta_1,B_1),F_1,S_1}$. If $|\beta_1|>|\beta|\geq 1$ then $D_{(\beta_1,B_1),F_1,S_1}\subseteq D_{(\beta,B)}$. But if $|\beta_1|=|\beta|\geq 1$ then we have $\beta_1=\beta$ and, by minimality of $B_1$, we have $|B\cap B_1|<\infty$ or $B_1\subseteq B$. If $|B\cap B_1|<\infty$, then $x\in D_{(\beta,B\cap B_1),F_1,S_1}\subseteq D_{(\beta,B)}\cap D_{(\beta_1,B_1),F_1}$, and if $B_1\subseteq B$ then $D_{(\beta_1,B_1),F_1,S_1}\subseteq D_{(\beta,B)}$.
		
	Finally, if $x\in D_{(\beta_1,B_1)}\cap D_{(\beta_2,B_2)}$ and $|\beta_1|>|\beta_2|$ then $D_{(\beta_1,B_1)}\subseteq D_{(\beta_2,B_2)}$. If $|\beta_1|=|\beta_2|>0$ then $\beta_1=\beta_2$ and $x\in D_{(\beta_1,B_1\cap B_2)}\subseteq D_{(\beta_1,B_1)}\cap D_{(\beta_2,B_2)}$. 
		
	For $x=(\alpha,A)\in X_{sin}$ with $A$ a minimal sink, the proof is handled similarly. Just change, eventually, $|A\cap B|<\infty$ by $|\varepsilon(A\cap B)|<\infty$.
	
	For $x=(\alpha,\{v\})\in X_{sin}$, with $v$ sink, just note that $(\alpha,\{v\})\in\{(\alpha,\{v\})\}=D_x\subseteq A\cap B$, for every cylinder $A$ and $B$.
	
	For the second part, let $\gamma\in \p^\infty$ so that $\gamma=e_1e_2\ldots$. Suppose that $\gamma\in D_{(\beta,B),F,S}$. Then we have $(\beta,B)=(e_1\ldots e_n,B)$ or $(\beta,B)=(B,B)$, where $B$ is a minimal infinite emitter or a minimal sink. If $(\beta,B)=(e_1\ldots e_n,B)$ then $e_{n+1}\notin F$, so $\gamma\in D_{(\beta e_{n+1},r(e_{n+1}))}\subseteq D_{(\beta,B),F,S}$. If $(\beta,B)=(B,B)$ then $e_1\notin F$, and hence $\gamma\in D_{( e_1,r(e_1))}\subseteq D_{(\beta,B),F}$. On the other hand, if $\gamma\in D_{(\beta,B)}$ then obviously $(\beta,B)=(e_1\ldots e_n,B)$, and hence $\gamma\in D_{(\beta e_{n+1},r(e_{n+1}))}\subseteq D_{(\beta,B)}$.
		
	For $x=(\alpha,A)\in X_{fin}$ with $A$ a minimal infinite emitter or a minimal sink, first suppose that $x\in D_{(\beta,B),F,S}$. Obviously $|\alpha|\geq|\beta|\geq 0$. If $|\alpha|=|\beta|= 0$ then $\alpha=A$ and $\beta=B$, furthermore, $A\subseteq B$. Since $A$ is minimal and $B$ is minimal, then $A=B$. Therefore, we can consider $F'=F$ and $S'=S$, so we have $x\in D_{(\alpha,A),F',S'}\subseteq D_{(\beta,B),F,S}$. If $|\alpha|>|\beta|= 0$ or $|\alpha|>|\beta|> 0$ then $\beta$ is an initial segment of $\alpha$. In both cases, we have $x\in D_{(\alpha,A)}\subseteq D_{(\beta,B),F,S}$. Finally, if $|\alpha|=|\beta|> 0$, we have $\alpha=\beta$ and $A\subseteq B$. Since both sets are minimal, $A=B$. Thus, $x\in D_{(\alpha,A),F,S}\subseteq D_{(\beta,B),F,S}$. On the other hand, if $x\in D_{(\beta,B)}$ then $|\alpha|\geq|\beta|\geq 1$. If $|\alpha|>|\beta|$ then $\beta$ is an initial segment of $\alpha$ and hence $x\in D_{(\alpha,A),\emptyset}\subseteq D_{(\beta,B)}$. If $|\alpha|=|\beta|$, then $x\in D_{(\alpha,A),\emptyset}\subseteq D_{(\beta,B)}$.	 
\end{proof}
	
\begin{obs}
	Notice that if $\G$ is an ultragraph that has no sinks, our proposition coincides with Proposition 3.4 in \cite{danidani} and the topological space here coincides with the topological space defined there. Also, if $\G$ is a graph, the only minimal infinite emitters are sets consisting of only one singular vertex (in graph terminology), and our topological space coincides with the boundary path space of Definition 2.1 in \cite{webster}.
\end{obs}

The topology described in Proposition~\ref{prop_base_top} has some interesting properties that we describe below.

\begin{prop}
\label{cilindros_fechados} 
    Each cylinder that forms the basis of the topology of $X$ described in Proposition~\ref{prop_base_top} is closed.
\end{prop}
\begin{proof}
	First we show for a cylinder $D_{(\beta,B)}$ with $|\beta|\geq 1$. Suppose that $\beta=\beta_1\ldots\beta_n$. Let $\delta\in D_{(\beta,B)}^C$. We have 4 cases:
\begin{itemize}
    \item If $|\delta|=\infty$, say $\delta=\delta_1\delta_2\ldots$, then $s(\delta_{n+1})\notin B$ or $\delta_1\ldots\delta_n\neq \beta$. In both, $\delta\in D_{(\delta_1\ldots\delta_{n+1},r(\delta_{n+1}))}\subseteq D_{(\beta,B)}^C$.
    \item If $|\delta|<\infty$ and $\delta=(\delta,\{v\})\in X_{sin}$, we have $\delta\in \{\delta\}=D_\delta\subseteq D_{(\beta,B)}^C$.
    \item If $|\delta|<\infty$ and $\delta=(\delta,A)\in X_{min}$, we have 3 sub cases:
\begin{enumerate}[$\circ$]
	\item $|\delta|<|\beta|$: Of course we only have to deal with the case where $\delta$ is the initial segment of $\beta$, otherwise the result is obvious. Note that $\delta_1\ldots\delta_k=\beta_1\ldots\beta_k$, where $|\delta|=k<n$. If $s(\beta_{k+1})\in A$ then $\delta\in D_{(\delta,A),\{\beta_{k+1}\}}\subseteq D_{(\beta,B)}^C$. If $s(\beta_{k+1})\notin A$ then $\delta\in D_{(\delta,A)}\subseteq D_{(\beta,B)}^C$.
	
	\item $|\delta|=|\beta|$: We only deal with the case where $\beta=\delta$ because otherwise clearly $\delta\in D_\delta\subseteq D_{(\beta,B)}^C$. So suppose that $\beta=\delta$. We have 2 cases: If $|A|<\infty$ then by Lemma~\ref{Lema_1}, $|A|=1$ and therefore $A=\{u\}$, where $u$ is a minimal infinite emitter. So $u\notin B$, since $\delta\in D_{(\beta,B)}^C$. Thus $\delta \in D_{\delta}\subseteq D_{(\beta,B)}^C$. If $|A|=\infty$ then $|A\cap B|<\infty$, otherwise, by the minimality of $A$ we would have $A\cap B=A$, i.e $A\subseteq B$, which would make $\delta\in D_{(\beta,B)}$, a contradiction. Also note that $A\cap B$ contains no infinite emitter as it would contradict the minimality of $A$. So consider $F=\{e\in \G^1: s(e)\in A\cap B \}$. Also, consider $S=\{v\in G^0_s: v\in A\cap B\}$. Notice that  the sets $F$ and $S$ are finite and $\delta\in D_{(\delta,A),F,S}\subseteq D_{(\beta,B)}^C$.
	
	\item $|\delta|>|\beta|$: Again, we only deal with the case where $\delta_1\ldots\delta_n=\beta_1\ldots\beta_n$. But note that for $\delta\in D_{(\beta,B)}^C$ we must have $s(\delta_{n+1})\notin B$, and so $\delta\in D_\delta \subseteq D_{(\beta,B)}^C$.
\end{enumerate} 
\item If $|\delta|<\infty$ and $\delta=(\delta,A)\in X_{sin}$, with $A$ a minimal sink, we have 3 sub cases:
\begin{enumerate}[$\circ$]
	\item $|\delta|<|\beta|$: We only deal with the case where $\delta$ is the initial segment of $\beta$. Note that $\delta_1\ldots\delta_k=\beta_1\ldots\beta_k$, where $|\delta|=k<n$. If $s(\beta_{k+1})\in A$ then $\delta\in D_{(\delta,A),\{\beta_{k+1}\}}\subseteq D_{(\beta,B)}^C$. If $s(\beta_{k+1})\notin A$ then $\delta\in D_{(\delta,A)}\subseteq D_{(\beta,B)}^C$ 
	
	\item $|\delta|=|\beta|$: we only deal with the case where $\beta=\delta$. 
Since $|A|=\infty$, we have by the minimality of $A$ that $|A\cap B|<\infty$. 
Furthermore, since $A$ is a minimal sink, $|\varepsilon(A\cap B)|<\infty$. Thus, we have $\delta\in D_{(\delta,A),\varepsilon(A\cap B),A\cap B}\subseteq D_{(\beta,B)^C}$.
	
	\item $|\delta|>|\beta|$: Again, we only deal with the case where $\delta_1\ldots\delta_n=\beta_1\ldots\beta_n$. But note that for $\delta\in D_{(\beta,B)}^C$, we must have $s(\delta_{n+1})\notin B$, and so $\delta\in D_\delta \subseteq D_{(\beta,B)}^C$.
\end{enumerate} 
\end{itemize}

	Now we show that each $D_{(B,B)}$, with $B$ a minimal infinite emitter, is closed. Let $\delta\in D_{(B,B)}^C$. So if $|\delta|>0$, we have $\delta \in D_\delta\subseteq D_{(B,B)}^C$. On the other hand, suppose that $|\delta|= 0$. Then $\delta=(A,A)\in X_{min}\sqcup X_{sin}$, with $A$ minimal or $\delta=(\{v\},\{v\})\in X_{sin}$. If $\delta=(\{v\},\{v\})$ then clearly $v\notin B$ and thus $\delta \in D_\delta\subseteq D_{(B,B)}^C$. On the other hand, if $\delta=(A,A)$ with $A$ minimal (a sink or an infinite emitter), then $0\leq|A\cap B|<\infty$, (otherwise, since $A$ and $B$ are minimal, we would have $A=B$, a contradiction). Furthermore, $|\varepsilon(A\cap B)|<\infty$, otherwise we get a contradiction with the minimality of $A$ and $B$. Thus, if $F=\{e\in \G^1:s(e)\in A\cap B\}$ and $S=A\cap B$, we have $\delta\in D_{(\delta,A)F,S} \subseteq D_{(B,B)}^C$. 

	For $D_{(\{v\},\{v\})}$ with $v$ sink, if necessary, just consider $S=\{v\}$ and observe that $\delta\in D_{\delta,S}\subseteq D_{(\{v\},\{v\})}^C$. The case $D_{(B,B)}$, with $B$ a minimal sink is handled similarly to $B$ a minimal infinite emitter. 
	
	Finally, for $D_{(\beta,B),F,S}$ with $B$ a minimal infinite emitter or minimal sink, note that $$D_{(\beta,B),F,S}=D_{(\beta,B)}\displaystyle\bigcap_{\gamma \in F} D_{(\beta\gamma,r(\gamma))}^C\displaystyle\bigcap_{v \in S} D_{(\beta,\{v\})}^C.$$
\end{proof}

Using techniques similar to the one presented above one can show that the boundary ultrapath space is Hausdorff. We state this below and leave the proof to the reader (who can also check more details in \cite{tesefelipe}).
	
\begin{prop} 
The boundary ultrapath space $X$, with the topology of Proposition~\ref{prop_base_top}, is Hausdorff.
\end{prop}

From the above, using Urysohn's Metrization Theorem, we conclude that the boundary ultrapath space is metrizable. Therefore, it is interesting to describe a notion of convergence of sequences in $X$. By Proposition~\ref{prop_base_top} we have the following:

\begin{corolario}
\label{corolario_convergencia}
Let $(x^n)_{n=1}^{\infty}\subseteq X$ be a sequence and $x\in X$. Clearly we have $x^n=(\gamma_1^n\ldots \gamma_{k_n}^n,A_n)$ with $A_n\in A_\infty\sqcup A_s$ or $x^n=(\gamma_1^n\ldots\gamma_{k_n}^n,\{v_n\})$, with $v_n$ sink, or $x^n=\gamma_1^n\gamma_2^n\ldots$. Thus, we have:
\begin{enumerate}[i.]
	\item If $|x|=\infty$, say $x=\gamma_1\gamma_2\ldots$, then $(x^n)_{n=1}^{\infty}$ converges to $x$ if, and only if, for each $M\in \N$ there exists $N\in \N$ such that $n>N$ implies that $|x^n|\geq M$ and $\gamma_i^n=\gamma_i$ for all $i=1,\dots,M$.
	\item If $|x|<\infty$, say $x=(\gamma_1\ldots\gamma_{k},A)$ with $A$ a minimal infinite emitter, then $(x^n)_{n=1}^{\infty}$ converges to $x$ if, and only if, for each finite subset $F\subseteq \varepsilon(A)$ and each finite subset $S\subseteq A$, there exists $N\in\N$ such that $n>N$ implies that $x^n=x$ or $|x^n|>k$, $\gamma_{k+1}^n\in\varepsilon(A)\backslash F$, and $\gamma_i^n=\gamma_i$ for all $i=1,\dots,k$;
	\item If $|x|<\infty$, say $x=(\gamma_1\ldots\gamma_{k},A)$ with $A$ a minimal sink, then $(x^n)_{n=1}^{\infty}$ converges to $x$ if, and only if, for each finite subset $S\subseteq A\cap G_s^0$, there exists $N\in\N$ such that $n>N$ implies that $x^n=x$ or $|x^n|=k$, $\gamma_i^n=\gamma_i$ for all $i=1,\dots,k$, and $v_n\in A\backslash S$ (i.e, in this last case, $x^n=(\gamma_1,\ldots,\gamma_k,\{v_n\})$, with $v_n$ sink and $v_n\in A\backslash S$).
	\item If $|x|<\infty$, say $x=(\gamma_1\ldots\gamma_{k},\{v\})$ with $v$ sink, then $(x^n)_{n=1}^{\infty}$ converges to $x$ if, and only if, there exists $N\in \N$ for which $n\geq N$ implies $x^n=x$.
\end{enumerate}
\end{corolario}
	
\begin{obs}
\label{Obs25}
	The cylinders $D_{(\beta,B)}$ with $(\beta,B)\in \p$ are not, in general, compact. Indeed, just note that if we consider an ultragraph $\G$ with an edge $e$ such that its range contain an infinite number of minimal infinite emitters then $D_{(e,r(e))}$ is not compact. 
\end{obs}

In order to obtain compact basic sets, we will impose a condition on the ultragraphs we are working with. The extra hypothesis required will be called (RFUM2). Such a condition is based on the Condition~(RFUM) required in \cite{danidani}. Note that if the ultragraph does not contain sinks then (RFUM) and (RFUM2) below coincide.
	
{\textbf{Condition (RFUM2):}}\label{RFUM2} We say that an ultragraph $\G$ satisfies Condition~(RFUM2) if for each edge $e\in \G^1$ its range can be written as $$r(e)=\bigcup_{n=1}^k A_n,$$
where $A_n$ is either a minimal infinite emitter, or a minimal sink, or a singleton formed by a sink or a regular vertex.

\begin{prop}
\label{cilindros_compactos}
	Let $\G$ be an ultragraph that satisfies Condition~(RFUM2). Then each basis element of the topology on $X$, as in Proposition \ref{prop_base_top}, is compact.  
\end{prop}
\begin{proof}
	First we show that $D_{(\beta,B)}$, with $|\beta|\geq 1$ and $(\beta,B)\in\p$, is sequentially compact. Let $(x^n)$ be a sequence in $D_{(\beta,B)}$. Then  $x^n=(\beta\alpha_1^n\ldots\alpha_{k_n}^n,A_n)\in X_{fin}$ with $A_n$ minimal, or  $x^n=(\beta\alpha_1^n\ldots\alpha_{k_n}^n,\{v_n\})\in X_{sin}$, or $x^n=\beta\alpha_1^n\alpha_2^n\ldots\in\p^\infty$. We need to find a convergent subsequence. We have two cases:
\begin{itemize}
\item If there exists a subsequence of  $x^n$ such that $|x^n|=|\beta|$: Since $\G$ satisfies Condition~(RFUM2) and $B(\subseteq r(\beta))$, we have that $B$ contains only a finite number of $A_n$'s (whether they are minimal infinite emitters or minimal sinks) and vertices $v_n$ that are regular or a sink and do not belong to any $A_n$. So if we have infinite terms of the form $x^n=(\beta,A_n)$, then we get a subsequence that converges to some $(\beta,A_{n_0})$, for some $n_0$. Now, if we have infinite terms of the form $(\beta,\{v_n\})$, where $v_n$ is not an element of any minimal sink, since we have a finite number of $v_n$'s, there is a constant subsequence of $x^n$ converging to some $(\beta,\{v_{n_0}\})$, for some $n_0$. If the $v_n$'s are elements of minimal sinks $A_n$'s, since we have a finite number of $A_n$'s, there exists a minimal sink $A\subseteq B$ such that $|\{n:v_n\in A\}|=\infty$. Thus, there is a subsequence of $x^n$ that converges to $(\beta,A)$ (when $x^n$'s are distinct) or to $(\beta,\{v_{n_0}\})$ for some $v_{n_0}$ (when $x^n$'s repeat an infinite number of times).

\item On the other case, there is a subsequence of $x^n$ such that $|x^n|>|\beta|$:

If there is no infinite number of indices such that $\alpha_1^n$ coincide, then we can assume, without loss of generality, that $\alpha_1^n\neq \alpha_1^m$ for all $n\neq m$, passing maybe, to a subsequence. So by the description of $\G^0$, $B$ must contain at least one infinite emitter, and since  $B\subseteq r(\beta)$ and $\G$ satisfy Condition~(RFUM2), $B$ can be written as finite union of minimal infinite emitters, minimal sinks, and singletons formed by sinks or regular vertices. So, since the number of regular vertices in the description of $B$ is finite, there is an $A\in A_\infty$, $A\subseteq B$, such that $|\{n:s(\alpha_1^n)\in  A\}|=\infty$, and hence we obtain a subsequence converging to $(\beta,A)$.

If there exists an infinite number of indexes for which $\alpha_1^n$ coincide, say $\alpha_1^n=\gamma_1$, we pass to a subsequence such that $x^n=(\beta\gamma_1\alpha_2^n\ldots\alpha_{k_n}^n,A_n)$, or $x^n=(\beta\gamma_1\alpha_2^n\ldots\alpha_{k_n}^n,\{v_n\})$, or $x^n=\beta\gamma_1\alpha_2^n\alpha_3^n\ldots$. Notice that $x^n\in D_{(\beta\gamma_1,r(\gamma_1))}$. Now, we repeat the procedure performed in the previous paragraph (where we split the proof in two) and hence we either obtain a subsequence converging to a finite sequence or we pass to a subsequence such that $x^n=(\beta\gamma_1\gamma_2\alpha_3^n\ldots\alpha_{k_n}^n,A_n)$, or $x^n=(\beta\gamma_1\gamma_2\alpha_3^n\ldots\alpha_{k_n}^n,\{v_n\})$, or $x^n=\beta\gamma_1\gamma_2\alpha_3^n\alpha_4^n\ldots$. Proceeding inductively, we either obtain at step $k$, a subsequence converging to a finite sequence or, through a Cantor diagonal argument, we obtain a subsequence converging to an infinite sequence $\beta\gamma_1\gamma_2\ldots$. 
\end{itemize}
		
Thus, we conclude that $D_{(\beta,B)}$ is sequentially compact. The proof that the cylinders of the form $D_{(\beta,B),F,S}$ are sequentially compact is completely analogous.
\end{proof}  
	
\begin{corolario}
	Let $X$ be the boundary ultrapath space associated to an ultragraph that satisfies Condition~(RFUM2). Then $X$ is a locally compact Hausdorff space. 
\end{corolario}
	
To end this section, we describe  a dense subset of $X$. 

\begin{prop}
    The set $\p^\infty\cup X_{\G^0_s}$, where $X_{\G^0_s}:=\{ (\alpha,A)\in\p:|\alpha|\geq 1, A\in \G^0_s\}\cup \{ (A,A)\in \p^0:A\in \G_s^0\}$ (the set of ultrapaths whose range is a (vertex) sink) is dense in $X$.
\end{prop}
\begin{proof}
	Indeed, if $x=(\alpha,A)\in X_{min}$ then $|\varepsilon(A)|=\infty$ and, for each $e_n\in\varepsilon(A)$, either there is an infinite path $\gamma^n=\gamma_1^n\gamma_2^n\ldots$ with $e_n=\gamma_1^n$, or there is a finite ultrapath $\gamma^n=(\gamma_1^n\ldots\gamma_{k_n}^n,\{v_n\})$ such that $v_n$ is a sink and $\gamma_1^n=e_n$. Thus, $x$ is the limit of the sequence $\{\alpha.\gamma^n\}$. On the other hand, if $x=(\alpha,A)\in X_{sin}$, with $A$ a minimal sink, then $|A|=\infty$, say $A=\{v_1,v_2,v_3,\ldots\}$. Then, $\{(\alpha,\{v_n\})\}_{n=1}^\infty$ converges to $x$. 
\end{proof}


\section{Ultragraph C*-algebra as partial crossed product}\label{ultracross}

In this section we define a partial action (of the free group on the edges) on the boundary path space, and use it to realize C*-algebras associated to ultragraphs that satisfy Condition~(RFUM2) as partial crossed products. This partial action will be important in the understanding of the ultragraph groupoid, see Section~\ref{gpd}, and the partial crossed product realization is key in the study of KMS states done in Section~\ref{KMS}.
We notice that the results we present completely generalize the description of graph C*-algebras as partial crossed products given in \cite[Theorem 3.1]{CarLar} (previous generalizations, see  \cite[Theorem 4.12]{danidani},  \cite[Theorem 4.11]{Dani&Dani2017}, excluded sinks) and also generalize the constructions in \cite[Theorem 4.12]{danidani} and \cite[Theorem 4.11]{Dani&Dani2017}. The results we present also exclude the use of filters and ultra filters in the topology (differently from what is done in \cite{GilDanie}). 	

\begin{obs}
	Unless otherwise noted, from now on all ultragraphs satisfy Condition~(RFUM2).
\end{obs}

Before we define the relevant partial action we state a auxiliary result, whose proof is similar to the proof of \cite[Lemma 4.1]{danidani}.

\begin{lema}
		\label{vertice_generalizado_clopen_compacto}
		The set $X_v:=\{(\alpha,A)\in X_{fin}:s(\alpha)=v\}\cup\{\gamma\in\p^\infty:s(\gamma)=v\}=\{x\in X:s(x)=v\}$ is nonempty, clopen, and compact, for each $v\in G^0$.
\end{lema}

Now we define the sets that will be used to construct the partial action.
	
\begin{defi}
		\label{def_X_e}
		Let $\F$ be the free group generated by $\G^1$. Let $P\subseteq \F$ be defined by: $$P:=\{e_1\ldots e_n\in\F:e_i\in\G^1:n\geq 1\};$$
		and for each element of $\F$ define subsets as follows.
		\begin{itemize}
			\item $X_0=X$, where 0 is the neutral element of group $\F$.
			\item If $a\in P$ then 
			$X_a=\{(\beta,B)\in X_{fin}:\beta_1\ldots\beta_{|a|}=a\}\cup\{\gamma\in\p^\infty:\gamma_1\ldots\gamma_{|a|}=a\}$;
			\\
			$X_{a^{-1}}=\{(A,A)\in X_{fin}:A\subseteq r(a)\}\cup\{(\beta,B)\in X_{fin}:s(\beta)\in r(a)\}\cup\{\gamma\in\p^\infty:s(\gamma)\in r(a)\}$;
			\item If $a,b\in P$, with  $ab^{-1}$ in its reduced form, then\\
			$X_{ab^{-1}}=\{(a,A)\in X_{fin}:A\subseteq r(a)\cap r(b)\}\cup\{(\beta,B)\in X_{fin}:\beta_1\ldots\beta_{|a|}=a$ and $s(\beta_{|a|+1})\in r(a)\cap r(b)\}$ $\cup\{\gamma\in\p^\infty:\gamma_1\ldots\gamma_{|a|}=a$ and $s(\gamma_{|a|+1})\in r(a)\cap r(b)\}$;
			\item For all other $c\in\F$, define $X_c=\emptyset$.
		\end{itemize}
\end{defi}
	
\begin{obs}
		Notice that if $a\in P$ is not a path in $X$, then $X_a$ is empty. Analogously, if $a,b\in P$ are such that $r(a)\cap r(b)=\emptyset$, then $X_{{ab}^{-1}}=\emptyset$.
\end{obs}

An important property of the above subsets is that they are compact and clopen, as we show below.

	\begin{prop}
		\label{X_c_clopen_compacto}
		The subsets $X_c$, with $c\in\F$, are clopen and compact in $X$.
	\end{prop}
	\begin{proof}
		For $a\in P$ (with $a$ a path), note that $X_a=D_{(a,r(a))}$, hence is open. Furthermore, by Propositions~\ref{cilindros_fechados} and \ref{cilindros_compactos} $X_a$ is closed and compact. For $X_{a^{-1}}$, by Condition~(RFUM2), we have $r(a)=\displaystyle\bigcup_{n=1}^{k} A_n$, where each $A_n$ is a minimal infinite emitter, or a minimal sink, or a regular vertex, or a sink. If $A_n$ is a minimal infinite emitter or a minimal sink, then $D_{(A_n,A_n)}$ (open) is closed and compact by Propositions~\ref{cilindros_fechados} and \ref{cilindros_compactos}. If $A_n$ is a regular vertex or a sink then, by Lemma~\ref{vertice_generalizado_clopen_compacto}, $D_{(A_n,A_n)}$ is clopen and compact. Hence, $X_{a^{-1}} =\displaystyle\bigcup_{n=1}^{k} D_{(A_n,A_n)}$ is clopen and compact.
		
		For $X_{ab^{-1}}$ with $a,b\in P$, where $a$ and $b$ are paths such that $r(a)\cap r(b)\neq \emptyset$. First, we show that $X_{ab^{-1}}$ is open, i.e, if $\beta\in X_{ab^{-1}}$, then there is an open $U$ such that $\beta\in U\subseteq X_{ab^{-1}}$. Let $\beta \in X_{ab^{-1}}$. 
		\begin{itemize}
			\item If $\beta\in \p^\infty$ then $\beta=a\beta'$ with $s(\beta')\in r(a)\cap r(b)$ and $\beta\in D_{(a\beta_1',r(\beta_1'))}\subseteq X_{ab^{-1}}$; 
			\item If $\beta=(\alpha,A)\in X_{fin}$ with $A$ minimal infinite emitter or minimal sink, we have 2 cases:
			\begin{enumerate}[$\circ$]
				\item If $|\alpha|>|a|$ then $\alpha=a\alpha'$ with $s(\alpha')\in r(a)\cap r(b)$, and hence $\beta\in D_{(a\alpha',r(\alpha'))}\subseteq X_{ab^{-1}}$;
				\item If $|\alpha|=|a|$ then $\alpha=a$ and hence $\beta=(a,A)$, with $A\subseteq r(a)\cap r(b)$, so that $\beta\in D_{(a,A)}\subseteq X_{ab^{-1}}$.
			\end{enumerate}
			\item If $\beta=(\alpha,\{v\})\in X_{sin}$, with $v$ sink, we have again 2 cases:
			\begin{enumerate}[$\circ$]
				\item If $|\alpha|>|a|$ then $\alpha=a\alpha'$ with $s(\alpha')\in r(a)\cap r(b)$, and hence $\beta\in D_{(a\alpha',r(\alpha'))}\subseteq X_{ab^{-1}}$;
				\item If $|\alpha|=|a|$ then $\alpha=a$ and hence $\beta=(a,\{v\})$, with $v\in r(a)\cap r(b)$, and so $\beta\in D_{(a,\{v\})}\subseteq X_{ab^{-1}}$.
			\end{enumerate}
		\end{itemize}
		Next we show that $X_{ab^{-1}}$ is closed and compact. Let $(\beta^n)_{n\in\N}\subseteq X_{ab^{-1}}$ be a sequence that converges to $\beta\in X$. Since $X_{ab^{-1}}\subseteq X_a$ and $X_a$ is closed, we have $\beta \in X_a$. Thus, we have 2 cases:
		\begin{itemize}
			\item If $|\beta|>|a|$, let $n_0\in \N$ be such that $\beta^n\in D_{(\beta_1\ldots\beta_{|a|+1},r(\beta_{|a|+1}))}$, for all $n\geq n_0$. Then $a=\beta_1^{n_0}\ldots\beta_{|a|}^{n_0}=\beta_1\ldots\beta_{|a|}$ and $\beta_{|a|+1}^{n_0}=\beta_{|a|+1}$. Since $s(\beta_{|a|+1}^{n_0})\in r(a)\cap r(b)$, we have $\beta\in X_{ab^{-1}}$;  
			\item If $|\beta|=|a|$ then $\beta= (a,A)$, with $A\in r(a)$ a minimal infinite emitter, or a minimal sink, or a sink. 
			\begin{enumerate}[$\circ$]
				\item If $A$ is sink then there exists $N\in\N$ such that $\beta^n=\beta=(a,A)$, $\forall n\geq N$. Hence $A=r(\beta)=r(\beta^n)\in r(a)\cap r(b)$, $\forall n\geq N$, because $(\beta^n)_{n\in\N}\subseteq X_{ab^{-1}}$. Therefore $\beta\in X_{ab^{-1}}$. 
				\item If $A$ is a minimal infinite emitter, we need to show that $A\subseteq r(b)$. Since $\beta^n\rightarrow \beta$, there exists $n_0$ such that if $n\geq n_0$ then $\beta^n=\beta$, or $|\beta^n|>|a|$ with $\beta_1\ldots \beta_{|a|}=a$ and $\beta_{|a|+1}^n\in \varepsilon (A)$. If $\beta^n=\beta$ for some $n$ the result follows. Hence suppose, without loss of generality, that $|\beta^n|>|a|$, $\forall n\geq n_0$. Suppose that $F:=\{\beta_{|a|+1}^n:n\geq n_0\}$ is a finite set. Thus, $\beta^n\notin D_{(a,A),F}$, which is a contradiction by Corollary \ref{corolario_convergencia}, since $\beta^n\rightarrow (a,A)\in X_{min}$. Hence the set $F$ is infinite. Since $\{s(\beta_{|a|+1}^n):n\geq n_0\}\subseteq r(a)\cap r(b)$ and $\{s(\beta_{|a|+1}^n):n\geq n_0\}\subseteq A$, we have $\{s(\beta_{|a|+1}^n):n\geq n_0\}\subseteq A\cap r(b)$. So, $A\cap r(b)$ is an infinite emitter. Since $A$ is minimal, we have $A=A\cap r(b)$, i.e, $A\subseteq r(b)$ and hence $\beta=(a,A)\in X_{ab^{-1}}$.
				\item If $A$ is a minimal sink, since $\beta^n\rightarrow \beta$, there exists $n_0$ such that if $n\geq n_0$ then $\beta^n=\beta$, or $|\beta^n|=|a|$ with $\beta_1\ldots \beta_{|a|}=a$ and $\beta^n=(a,\{v_n\})$ where $v_n\in A$. If $\beta^n=\beta$ for some $n$ the result follows. Hence suppose, without loss of generality, that $\forall n\geq n_0$, $|\beta^n|=|a|$, $\beta_1\ldots \beta_{|a|}=a$, and $\beta^n=(a,\{v_n\})$ where $v_n\in A$. Suppose that $H:=\{v_n:n\geq n_0\}$ is a set with finite cardinality. Hence, $\beta^n\notin D_{(a,A),H}$, which is a contradiction by Corollary \ref{corolario_convergencia}, since $\beta^n\rightarrow (a,A)\in X_{sin}$. Thus, the set $H$ has infinite cardinality. Clearly $H\subseteq A$ and, since $H=\{v_n:n\geq n_0\}\subseteq r(a)\cap r(b)$, we have that $H\subseteq A\cap r(b)$. Thus, $A\cap r(b)$ is a set which has infinite cardinality and such that $|\varepsilon(A\cap r(b))|<\infty$, since $|\varepsilon(A)|<\infty$. Since $A$ is minimal sink, we have $A=A\cap r(b)$, i.e, $A\subseteq r(b)$ and so $\beta=(a,A)\in X_{ab^{-1}}$.
			\end{enumerate}
		\end{itemize}
			We conclude that $X_{ab^{-1}}$ is closed. Since $X_a$ is compact and $X_{ab^{-1}}\subseteq X_a$ is closed, we obtain that $X_{ab^{-1}}$ is also compact.
	\end{proof} 
	
	When identifying $C^*(\G)$ with a partial crossed product we will associate to the projections $P_A$, $A\in \G^0$, the characteristic function of the sets defined below.
	
	\begin{defi}
		\label{def_X_A}
		For each $A\in \G^0$, let $X_A\subseteq X$ be defined by $X_A=\{x\in X: s(x)\subseteq A\}$.
	\end{defi}
	
	The following lemma shows that the sets defined above have the the correct properties to be associated with projections. 
	
	\begin{lema}
		\label{lema_4.8}
		For each $A,B\in \G^0$, we have $X_{A\cap B}=X_A\cap X_B$ and $X_{A\cup B}=X_A\cup X_B$.
	\end{lema}
	\begin{proof}
		Let $A,B\in\G^0$, then $$X_{A\cap B}=\{x\in X: s(x)\subseteq A\cap B\}=\{x\in X: s(x)\subseteq A\}\cap \{x\in X: s(x)\subseteq B\}=X_A\cap X_B.$$
		On the other hand, it is clear that $X_A\subseteq X_{A\cup B}$ and $X_B\subseteq X_{A\cup B}$, hence $X_A\cup X_B\subseteq X_{A\cup B}$. Suppose that $x\in X_{A\cup B}$. Then $s(x)\subseteq A\cup B$. If $|x|\geq 1$ then $x=(\gamma,C)$ and $s(\gamma)=s(x)$, where $s(x)$ is a a vertex in $G^0$, so $s(x)\subseteq A$ or $s(x)\subseteq B$, i.e $x\in X_A\cup X_B$. If $|x|=0$ and $s(x)$ is a sink (not minimal), then $s(x)$ is again a vertex in $G^0$, and therefore $x\in X_A\cup X_B$. If $|x|=0$ and $s(x)$ is not a sink, then $x=(C,C)$, with $C\in \G^0$ a minimal infinite emitter or a minimal sink. If $C$ is minimal infinite emitter, since $C\subseteq A\cup B$, we have that $C\cap A$ or $C\cap B$ is an infinite emitter. Suppose, without loss of generality, that $C\cap A$ is an infinite emitter. Since $C$ is minimal, we have $C\cap A=C$ and hence $C\subseteq A$. So $x\in X_A$, i.e $X_{A\cup B}\subseteq X_A\cup X_B$. On the other hand, if $C$ is a minimal sink, since $C\subseteq A\cup B$, we have $|C\cap A|=\infty$ or $|C\cap B|=\infty$. Suppose, without loss of generality, that $|C\cap A|=\infty$. Since $C$ is minimal we must have $C\cap A=C$, and hence $C\subseteq A$. Then $x\in X_A$. Therefor $X_{A\cup B}\subseteq X_A\cup X_B$.
	\end{proof}
	
	\begin{corolario}
		\label{corolario_4.9}
		For each $A\in \G^0$, the set $X_A$ is clopen and compact in $X$. 
	\end{corolario}
	\begin{proof}
		For vertices of the form $v\in G^0$ we have, by Lemma~\ref{vertice_generalizado_clopen_compacto} that $X_v$ is clopen and compact. For a generalized vertex of the form $r(e)$, notice that $X_{r(e)}=X_{e^{-1}}$ for all $e\in \G^1$. Thus, by Proposition \ref{X_c_clopen_compacto}, $X_{r(e)}=X_{e^{-1}}$ is clopen and compact. Finally, by Proposition \ref{prop_vert_gen} and by Lemma \ref{lema_4.8}, we have that $X_A$ is clopen and compact for every generalized vertex $A\in \G^0$. 
	\end{proof}
	
	\subsection*{The partial action}
	
	With the sets defined above, we have almost all ingredients necessary to define a topological partial action. We are missing homeomorphisms between nonempty subsets of the form $X_c$, with $c\in \F$. We define these maps below. 
	
	\begin{defi}
		\label{def_acao_par}
		For $a,b\in P$, where $X_a$ and $X_{ab^{-1}}$ are nonempty sets, and $A$ is either a minimal infinite emitter, or a minimal sink, or a singleton that contains only a sink, define: 
		\vspace{-1cm}
		\begin{multicols}{3} 
			\begin{align*}
			\theta_a: X_{a^{-1}}&\rightarrow X_a\\ 
			x&\mapsto ax\\
			(A,A)&\mapsto (a,A)\\
			(\beta,A)&\mapsto (a\beta,A)\\
			\gamma &\mapsto a\gamma
			\end{align*} \vfill
			\begin{align*}
			\theta_{a^{-1}}: X_a&\rightarrow X_{a^{-1}}\\
			ay&\mapsto \widehat{a}y=y\\
			(a,A)&\mapsto (A,A)\\
			(a\beta,A)&\mapsto (\beta,A)\\
			a\gamma &\mapsto \gamma
			\end{align*} 
			\vfill
			\begin{align*} 
			\theta_{ab^{-1}}: X_{ba^{-1}}&\rightarrow X_{ab^{-1}}\\
			bx&\mapsto a \widehat{b}x=ax\\
			(b,A)&\mapsto (a,A)\\
			(bc,A)&\mapsto (ac,A)\\
			b\gamma &\mapsto a\gamma
			\end{align*}
		\end{multicols} 
	where $x,y,\gamma \in X$.
	\end{defi} 
	Notice that we are denoting $a$ instead of $(a,r(a))\in \p$ and hence $\theta_a(x)=a.x=(a,r(a)).x$ for $x\in X_{a^{-1}}$. Moreover, the symbol ``\textasciicircum'' (hat) over a letter $a$ is understood as notation for ``remove'' the letter $a$ from the beginning of the path. Also notice that $\theta_a^{-1}=\theta_{a^{-1}}$ and $\theta_{ab^{-1}}^{-1}=\theta_{ba^{-1}}$.
	
	\begin{prop}
		For each $c\in\F$, the map $\theta_c: X_{c^{-1}}\rightarrow X_c$, as defined above, is a homeomorphism.
	\end{prop}
	\begin{proof}
		Let $a,b\in P$. Let us show that $\theta_{ab^{-1}}: X_{ba^{-1}}\rightarrow X_{ab^{-1}}$ is continuous. Let $(x^n)_{n\in \N}\subseteq X_{ba^{-1}}$ be a sequence such that $x^n\rightarrow x\in X_{ab^{-1}}$. Since $x^n\in X_{ba^{-1}}\subseteq X_b$ we have $x^n=b\gamma_1^n\gamma_2^n\ldots$, or $x^n=(b\gamma_1^n\ldots\gamma_{k_n}^n,A_n)$, or $x^n=(b\gamma_1^n\ldots\gamma_{k_n}^n,\{v_n\})$, for each $n\in \N$. 
		\begin{itemize}
			\item If $x\in \p^\infty$, say $x=b\gamma_1\gamma_2\ldots$, then for each $M\in \N$ there exists an $n_0\in\N$ such that $|x^n|>M+|b|$ and $b\gamma_1^n\ldots\gamma_M^n=b\gamma_1\gamma_2\ldots\gamma_M$, $\forall n\geq n_0$. For each $n\in \N$, consider $y^n=\theta_{ab^{-1}}(x^n)$ and $y=\theta_{ab^{-1}}(x)$. Then, for each $n$, we have $y^n=a\gamma_1^n\gamma_2^n\ldots$, $y^n=(a\gamma_1^n\ldots\gamma_{k_n}^n,A_n)$, or $y^n=(a\gamma_1^n\ldots\gamma_{k_n}^n,\{v_n\})$. Since $a\gamma_1^n\ldots\gamma_M^n=a\gamma_1\gamma_2\ldots\gamma_M$ for each $n\geq n_0$, by Corollary~\ref{corolario_convergencia} we have $y^n\rightarrow a\gamma_1\gamma_2\ldots =y$. 
			\item If $x\in X_{min}$, or $x\in X_{sin}$ with minimal range, then $x=(b,B)$ or $x=(b\gamma_1\ldots\gamma_{|x|-|b|},B)$:
			\begin{enumerate}[$\circ$]
				\item If $x=(b,B)$ then, for each finite set $F\subseteq \varepsilon(B)$ and for each finite set $S\subseteq B$, there exists $N\in \N$ such that if $n>N$ then $x^n=x$, or $|x^n|>|x|=|b|$ and  $\gamma_{|b|+1}^n \in \varepsilon(B)\backslash F$, or $|x^n|=|x|$ with $x^n=(b,\{v_n\})$ such that $v_n\in B\backslash S$. For each $n\in \N$, let $y^n=\theta_{ab^{-1}}(x^n)$ and $y=\theta_{ab^{-1}}(x)=(a,B)$. Hence, for each $F\subseteq \varepsilon(B)$ and for each finite set $S\subseteq B$, there exists $N\in \N$ such that if $n>N$ then $y^n=y$, or $|y^n|>|y|=|a|$ and  $\gamma_{|a|+1}^n \in \varepsilon(B)\backslash F$, or $|y^n|=|y|$ with $y^n=(a,\{v_n\})$ such that $v_n\in B\backslash  S$. Therefore, by Corollary~\ref{corolario_convergencia}, $y^n\rightarrow y$. 
				\item If $x=(b\gamma_1\ldots\gamma_{|x|-|b|},B)$, one proceeds analogously to the previous sub item, changing $b$ for $b\gamma_1\ldots\gamma_{|x|-|b|}$.
			\end{enumerate}
			\item If $x\in X_{sin}$, where $x=(b,\{v\})$ or $x=(b\gamma_1\ldots\gamma_{|x|-|b|},B)$, then $x^n$ is eventually constant. Thus, $y^n:=\theta_{ab^{-1}}(x^n)$ will also be eventually constant and obviously $y^n\rightarrow y:=\theta_{ab^{-1}}(x)$.
		\end{itemize} 
		We conclude that $\theta_{ab^{-1}}$ is continuous. Since $\theta_{ab^{-1}}^{-1}=\theta_{ba^{-1}}$, it follows that $\theta_{ab^{-1}}$ is a homeomorphism. The proof that $\theta_a$ is a homeomorphism is analogous.
	\end{proof}
	
	A topological partial action induces a partial action on the C* level (see \cite{Exel_livro}), as we describe in the following corollary.
	
	\begin{corolario} 
		For each $t\in\F$, the map 
		\begin{align*}
		\alpha_t: C(X_{t^{-1}})&\rightarrow C(X_t)\\
		f&\mapsto f\circ \theta_{t^{-1}}
		\end{align*} is a $*-$isomorphism. 
		Moreover, the collection $(\{X_t\}_{t\in\F},\{\theta_t\}_{t\in\F})$ is a topological partial action. Consequently, $(\{C(X_t)\}_{t\in\F},\{\alpha_t\}_{t\in\F})$ is a  C*-algebraic partial action of $\F$ in $C_0(X)$.
	\end{corolario}
	
	\begin{obs}
		When $\G$ is an ultragraph without sinks (but satisfies Condition~(RFUM2)), the partial action defined above coincides with the one defined in \cite{danidani}. Moreover, if $\G$ is a directed graph then the partial action above coincides with the one introduced in \cite{CarLar}. 
	\end{obs}
	
	\subsection*{$C^*(\G)$ as a partial crossed product} 
	
We finish the section describing the isomorphism between an ultragraph C*-algebra and the partial crossed product induce by the partial action defined above.  Before we proceed we make a couple of observations. 
	
\begin{obs} Notice that, for every $A\in \G^0$, the characteristic map $1_A$ of the set $X_A$ is an element of $C(X)$ (recall that by Corollary \ref{corolario_4.9} $X_A$ is open and compact). Furthermore, for each $c\in \F$, by Proposition~\ref{X_c_clopen_compacto}, the set $X_c$ is open and compact and so the characteristic map $1_c$ of $X_c$ is also an element of $C(X)$.
\end{obs}
	
	\begin{lema}
	\label{lema3}
	The subalgebra $D\subseteq C_0(X)$ generated by elements of the form $1_c$, $1_A$, and $\alpha_c(1_{c^{-1}} 1_A)$, with $c\in \displaystyle\bigcup_{n=1}^{\infty} (\G^1)^n$ and $A\in \G^0$, is dense in $C_0(X)$. Moreover, for each $0\neq g\in \F$, the subalgebra $D_g\subseteq C(X_g)$ generated by all the maps $1_g1_c$, $1_g1_A$, and $1_g\alpha_c(1_{c^{-1}}1_A)$ is dense in $C(X_g)$.    
	\end{lema}
	
\begin{proof}
Clearly $D$ is a self adjoint subalgebra. Density follows from the Stone-Weierstrass Theorem. To check that $D$ vanishes nowhere and separate points one proceeds as in \cite[Lemma 4.11]{danidani}, using characteristic functions $1_\alpha$, $1_v$ and $1_A$. An analogous proof holds for $D_g$.

	\end{proof}

    The following theorem is the main result of this section. We have set up the ground so that the proof of it is similar to the one presented in \cite[Thm. 4.12]{danidani} (applying Lemmas~\ref{lema3} and \ref{lema_4.8}). We therefore omit the proof, but more details can be found in \cite[Thm. 2.15]{tesefelipe}.   
	
	\begin{teorema}
		\label{teo_1}
		Let $\G$ be an ultragraph that satisfies Condition~(RFUM2). Then there exists a $*$-isomorphism $\Phi:C^*(\G)\rightarrow C_0(X)\rtimes_\alpha\F$ such that $\Phi(s_e)=1_e\delta_e$, for each $e\in \G^1$, and $\Phi(p_A)=1_A\delta_0$, for each $A\in\G^0$. 
	\end{teorema}

\section{Shift space and ultragraph groupoid}\label{gpd}
	
It is known that every partial crossed product can be seen as a groupoid algebra, but the concrete realization is important for many applications. So in this section we use Theorem~\ref{teo_1} to realize ultragraph C*-algebras as groupoid C*-algebras and define the shift space associated to an ultragraph as the unit space of the associated groupoid, attached with the shift map (which we define below). 

\begin{obs} In \cite{mu&mu} ultragraph C*-algebras are described as groupoid C*-algebras in general, but the lack of Condition~(RFUM2) implies in a more technical topology, which is more difficult to use in applications such as the ones studied in the following sections.
\end{obs}

	
	
	\subsection*{Ultragraph Shift Space}

In \cite[Section 3.2]{danidani} the shift map is defined even for elements of lentgh zero of the boundary path space $X$. For the applications that follow it will be better to define the shift map only for elements of length greater than zero. To make this precise we need the following notation. 

\begin{defi} 
 For $n\in\{0,1,2,\ldots\}$, define the following subsets of $X=X_{fin}\cup\p^\infty$:
$$X^n=\p^n\cap X_{fin}=\{x\in X_{fin}:|x|=n\}; 
\ X^{\geq n}=\displaystyle\bigcup_{k\geq n} X^k
\text{ and }
 X_{\infty}^{\geq n}=X^{\geq n}\cup \p^\infty .$$
	\end{defi}

Notice that $X_{\infty}^{\geq 1}=X\backslash X^0$ and that $X_{\infty}^{\geq 0}=X$.
	
	\begin{defi}
	Let $X$ be the boundary ultrapath space associated to an ultragraph $\G$. We define the shift map $\sigma:X_{\infty}^{\geq 1}\rightarrow X$ by:
		
		$\sigma(x)=
		\left\{
		\begin{array}{ll}
		\gamma_2\gamma_3\ldots, & \text{if $x=\gamma_1\gamma_2\ldots \in \p^\infty$};\\
		(\gamma_2\ldots\gamma_n,A), & \text{if $x=(\gamma_1\ldots\gamma_n,A) \in X^{\geq 2}$};\\
		(A,A), & \text{if $x=(\gamma_1,A) \in X^1$}.\\
				\end{array}
		\right.$  
		
		For $n\geq 1$ we define $\sigma^n$ as the composition $n$ times of $\sigma$, and for $n=0$ we define $\sigma^0$ as the identity. When we write $\sigma^n(x)$ we are implicitly assuming that $x\in X_{\infty}^{\geq n}$. 
		
	\end{defi}
	
Next we show that the shift map is a local homeomorphism on elements of length greater than zero. This will be important when we see the ultragraph groupoid as a Deaconu-Renault groupoid.
	
	\begin{prop}
		The shift map defined above is continuous at all points of $X_{\infty}^{\geq 1}$. In addition, if $x\in X\backslash X^0$ then there exists an open set $U$ of $X$ that contains no elements of length zero, and such that $x\in U$, $\sigma(U)$ is an open subset of $X$, and $\sigma|_U :U\rightarrow \sigma(U)$ is a homeomorphism. 
	\end{prop}
	\begin{proof}
		Let $\{x^n\}_{n\in\N}\subseteq X$ be a sequence converging to $x\in X$. Then it is clear that $x^n=(\gamma_1^n\ldots\gamma_k^n,A_n)$, or $x^n=(\gamma_1^n\ldots\gamma_k^n,\{v_n\})$ or $x^n=\gamma_1^n\gamma_2^n\ldots$, for all $n\in\N$.
		
		If $|x|=\infty$ it is obvious that $\sigma(x^n)$ converges to $\sigma(x)$ by the characterization of convergence given in the Corollary \ref{corolario_convergencia}. If $|x|<\infty$ and $x\in X_{sin}$ with $r(x)\in G_s^0$, we have that $x^n$ is eventually constant , and thus $\sigma(x^n)$ is also eventually constant. Therefore $\sigma(x^n)$ converges to $\sigma(x)$. 
		
		If $|x|<\infty$ and $x=(\gamma_1\ldots\gamma_k,A)\in X_{min}$, by item $(ii)$ of Corollary \ref{corolario_convergencia}, given $F\subseteq \varepsilon(A)$ and $S\subseteq A$ there exists $N\in\N$ such that: if $n\geq N$ then either $x^n=x$ or  $|x^n|>k$, $\gamma_i^n=\gamma_i$ for all $1\leq i \leq k$, and $\gamma_{k+1}^n\in\varepsilon(A)\backslash F$. Hence, for all $n\geq N$, $\sigma(x^n)= \sigma(x)$ or $|\sigma(x^n)|>k-1$ and $\sigma(x^n)_k=\gamma_{k+1}^n\in \varepsilon(A)\backslash F$, i.e $\sigma(x^n)$ converges to $\sigma(x)$ and so $\sigma$ is continuous. 
		
		If $|x|<\infty$ and $x=(\gamma_1\ldots\gamma_k,A)\in X_{sin}$ with $r(x)\in G_s^0$, by item $(iii)$ of Corollary \ref{corolario_convergencia}, for all $S\subseteq A\cap G_s^0$ there exists $N\in\N$ such that: if $n\geq N$ then either $x^n=x$ or we have that $|x^n|=k$, $\gamma_i^n=\gamma_i$ for all $1\leq i \leq k$, and $r(x^n)=\{v_k\}\in A\backslash S$. Hence, for all $n\geq N$, $\sigma(x^n)= \sigma(x)$ or $|\sigma(x^n)|=k-1$ and $r(\sigma(x^n)_{k-1})=r(\gamma_{k}^n)=\{v_k\} \in A\backslash S$, i.e $\sigma(x^n)$ converges to $\sigma(x)$ and so $\sigma$ is continuous. 
		
		Finally, notice that if $|x|\geq 1$ and $U$ is one of the basic neighborhoods of Proposition~\ref{prop_base_top} that contains $x$, then $\sigma$ is a homeomorphism between $U$ and $\sigma(U)$.
	\end{proof}
	
	We now make a precise definition of the shift space associated to an ultragraph.
	
	\begin{defi}\label{esp_shift}	
	Let $\G$ be an ultragraph. We define the pair $(X,\sigma)$ to be the \emph{shift space} associated to the ultragraph $\G$, where $X$ is the boundary ultrapath space (with the topology given in Propostion~\ref{prop_base_top}) and $\sigma$ is the shift map defined above.
	\end{defi}
	
	\begin{obs}
		If $\G$ is a finite directed graph without sinks then $\G$ has no infinite emitters, and so $X$ is the usual edge shift space associated with the graph. If $\G$ is any graph, our definition of shift space coincides with that in \cite{BCW}, and if $G$ is an ultragraph without sinks our definition coincides with the one in \cite{danidani}.
	\end{obs}
	

	\subsection*{Ultragraph Groupoids}
In this subsection we define a groupoid associated to an ultragraph (that satisfies Condition~(RFUM2)) such that the associated C*-algebra is isomorphic to the ultragraph C*-algebra. There are two ways to prove the later result, once the groupoid is defined: one is to give a proof via core subalgebras, analogous to the graph case (see for example \cite{KPRR}, \cite{Ark_Eil_Rui}, and mainly \cite{BCW}). The other way is to notice that the groupoid we define is isormorphic to the transformation groupoid given by the partial action of Definition~\ref{def_acao_par}. By \cite{Abadie}, we know that the C*-algebra of the transformation groupoid is isomorphic to the partial crossed product (which in turn is isomorphic to the ultragraph C*-algebra by Theorem \ref{teo_1}). We need the following definitions:

	\begin{defi}\label{grupoide}
		Let $\G$ be an ultragraph (satisfying (RFUM2)). We define the \emph{groupoid associated to $\G$} by $G(X,\sigma):=\{(x,m-n,y):x,y\in X;m,n\in \N;\sigma^m(x)=\sigma^n(y) \}$. The source function is defined by $s(x,k,y)=y$ and the range by $r(x,k,y)=x$. The multiplication is given by $(x,k,y)(y,l,z)=(x,k+l,z)$ and the inversion is given by $(x,k,y)^{-1}=(y,-k,x)$. Thus, the unit space is $\{(x,0,x):x\in X\}\cong X$ and the topology considered is that of cylinders (as in \cite{BCW} and in \cite[8.3]{Ort_Nyl}), i.e the topology generated by the basic sets 
			$$Z(U,m,n,V)=\{(x,k,y)\in G(X,\sigma): x\in U, k=m-n, y\in V; \sigma^m(x)=\sigma^n(y)\},$$
		where $U\subseteq X^{\geq m}_\infty$, $V\subseteq X^{\geq n}_\infty$ are open such that $\sigma_{|_U}^m$ and $\sigma_{|_V}^n$ are injective and $\sigma^m(U)=\sigma^n(V)$. 
	\end{defi}
	\begin{obs}
	As we will see in Section~\ref{sec3.3.2}, the groupoid above coincides with the Deaconu-Renault groupoid associated to the shift space $X$. Also note that this groupoid (for ultragraph without sinks) has been studied recently in connection with full groups, see \cite{danielgillesdanie}.
	\end{obs}
	
	To define the transformation groupoid, we will rely in \cite{Abadie}. 
	
	\begin{defi}
		Consider $\G$ an ultragraph that satisfies (RFUM2) and the partial action of Definition \ref{def_acao_par}. \emph{The transformation groupoid associated with $\G$} is given by  $\G_T:=\{(x,c,y)\in X\times\F\times X;x=\theta_c (y)=cy\}$. The source function is given by $s(x,c,y)=y$ and the range by $r(x,c,y)=x=cy$. The multiplication $(x,c,y)(w,d,z)$ is only defined when $s(x,c,y)=r(w,d,z)\iff y=w=dz$, and thus $(x,c,y) (y,d,z):=(x,cd,z)$. The inversion is given by $(x,c,y)^{-1} = (y,c^{-1},x)$. The unit space is $X$. The topology is the product topology inherited from $X\times\F\times X$. As $\F$ is discrete we have that the groupoid is étale.
	\end{defi}
	
	\begin{obs}
		Notice that $\G_T\cong	X\rtimes_\theta \F=\{(c,x)\in \F\times X;c\in \F,x\in X_c\}$, where the source and range are  defined by $s(c,x)=x$ and $r(c,x)=cx$. The multiplication is given by $(a,bx)(b,x)=(ab,x)$, and the inversion by $(a,x)^{-1}=(a^{-1},ax)$. The transformation groupoid built in \cite{Abadie} has this format.  
	\end{obs}
	
	Let us then show that the two groupoids defined above are isomorphic, thus achieving an isomorphism between the ultragraph C*-algebra and the groupoid C*-algebra in a different way. Before we give the proof we need one auxiliary result.
	
\begin{lema}\label{lemautil}
	An element in $\G_T$ is always of the form $\big((\alpha,A).x,\alpha\beta^{-1},(\beta,A).x\big)$, with $x\in X$.
\end{lema}		
\begin{proof}
 Indeed, suppose that $\big((\alpha,A),\gamma\delta^{-1},(\beta,B)\big)\in \G_T$. Then we have $\theta_{\gamma\delta^{-1}}(\beta,B)=(\alpha,A)$, i.e $(\gamma\widehat{\delta}\beta,B)=(\alpha,A)$. Hence $B=A$ and $\gamma\widehat{\delta}\beta=\alpha$. Thus, by the definition of elements of $\G_T$, we have that $\alpha=\gamma\alpha'$ and $\beta=\delta\beta'$. So  $\gamma\alpha'=\alpha=\gamma\widehat{\delta}\beta=\gamma\widehat{\delta}\delta\beta'=\gamma\beta'$, i.e $\alpha'=\beta'$, and the statement follows since $\alpha\beta^{-1}=\gamma\alpha'(\delta\beta')^{-1}=\gamma\alpha'(\beta')^{-1}\delta^{-1}=\gamma\delta^{-1}$.
 Therefore $\theta_{\alpha\beta^{-1}}(\beta,A)=(\alpha,A)$.
\end{proof}
	\begin{teorema}
\label{teo_grupoides}
		Let $\G$ be an ultragraph that satisfies (RFUM2). Consider $G(X,\sigma)$ and $\G_T$ the groupoids defined above. If $ab^{-1}$ is in reduced form, where $x=ax'$, $y=bx'$, and $x=\theta_{ab^{-1}}(y)$, then the map 
		\begin{align*}
		\Psi :\G_T&\rightarrow G(X,\sigma) \\
		(x,ab^{-1},y)&\mapsto (x,|a|-|b|,y)
		\end{align*} 
 is a topological $*$-isomorphism of groupoids. 
	\end{teorema}
	\begin{proof}
		First notice that, by the previous lemma, we need $ab^{-1}$ to be in reduced form for $\Psi $ to be well defined. 
		Let $(ax,ab^{-1},bx)$, $(bx,bc^{-1},cx)\in \G_T$ (compare with Lemma~\ref{lemautil}). Then
		\begin{align*}
		\Psi \big((ax,ab^{-1},bx).&(bx,bc^{-1},cx)\big) = \Psi(ax,ab^{-1}bc^{-1},cx) = (ax,|a|-|c|,cx)= \\
		&=(ax,|a|-|b|+(|b|-|c|),cx)= (ax,|a|-|b|,bx).(bx,|b|-|c|,cx)=\\
		&=\Psi(ax,ab^{-1},bx).\Psi(bx,cd^{-1},cx),
		\end{align*} and hence $\Psi$ is multiplicative.
		Also, $\Psi$ preserves inverse, since $$\Psi\big((x,ab^{-1},y)^{-1}\big)=\Psi(y,ba^{-1},x)=(y,|b|-|a|,x)=(x,|a|-|b|,y)^{-1}=\Psi\big((x,ab^{-1},y)\big)^{-1}.$$ 
		
	Next we show that $\Psi$ is injective. Let $(x,ab^{-1},y)$ and $(w,cd^{-1},z)$ be elements of $\G_T$ such that $\Psi(x,ab^{-1},y)=\Psi (w,cd^{-1},z)$, i.e $(x,|a|-|b|,y)=(w,|c|-|d|,z)$. Thus, $x=w$, $y=z$, $|a|-|b|=|c|-|d|$, $\sigma^{|a|}(x)=\sigma^{|b|}(y)$ and $\sigma^{|c|}(w)=\sigma^{|d|}(z)$. 
		We need to show that $ab^{-1}=cd^{-1}$. But we know that $\theta_{ab^{-1}}(y)=a\widehat{b}y=x=w=c\widehat{d}z=c\widehat{d}y$. Thus either $a=ca'$ or $c=ac'$ or $a=c$. Suppose, without loss of generality, that $c=ac'$. Hence, 
		\begin{equation}
		\widehat{b}y=\sigma^{|a|}(a\widehat{b}y)=\sigma^{|a|}(x)=\sigma^{|a|}(w)=\sigma^{|a|}(c\widehat{d}y)=\sigma^{|a|}(ac'\widehat{d}y)=c'\widehat{d}y.\label{eq3.1}
		\end{equation} 
		But since $y\in X_{ba^{-1}}$, we have that $y=by'$. Thus, by (\ref{eq3.1}), we have 
		\begin{equation}
		y'=\widehat{b}by'=\widehat{b}y=c'\widehat{d}by'.\label{eq3.2}
		\end{equation} 
		On the other hand, since $c=ac'$, we have $|c|\geq |a|$ and  so $|d|-|b|=|c|-|a|\geq 0$, which implies that $|d|\geq |b|$. Since $y=z\in X_{dc^{-1}}$, we have that $d=bd'$. Hence, by (\ref{eq3.2}), 
		\begin{equation}
		y'=c'\widehat{d}by'\underbrace{=}_{d=bd'}c'\widehat{d'}\widehat{b}by'=c'\widehat{d'}y',\label{eq3.3}
		\end{equation} 
		i.e, $y'=c'\widehat{d'}y'=\theta_{c'd'^{-1}}(y')$. Now, note that if $c'\neq d'$ then $c'd'^{-1}$ is irreducible, and $\theta_{c'd'^{-1}}$ cannot contain fixed points. Thus the existence of such $y'$ is a contradiction. Therefore we must have $c'=d'$, and so $cd^{-1}=ac'd^{-1}=ac'd'^{-1}b^{-1}=ab^{-1}$, and $\Psi$ is injective.

		Now, we prove that $\Psi$ is surjective. Let $(x,n-m,y)\in G(X,\sigma)$, where $x,y\in X$ are such that $\sigma^n(x)=\sigma^m(y)$. Suppose, without loss of generality, that $n\geq m$. Then $x=x_1\ldots x_m \ldots x_n e_1e_2\ldots$ and $y=y_1\ldots y_m e_1e_2\ldots$. Consider $a=x_1\ldots x_n$ and $b=y_1\ldots y_m$. It is clear that $(x,ab^{-1},y)\in \G_T$ and $\Psi(x,ab^{-1},y)=(x,|a|-|b|,y)=(x,n-m,y)$. 
		
		It remains to be shown that $\Psi$ is continuous and open. Before we do this, we describe the elements in the ultragraph groupoid: Let $(x,m-n,y)\in G(X,\sigma)$. Notice that if $x\in X_{fin}$ then $y\in X_{fin}$ and $ (x,m-n,y)=((\alpha,A),m-n,(\beta,B))$, and also $\sigma^m(\alpha,A)=\sigma^n(\beta,B)$, i.e, $(\alpha_{m+1}\ldots \alpha_{|\alpha|},A)=(\beta_{n+1}\ldots \beta_{|\beta|},B)$, that is, $A=B$ and $\alpha_{m+1}\ldots \alpha_{|\alpha|}=\beta_{n+1}\ldots \beta_{|\beta|}$. Hence $|\alpha|-(m+1)+1=|\beta|-(n+1)+1$, i.e, $|\alpha|-|\beta|=n-m$. This happens with $A,B $ both being minimal infinite emitters, or with both being sinks. On the other hand, if $x\in \p^{\infty}$ then we must have that $y\in\p^{\infty}$, and hence $(x,m-n,y)=(\gamma_1\gamma_2\ldots,m-n,\delta_1\delta_2\ldots)$ and also $\sigma^m(\gamma_1\gamma_2\ldots)=\sigma^n(\delta_1\delta_2\ldots)$, i.e, $\gamma_{m+1}\gamma_{m+2}\ldots=\delta_{n+1}\delta_{n+2}\ldots$. Therefore, elements in $G(X,\sigma)$ are either of the form $((\alpha,A)y,m-n,(\beta,A)y)$, with $y\in X$ and $|\alpha|-|\beta|=n-m$, or of the form $(\gamma_1\ldots\gamma_my,m-n,\delta_1\ldots \delta_n y)$ with $y\in X$. 
		
	Before we proceed, we recall the basic open sets of the topologies of each space. We know by Definition~\ref{grupoide} that the basic open sets of $G(X,\sigma)$ are of the form: 
		$$Z(U,m,n,V)=\{(x,k,y)\in G(X,\sigma): x\in U, k=m-n, y\in V; \sigma^m(x)=\sigma^n(y)\},$$
		where $U\subseteq X^{\geq m}_\infty$, $V\subseteq X^{\geq n}_\infty$ are open such that $\sigma_{|_U}^m$ and $\sigma_{|_V}^n$ are injective and $\sigma^m(U)=\sigma^n(V)$. But, like in \cite[8.3]{Ort_Nyl}, we will use a coarser subcollection of this basis, which is still a basis for the ultragraph groupoid. Namely, we consider the collection $\{Z(U,|\alpha|,|\beta|,V):\sigma^{|\alpha|}(U)=\sigma^{|\beta|}(V)\}$, where $(\alpha,A)$ and $(\beta,A)$ in $\p$, $U\subseteq D_{(\alpha,A)}$ are compact open,  $V\subseteq D_{(\beta,A)}$, 
		and $Z(U,|\alpha|,|\beta|,V)=\{((\alpha,A)y,|\alpha|-|\beta|,(\beta,A)y):y\in X\}$ (for more details see \cite[Lemma~3.4]{danielgillesdanie}). 
	
		On the other hand, to analyze the basic opens in $\G_T$, we use the characterization of its elements made in the previous lemma. Since the topology in $\G_T$ is inherited from $X\times \F\times X$, 
		given $(\alpha,A),(\beta,A)\in \p$, an open basic set is $$D_{\alpha,\beta}^A=\big\{\big( (\alpha,A)x,\alpha\beta^{-1},(\beta,A)x\big): s(x)\in A, x\in X \big\}.$$
		Other open sets, when $F\subseteq\varepsilon(A)$ and $S\subseteq A\cap G_s^0$ are finite, are of the form
		$$D_{\alpha,\beta}^{A,F,S}=\big\{\big( (\alpha,A)x,\alpha\beta^{-1},(\beta,A)x\big): s(x)\in A, x\in X, x_1\in \varepsilon(A)\backslash F,\big\}\cup\big\{\big( (\alpha,v),\alpha\beta^{-1},(\beta,v)\big);v \notin S\big\}.$$
		That is, these open set depend on $(\alpha,A)$ and $(\beta,A)$, and eventually on $F$ and $S$. We show that $\psi$ takes open sets of the first kind to open sets in $G(X,\sigma)$ (the other cases are analogous). To do this, notice that
		\begin{align*}
		\Psi\Big(D_{\alpha,\beta}^A\Big)&=\Psi\Big(\big\{\big( (\alpha,A)x,\alpha\beta^{-1},(\beta,A)x\big): s(x)\in A, x\in X\big\}\Big)=\\
		&=\big\{\Psi\big( (\alpha,A)x,\alpha\beta^{-1},(\beta,A)x\big): s(x)\in A, x\in X\big\}=\\
		&=\big\{\big( (\alpha,A)x,|\alpha|-|\beta|,(\beta,A)x\big): x\in X\big\}=Z\big(U,|\alpha|,|\beta|,V\big),
		\end{align*}
		with $U\subseteq D_{(\alpha,A)}$ open compact, $V\subseteq D_{(\beta,A)}$ open compact, and $\sigma^{|\alpha|}(U)=\sigma^{|\beta|}(V)$. Thus, $\Psi\Big(D_{\alpha,\beta}^A\Big)$ is a basic open set in $G(X,\sigma)$. 
		
		On the other hand, let $(\alpha,A),(\beta,A)\in \p$, $U\subseteq D_{(\alpha,A)}$ be open compact and $V\subseteq D_{(\beta,A)}$ be open compact. We describe the inverse image of the cylinder $Z(U,|\alpha|,|\beta|,V)$, i.e $\Psi^{-1}(Z(U,|\alpha|,|\beta|,V))$. For this, note that $$Z(U,|\alpha|,|\beta|,V)=\{((\alpha,A)y,|\alpha|-|\beta|,(\beta,A)y):y\in X:\sigma^{|\alpha|}(U)=\sigma^{|\beta|}(V)\}.$$
		Also, notice that $((\alpha,A)y,\alpha\beta^{-1},(\beta,A)y)\in \G_T$ is such that $$\Psi ((\alpha,A)y,\alpha\beta^{-1},(\beta,A)y) =((\alpha,A)y,|\alpha|-|\beta|,(\beta,A)y)$$. Since $\Psi$ is injective, $((\alpha,A)y,\alpha\beta^{-1},(\beta,A)y)$ is the unique element that is taken to $((\alpha,A)y,|\alpha|-|\beta|,(\beta,A)y)$. Thus, we have that $$\Psi^{-1}(Z(U,|\alpha|,|\beta|,V))=\{((\alpha,A)y,\alpha\beta^{-1},(\beta,A)y):y\in X\}\subseteq X\times\F\times X$$
		is a basic open set in $\G_T$, and it follows that $\Psi$ is a homeomorphism, and hence an isomorphism of groupoids.
	\end{proof}

\section{Continuous orbit equivalence and eventual conjugacy between ultragraph shift spaces}	

In this section we describe both continuous orbit equivalence and eventual conjugacy between ultragraph shift spaces. We will rely on the general results regarding Deaconu-Renault systems of \cite{DeaRen}. In fact, while developing this section, our first goal was to generalize to ultragraphs Theorem~4.2 of \cite{CarlsenLie}, which states that two graph groupoids $\G_E$ and $\G_F$ are isomorphic if, and only if, there is a continuous orbit equivalence $h:\partial E \rightarrow \partial F$ that preserves isolated eventually periodic points. As we developed our proof the paper \cite{DeaRen} on orbit equivalence for Deaconu-Renault systems was posted on arXiv. We then changed our strategy and decided to apply the results of \cite{DeaRen} to our context. We notice that the application of the results in \cite{DeaRen} to the ultragraph case are not straightforward, and in fact we need to recall and develop a few concepts regarding ultragraphs. We do this below.

Notice that from the definition of initial segment made on page \pageref{seg_ini}, we can decompose an ultrapath of length $n$ in $n$ ultrapaths of length $1$, that is, if $|(\alpha,A)|=n$, we can write $(\alpha,A)=(\alpha_1,A_1)(\alpha_2,A_2)\ldots(\alpha_n,A)$, where $s(\alpha_{i+1})\in A_i\subseteq r(\alpha_i)$ and $|(\alpha_i,A_i)|=1$. We can now make the following definition.

\begin{defi}
A \emph{cycle} in $\G$ is an ultrapath $(\alpha,A) \in \p \backslash \p^0$ such that $s(\alpha) \in A$. An \emph{exit} to a cycle $(\alpha,A)=(\alpha_1, A_1)(\alpha_2, A_2) \ldots (\alpha_n, A_n) $ (where $|(\alpha_i,A_i)|=1$ for all $i$) is one of the following:
\begin{enumerate}[i)]
	\item an ultrapath $(e,E)\in\p^1$ such that $\exists$ $i$ for which $s(e)\in r(\alpha_i)$, but $e\neq \alpha_{i+1}$. (Note that it is not necessary that $s(e)\in A_i$ for some $i$).
	\item a sink $v\in\G^0$ such that $v\in r(\alpha_i)$ for some $i$.
\end{enumerate} 
In this case we say that $x=(\alpha,A)$ is a cycle of size $n$. When $n=1$, we call such cycle a loop.
\begin{obs}
Note that the concatenation of  cycles remains a cycle. Also we remark that in \cite{Tomforde2} the name loop was used for what we call a cycle.
\end{obs}
\begin{defi}
When a cycle is not the concatenation of  smaller cycles we say that this cycle is \emph{simple}.
\end{defi}

\end{defi}
\begin{exemplos} 
The following three ultragraphs contain cycles.
\vspace*{-0.8cm}
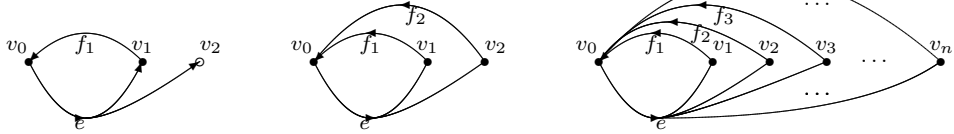
\begin{center}
\setlength{\unitlength}{1.5mm}
\begin{picture}(60,20)

\put(-10,11){\scriptsize$v_0$}
\put(-8,10){\circle*{0.7}}
\put(1,11){\scriptsize$v_1$}
\put(2,10){\circle*{0.7}}
\put(7,11){\scriptsize$v_2$}
\put(7,10){\circle{0.7}}

\put(-4,4){\scriptsize$e$}
\put(-3,5){\vector(1,0){0}}
\qbezier(-8,10)(-3,0)(2,10)
\qbezier(-3,5)(0,5)(7,10)
\qbezier(2,10)(-3,15)(-8,10)

\put(7,10){\vector(3,2){0}}
\put(1.9,10){\vector(1,3){0}}
\put(-8,10){\vector(-1,-1){0}}
\put(-4,11){\scriptsize$f_1$}

\put(15,11){\scriptsize$v_0$}
\put(17,10){\circle*{0.7}}
\put(26,11){\scriptsize$v_1$}
\put(27,10){\circle*{0.7}}
\put(32,11){\scriptsize$v_2$}
\put(32,10){\circle*{0.7}}

\put(21,4){\scriptsize$e$}
\put(22.5,5){\vector(1,0){0}}
\qbezier(17,10)(22,0)(27,10)
\qbezier(22,5)(25,5)(32,10)
\qbezier(27,10)(22,15)(17,10)
\qbezier(32,10)(25,20)(17,10)

\put(17,10){\vector(-1,-1){0}}
\put(21,11){\scriptsize$f_1$}
\put(25,13.6){\scriptsize$f_2$}
\put(21.5,12.5){\vector(-1,0){0}}
\put(24.5,15){\vector(-1,0){0}}

\put(40,11){\scriptsize$v_0$}
\put(42,10){\circle*{0.7}}
\put(52,11){\scriptsize$v_1$}
\put(52,10){\circle*{0.7}}
\put(56,11){\scriptsize$v_2$}
\put(57,10){\circle*{0.7}}
\put(61,11){\scriptsize$v_3$}
\put(62,10){\circle*{0.7}}
\put(71,11){\scriptsize$v_n$}
\put(72,10){\circle*{0.7}}
\put(65,10){\scriptsize$\ldots$}

\put(47,4){\scriptsize$e$}
\put(47.5,5){\vector(1,0){0}}
\qbezier(42,10)(47,0)(52,10)
\qbezier(47,5)(50,5)(62,10)
\qbezier(47,5)(50,5)(57,10)
\qbezier(47,5)(65,5)(72,10)		
\put(60,7){\scriptsize$\ldots$}
\put(60,15){\scriptsize$\ldots$}
\put(60,7){\scriptsize$\ldots$}

\qbezier(52,10)(47,15)(42,10)
\qbezier(57,10)(47,17)(42,10)
\qbezier(62,10)(50,20)(42,10)
\qbezier(72,10)(55,25)(42,10)	
\put(42,10){\vector(-1,-1){0}}
\put(46.5,12.5){\vector(-1,0){0}}
\put(48,13.5){\vector(-1,0){0}}
\put(50.5,15){\vector(-1,0){0}}
\put(55.5,17.5){\vector(-1,0){0}}

\put(46,11){\scriptsize$f_1$}
\put(50,12){\scriptsize$f_2$}
\put(52,13.6){\scriptsize$f_3$}
\put(55,16){\scriptsize$f_n$}
\end{picture}
\end{center}
\vspace*{-0.8cm}

Some examples of cycles in the first ultragraph are: $(ef_1,\{v_0\})$, $(f_1e,\{v_1\})$,$(f_1e,\{v_1,v_2\})$, while $(f_1e,\{v_2\})$ is not a cycle. The cycle $(ef_1,\{v_0\})$ is a simple cycle, while the cycle $(ef_1ef_1,\{v_0\})$ is not.

Some examples of cycles in the second ultragraph are: $(ef_1,\{v_0\})$, $(ef_1ef_2,\{v_0\})$, $(f_2e,\{v_1,v_2\})$, $(f_2ef_1e,\{v_2\})$, while $(f_2e,\{v_1\})$ is not a cycle. The cycle $(f_2ef_1e,\{v_2\})$ can be decomposed as follows: $(f_2ef_1e,\{v_2\})=(f_2e,\{v_1,v_2\})(f_1e,\{v_2\})$;

Some examples of cycles in the third ultragraph are: $(ef_i,\{v_0\})$, for all $i$, $(ef_1ef_2ef_3,\{v_0\})$, $(f_ie,\{v_i\})$, for all $i$, $(f_ie,r(e))$, for all $i$, while $(f_ie,\{v_j\})$ with $i\neq j$ is not a cycle. 

All cycles mentioned above contain exit. For example, in the first ultragraph $(e,\{v_2\})$ is an exit of $(f_1e,\{v_1\})$ and of $(f_1e,\{v_1,v_2\})$. In the second and third ultragraphs, $(e,\{v_j\})$ is always an exit of a cycle of the form $(f_ie,\{v_i\})$, when $j\neq i$.
\end{exemplos}

There is also a Condition (L) for ultragraphs:

\begin{defi}
An ultragraph $\G$ satisfies \emph{Condition~(L)} if all cycle in $\G$ has an exit.
\end{defi}

\subsection*{Isolated Points}

Isolated points play an important role in orbit equivalence. Therefore our next goal is to characterize the isolated points in ultragraph shift spaces. Recall that $(\alpha,A)\in X_{fin} $ is an isolated point in $X$ if, and only if, $A$ is a unitary set formed by a sink.

As with graphs, an infinite path $x\in\p^\infty$ is said to be a \emph{periodic point} if $x=\gamma\gamma\gamma\ldots$, where $\gamma\in\p\backslash\p^0$ is a cycle. Also, an infinite path $x\in\p^\infty$ is said to be an \emph{eventually periodic point} if $x=\mu\gamma\gamma\gamma\ldots$, where $\mu\in\p$, $\gamma\in\p\backslash\p^0$, and $\gamma$ is a cycle. Of course every periodic point is eventually periodic. 

Next we characterize when infinite paths that are not eventually periodic are isolated points.

\begin{prop}\label{isolados}
Let $x=x_1x_2x_2\dots\in \p^{\infty}$ be an infinite path that is not eventually periodic. Then $x$ is an isolated point if, and only if, the sets $S_1:=\{n\in\N:|\varepsilon(r(x_n))|\geq 2\}$ and $S_2:=\{n\in\N:|r(x_n)|\geq 2\}$ are finite. A path that is an eventually periodic point will be isolated if, and only if, its cycle has no exit.
\end{prop}
\begin{proof}
Suppose that $S_1$ or $S_2$ are infinite sets. Since $x$ is an infinite path, there is an initial segment $y$ of $x$ such that $x\in D_y$ (moreover, by definition of cylinders, any open $D_w$ which contains $x$ is such that $w$ is an initial segment of $x$). So we can consider $n\in S_1\cup S_2$, so that $n>|y|$. Then $|r(x_n)|\geq 2$ or $|\varepsilon(r(x_n))|\geq 2$, and therefore we have that while $x\in D_y$, either we have $(x_1x_2\ldots x_n,\{v_0\})\in D_y$ with $v_0\in r(x_n)\backslash s(x_{n+1})$, or we have $(x_1x_2\ldots x_ne,r(e))w\in D_y$ with $e \in \varepsilon(r(x_n))\backslash\{x_{n+1}\}$ and $s(w)\in r(e)$, $w\in X$. In other words, $D_y$ always contains more than one element. Since $y$ was arbitrarily chosen, there is no $y$ such that $D_y=\{x\}$ (or unions of $D_y$'s). Therefore $x$ is not isolated.

On the other hand, if $S_1$ and $S_2$ are finite, consider $n_0=\max\{S_1,S_2\}$. Then $y:=(x_1x_2\ldots x_{n_0+1},r(y))\in\p$ is an initial segment of $x$ and so $x\in D_y$. There is no other element in $D_y$ because if it existed we would have that either $|\varepsilon(r(x_m))|\geq 2$ or $|(r(x_m))|\geq 2$, for some $m>n_0$, a contradiction. Therefore $D_y=\{x\}$, i.e $x$ is isolated.

The second statement is demonstrated analogously to the above. Just note that an exit in a cycle $\gamma$ means that $|\varepsilon(r(\gamma_i))|\geq 2$ or $|(r(\gamma_i))|\geq 2$, for some $i=1,\ldots, |\gamma|$. 
\end{proof}
\begin{exemplo}
The points $x=e_1e_2e_3\ldots$ and $y=\mu_1\mu_2\mu_3\gamma_1\gamma_2\gamma_3\gamma_4\gamma_1\gamma_2\gamma_3\gamma_4\ldots$ are isolated in the ultragraphs below. \\

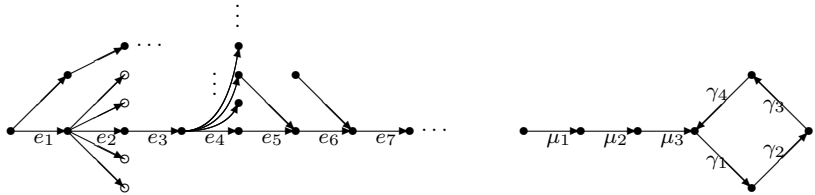
\begin{center}
\setlength{\unitlength}{1.5mm}
\begin{picture}(70,20)
	\put(2,10){\circle*{0.7}}
	\put(7,10){\circle*{0.7}}
	\put(12,10){\circle*{0.7}}		\put(7,15){\circle*{0.7}}
	\put(12,17.5){\circle*{0.7}}
	\put(12,15){\circle{0.7}}
	\put(12,12.5){\circle{0.7}}
	\put(12,5){\circle{0.7}}
	\put(12,7.5){\circle{0.7}}
	\put(17,10){\circle*{0.7}}
	\put(22,10){\circle*{0.7}}
	\put(27,10){\circle*{0.7}}
	\put(32,10){\circle*{0.7}}		\put(27,15){\circle*{0.7}}
	\put(22,15){\circle*{0.7}}
	\put(22,17.5){\circle*{0.7}}
	\put(22,12.5){\circle*{0.7}}		\put(37,10){\circle*{0.7}}
	
	\put(2,10){\vector(1,0){5}}
	\put(2,10){\vector(1,1){5}}		\put(7,15){\vector(2,1){5}}
	\put(7,10){\vector(1,0){5}}
	\put(7,10){\vector(1,1){5}}
	\put(7,10){\vector(2,1){5}}
	\put(7,10){\vector(1,-1){5}}
	\put(7,10){\vector(2,-1){5}}
	\put(12,10){\vector(1,0){5}}
	\put(17,10){\vector(1,0){5}}
	\put(22,10){\vector(1,0){5}}
	\put(27,10){\vector(1,0){5}}		\put(22,15){\vector(1,-1){5}}
	\put(27,15){\vector(1,-1){5}}
	\put(32,10){\vector(1,0){5}}

	\put(4,9){\scriptsize$e_1$}
	\put(9.5,9){\scriptsize$e_2$}
	\put(14,9){\scriptsize$e_3$}
	\put(19,9){\scriptsize$e_4$}
	\put(24,9){\scriptsize$e_5$}
	\put(29,9){\scriptsize$e_6$}
	\put(34,9){\scriptsize$e_7$}
	\put(38,10){\scriptsize$\ldots$}	\put(13,17.5){\scriptsize$\ldots$}

	\qbezier(17,10)(21,10)(22,17.5)
	\put(22.1,17.5){\vector(1,4){0}}
	\qbezier(17,10)(21,10)(22,15)
	\put(22,15){\vector(1,4){0}}	\qbezier(17,10)(21,10)(22,12.5)
	\put(22,12.5){\vector(1,2){0}}
	\begin{rotate}{90}
	\put(13,-20){\scriptsize$\ldots$}
	\put(19,-22){\scriptsize$\ldots$}
	\end{rotate}
	
\put(49,9){\scriptsize$\mu_1$}
\put(54,9){\scriptsize$\mu_2$}
\put(59,9){\scriptsize$\mu_3$}
\put(47,10){\circle*{0.7}}
\put(52,10){\circle*{0.7}}
\put(57,10){\circle*{0.7}}
\put(62,10){\circle*{0.7}}
\put(67,5){\circle*{0.7}}
\put(67,15){\circle*{0.7}}		
\put(72,10){\circle*{0.7}}
	
	\put(47,10){\vector(1,0){5}}	
	\put(52,10){\vector(1,0){5}}	
	\put(57,10){\vector(1,0){5}}	
	\put(62,10){\vector(1,-1){5}}	
	\put(67,5){\vector(1,1){5}}	
	\put(72,10){\vector(-1,1){5}}	
	\put(67,15){\vector(-1,-1){5}}	
	
	\put(63,7){\scriptsize$\gamma_1$}
	\put(68,8){\scriptsize$\gamma_2$}
	\put(68,12){\scriptsize$\gamma_3$}
	\put(63,13){\scriptsize$\gamma_4$}	
	\end{picture}
\end{center}
\vspace*{-1cm}

\end{exemplo}

Notice that an ultragraph $\G$ satisfies Condition~(L) if, and only if, there are no eventually periodic points that are isolated in $X$. Following
\cite{Ark_Eil_Rui}, we define non-wandering points and thus characterize all kinds of isolated points in an ultragraph shift space.

\begin{defi}\label{nomade}
A path $x\in \p^\infty$ is a \emph{non-wandering point} if $\varepsilon(r(x_i))=\{x_{i+1}\}$ and also $|r(x_i)|=1$, for all $i$. An \emph{eventually non-wandering point} is an element $x=\mu\gamma\in p^{\infty}$, where $\mu\in \p$ and $\gamma\in \p^{\infty}$ is non-wandering.
\end{defi}
 
Note that a path $\gamma^\infty$, where $\gamma$ is a cycle with no exit, is non-wandering. From this we have that eventually periodic isolated points are eventually non-wandering. Also note that there are points that are eventually non-wandering but are not periodic. Finally, directly from Definition~\ref{nomade} and Proposition~\ref{isolados}, we concluded that:

\begin{corolario}
The isolated points of $X$ are exactly the points that are eventually sinks (i.e whose range is a unitary set with a sink) or eventually non-wandering.
\end{corolario}


\subsection*{Continuous Orbit Equivalence} \label{sec3.3.2}

We already have almost all the ingredients needed to describe continuous orbit equivalence for ultragraph shift spaces. But, since we are applying general results for Deaconu-Renault systems and we would like the paper to be as self contained as possible, we will recall the relevant concepts of \cite{DeaRen} along the development of the ultragraph orbit equivalence theory.

\begin{defi} 
A \emph{Deaconu-Renault system} is a pair $(X,\sigma)$ consisting of a locally compact Hausdorff space $X$, and a local homeomorphism $\sigma: Dom(\sigma)\longrightarrow Im(\sigma)$ from an open set $Dom (\sigma)\subseteq X$ to an open set $Im(\sigma)\subseteq X$. Inductively define $Dom(\sigma^n):=\sigma^{-1}(\sigma^{n-1})$, so each $\sigma^n:Dom(\sigma^n)\longrightarrow Im(\sigma^n)$ is a local homeomorphism and $\sigma^m\circ \sigma^n=\sigma^{m+n}$ on $Dom(\sigma^{m+n})$.
\end{defi}

\begin{defi}\label{grup}
The \emph{Deaconu-Renault groupoid} $(X,\sigma)$ is $$G=G(X,\sigma)=\displaystyle\bigcup_{n,m\in\N} \{(x,n-m,y)\in X_\infty^{\geq n}\times \{n-m\}\times X_\infty^{\geq m}: \sigma^n(x)=\sigma^m(y)\},$$
under the topology with basic open sets $Z(U,n,m,V):=\{(x,n-m,y):x\in U,y\in V,\text{ and }\sigma^n(x)=\sigma^m(y)\}$ indexed by quadruples $(U,n,m,V)$, where $n,m\in \N$, $U\subseteq X_\infty^{\geq n}$ and $V\subseteq X_\infty^{\geq m}$ are open and $\sigma^n|_U$ and $\sigma^m|_V$ are homeomorphism.
$G$ is a locally compact, Hausdorff, amenable, étale groupoid, with unit space $\{(x,0,x):x\in X\}$ identified with $X$.
\end{defi}

Thus, given an ultragraph that satisfies Condition~(RFUM2), the shift space $(X,\sigma)$ associated to $\G$ is a Deaconu-Renault system. Furthermore, the Deaconu-Renault groupoid associated to this system is the same as the groupoid $G(X,\sigma)$ we associated with an ultragraph on page \pageref{grupoide}. 

Next we recall from \cite{DeaRen} the notion of stabilizers. Recall that when we write $\sigma^n(x)$ we are implicitly assuming that $x\in X_{\infty}^{\geq n}$ (this last set can be seen in the ultragraph system context or in the more general Deaconu-Renault systems context).

\begin{defi}
Let $x\in X^{\geq 1}_\infty$. The \emph{stabiliser group} at $x$ is defined by 
$$stab(x):=\{m-n:m,n\in\N,x\in X_\infty^{\geq m}\cap X_\infty^{\geq n}\text{ and }\sigma^m(x)=\sigma^n(x)\}\subseteq\Z.$$
The \emph{essential stabiliser group} of $x$ is defined by
$$stab^{ess}(x)=\{m-n:m,n\in\N\text{ and }\exists \ U \text{ neighborhood of }x;U\subseteq X_\infty^{\geq m}\cap X_\infty^{\geq n},\text{ such that }\sigma^m|_U=\sigma^n|_U\}.$$
\end{defi}
\begin{obs}
Note that, for each $x\in X^{\geq 1}_\infty$, we have $stab^{ess}(x)\subseteq stab(x)$.
\end{obs}
\begin{defi}
With the convention that $min\ \emptyset=\infty$, we define the \emph{minimal stabiliser} of $x\in X^{\geq 1}_\infty$ to be
$$stab_{min}(x)=min\{n:n\in stab(x);n\geq1\},$$
and the \emph{minimal essential stabiliser} to be
$$stab_{min}^{ess}(x)=min\{n:n\in stab^{ess}(x);n\geq1\}.$$
\end{defi}
\begin{obs}
	Note that for each $x\in X^{\geq 1}_\infty$, we have $stab_{min}(x)\leq stab_{min}^{ess}(x)$.
\end{obs}
\begin{exemplo}
Consider the following ultragraph:
\begin{center}
\vspace*{-1cm}
\setlength{\unitlength}{1.5mm}
\begin{picture}(50,20)
\put(9,9){\scriptsize$a_1$}
\put(14,9){\scriptsize$a_2$}
\put(19,9){\scriptsize$a_3$}
\put(24,9){\scriptsize$b$}
\put(29,9){\scriptsize$c$}
\put(30,11){\scriptsize$c$}
\put(34,9){\scriptsize$d$}
\put(39,9){\scriptsize$e$}
\put(46,13){\scriptsize$f$}

\put(7,10){\circle*{0.7}}
\put(12,10){\circle*{0.7}}
\put(17,10){\circle*{0.7}}
\put(22,10){\circle*{0.7}}
\put(27,10){\circle*{0.7}}		
\put(32,10){\circle*{0.7}}
\put(37,10){\circle*{0.7}}
\put(42,10){\circle*{0.7}}

\put(7,10){\vector(1,0){5}}	
\put(12,10){\vector(1,0){5}}	
\put(17,10){\vector(1,0){5}}	
\put(22,10){\vector(1,0){5}}	
\put(27,10){\vector(1,0){5}}	
\put(32,10){\vector(1,0){5}}	
\put(37,10){\vector(1,0){5}}

\qbezier(27,10)(35,15)(27,15)
\qbezier(27,15)(23,15)(22,10)
\qbezier(42,10)(50,15)(42,15)
\qbezier(42,15)(37,15)(32,10)
\put(22,10){\vector(-1,-3){0}}	
\put(32,10){\vector(-1,-1){0}}	
\end{picture}
\vspace*{-1.5cm}
\end{center}

An element that belongs to the ultrapath shift space is $x=a_1a_2a_3bcbcbcbcbcbc\ldots.$ Such point is eventually periodic and $stab(x)=\{5-3, 7-5,7-3,9-7,9-5,9-3,5-5,7-7,3-7,3-7,5-7,\ldots\}=2\Z$. Thus, $stab_{min}(x)=2$. But note that $stab^{ess}(x)=\{0\}$, because there is no $U\subseteq X_\infty^{\geq m}\cap X_\infty^{\geq n}=X_\infty^{\max\{m,n\}}$ such that $x\in U$ and $\sigma^m|_U=\sigma^n|_U$. Indeed, given $m,n$, $m\neq n$, for which $\sigma^m(x)=\sigma^n(x)$, if there is such $U$
there must be a cylinder $D_\alpha$ for which $x\in D_\alpha\subseteq U$. Thus $\alpha=a_1a_2a_3bcbcbc\ldots bc$. But note that if we take $z=a_1a_2a_3\underbrace{bc\ldots bc}_{|\alpha|-3+m+n} defdef\ldots\in D_\alpha\subseteq U$ then $\sigma^m(z)\neq \sigma^n(z)$, since $m\neq n$. So $stab^{ess}(x)=\{0\}$ and hence $stab^{ess}_{min}(x)=\infty$. Also note that $bc$ is a cycle that has an exit. 

Another element in the ultragraph shift space is $y=a_1a_2a_3bcdefdefdef\ldots$, for which we have $stab(y)=stab^{ess}(y)=3\Z$. The cycle of $f$ has no exit.
\end{exemplo}

In the next four results we will characterize periodicity in terms of the stabilizers. 

\begin{lema}
	\label{lema_64}
	If $x\in X$ is an eventually periodic point, then there exists a simple cycle $\gamma\in \p^n$, with $0<n<\infty$, such that $stab(x)=|\gamma|\Z$. Furthermore, $n$ is the length of the simple cycle of $x$.  
\end{lema}
\begin{proof}
	Indeed, since $x$ is eventually periodic, we have that $x=\mu_1\mu_2\ldots\mu_k\gamma_1\gamma_2\ldots\gamma_n\gamma_1\gamma_2\ldots\gamma_n\ldots.$
	So let $\gamma=\gamma_1\gamma_2\ldots\gamma_n$, (assuming that $\gamma_i\neq \gamma_j$, if $i\neq j$. Otherwise we rewrite $x$). Then $|\gamma|=n$, i.e $\gamma\in \p^n$, and  $stab(x)=n\Z=|\gamma|\Z$, where $n$ is the length of the cycle of $x$.  
\end{proof}

\begin{lema}
\label{lema_65}	If $x\in X_{fin}$ then $stab(x)=\{0\}$ and $stab_{min}(x)=\infty$. Consequently $stab^{ess}(x)=\{0\}$ and $stab_{min}^{ess}(x)=\infty$.
\end{lema}
\begin{proof}
	Indeed, if $x=(x_1x_2\ldots x_{|x|},r(x))$ then for every $0\leq m,n \leq |x|$, where $m\neq n$, we have $\sigma^m(x)\neq \sigma^n(x)$. 
\end{proof}
\begin{prop}
\label{prop_eve_per}
Let $x\in X$. Then $stab_{min}(x)<\infty$ if, and only if, $x\in \p^\infty$ and $x$ is eventually periodic.
\end{prop}
\begin{proof}
Suppose that $stab(x)_{min}<\infty$. Then $x\in\p^\infty$, since by~Lemma \ref{lema_65} if $x\in X_{fin}$ then $stab(x)=\{0\}$ and $stab_{min}(x)=\infty$. So there exist $m,n\in \N$, where $m\neq n$, such that $x\in X^{\geq m}\cap X^{\geq n}$ and $\sigma^m(x)=\sigma^n(x)$. Therefore, if $n<m$, we have that $x_{n+1}\ldots=x_{m+1}\ldots $, and hence $$x=x_1\ldots x_nx_{n+1}\ldots x_mx_{m+1}\ldots=x_1\ldots x_nx_{n+1}\ldots x_mx_{n+1}\ldots x_m\ldots=x_1\ldots x_n(x_{n+1}\ldots x_m)^\infty.$$ 
Fixing $\mu=x_1\ldots x_n$ and $\gamma=x_{n+1}\ldots x_m$, we have that $x=\mu\gamma\gamma\ldots$ and so $x$ is eventually periodic.
On the other hand, suppose that $x$ is eventually periodic. Then $x=\mu\gamma\gamma\ldots$ and, by Lemma~\ref{lema_64}, there exist $\gamma\in \p^n$, with $n>0$, such that $stab(x)=|\gamma|\Z=n\Z$. So $stab_{min}(x)=n<\infty$.
\end{proof}

We have not yet defined the notion of continuous orbit equivalence. We do this below. Notice that although the definition  is stated in the context of Deaconu-Renault systems (as in \cite{DeaRen}) it readly translated to the context of ultragraph shift spaces (which are examples of Deaconu-Renault systems). 

\begin{defi}\label{def_orb_equi}
Let $(X,\sigma)$ and $(Y,\tau)$ be Deaconu-Renault systems. We say that $(X,\sigma)$ and $(Y,\tau)$ are \emph{continuous orbit equivalent} if there exist a homeomorphism $h:X\rightarrow Y$ and continuous maps $k,l:Dom(\sigma)\rightarrow \N$, $k',l':Dom(\tau)\rightarrow \N$ such that
\begin{align*}
\tau^{l(x)}(h(x))&=\tau^{k(x)}(h(\sigma(x)))\\
\sigma^{l'(y)}(h^{-1}(y))&=\sigma^{k'(y)}(h^{-1}(\tau(y)),
\end{align*} for all $x,y$. 
We call $(h,l,k,l',k')$ a continuous orbit equivalence and we call  $h$ the \emph{underlying homeomorphism}.
We say that $(h,l,k,l',k')$ \emph{preserves stabilisers} if $stab_{min}(h(x))<\infty\Leftrightarrow stab_{min}(x)<\infty$ and 
\begin{equation}
\label{equa1}
\left| \displaystyle\sum_{n=0}^{stab_{min}(x)-1} l(\sigma^n(x))-k(\sigma^n(x))\right|=stab_{min}(h(x))
\end{equation}
and
\begin{equation}
\label{equa2}
\left| \displaystyle\sum_{n=0}^{stab_{min}(y)-1} l'(\tau^n(y))-k'(\tau^n(y))\right|=stab_{min}(h^{-1}(y))
\end{equation}
whenever $stab(x)\neq\{0\}$, $stab(y)\neq\{0\}$,  $\sigma^{stab_{min}(x)} (x)=x$ and $\sigma^{stab_{min}(y)} (y)=y$.

Likewise, we say $(h,l,k,l',k')$ \emph{preserves essential stabilisers} if $stab_{min}^{ess}(h(x))<\infty\Leftrightarrow stab_{min}^{ess}(x)<\infty$ and 
\begin{equation}
\label{equa3}
\left| \displaystyle\sum_{n=0}^{stab_{min}^{ess}(x)-1} l(\sigma^n(x))-k(\sigma^n(x))\right|=stab_{min}^{ess}(h(x))
\end{equation} 
and
\begin{equation}
\label{equa4}
\left| \displaystyle\sum_{n=0}^{stab_{min}^{ess}(y)-1} l'(\tau^n(y))-k'(\tau^n(y))\right|=stab_{min}^{ess}(h^{-1}(y))
\end{equation}
whenever $stab^{ess}(x)\neq\{0\}$, $stab^{ess}(y)\neq\{0\}$,  $\sigma^{stab_{min}^{ess}(x)} (x)=x$ and $\sigma^{stab_{min}^{ess}(y)} (y)=y$.
\end{defi}

As we have seen above, stabilizers are related to periodicity. We therefore make a definition regarding continuous orbit equivalence that preserves stabilizers in terms of periodicity.

\begin{defi}
We say that a continuous orbit equivalence $(h,l,k,l',k')$ between Deaconu-Renault systems $(X,\sigma)$ and $(Y,\tau)$ preserves isolated eventually periodic points if $x\in X$ is an isolated eventually periodic point $\iff$ $h(x)$ is an isolated eventually periodic point and, in addition, $y\in Y$ is an isolated eventually periodic point $\iff$ $h^{-1}(y)$ is an isolated eventually periodic point.
\end{defi}

\begin{obs}
In light of Proposition~\ref{prop_eve_per} we have that $stab(x)\neq \{0\}$ and $\sigma^{stab_{min}(x)}(x)=x \iff x$ is periodic. Thus it is clear that preserving stabilisers is the same as preserving eventually periodic points and, furthermore, whenever $x$ and $y$ are periodic the equations (\ref{equa1}) and (\ref{equa2}) are satisfied.
\end{obs}

We can now state and prove the main result of this section, which is a generalization of \cite[Theorem 4.2]{CarlsenLie}. But first we need a last definition.

\begin{defi}\cite[Page 26]{DeaRen}
Let $(X,\sigma)$ be a shift space. Define a continuous cocycle $c_X:G(X,\sigma)\rightarrow \Z$ by $c_X(x,n,y)=n$. 
\end{defi}

\begin{teorema}\label{teo_3}
Let $\G_1,\G_2$ be ultragraphs that satisfy Condition~(RFUM2). Let $(X_1,\sigma)$ and $(X_2,\tau)$ be the corresponding shift spaces. Let $h:X_1\rightarrow X_2$ be a homeomorphism, and $k,l:X_1^{\geq 1}\rightarrow \N$, $k',l':X_2^{\geq 1}\rightarrow \N$ be continuous maps. Then the following are equivalent:
\begin{enumerate}[(1)]
	\item There is a groupoid isomorphism $\theta :G(X_1,\sigma)\rightarrow G(X_2,\tau)$ such that $\theta|_{X_1}=h$, $\theta (x,1,\sigma(x))=(h(x),l(x)-k(x),h(\sigma(x)))$ and 
	$\theta^{-1} (y,1,\tau(y))=(h^{-1}(y),l'(y)-k'(y),h^{-1}(\tau(y)))$;
	
	\item $(h,l,k,l',k')$ is a stabilizer-preserving continuous orbit equivalence;
	
	\item $(h,l,k,l',k')$ is an essential-stabilizer-preserving continuous orbit equivalence; and
	
	\item $(h,l,k,l',k')$ is a continuous orbit equivalence that preserves isolated eventually periodic points and satisfies the equations (\ref{equa1}) and (\ref{equa2}).
\end{enumerate}
\end{teorema}
\begin{proof}
The equivalences $(1)\iff(2)\iff(3)$ are proved in \cite[Prop. 8.3]{DeaRen}. We have to show the equivalence of these statements with $(4)$. First we show $(2)\Rightarrow (4)$. Let $x\in X_1$ be an isolated eventually periodic point. Then, by Lemma~\ref{lema_64} and Proposition~\ref{prop_eve_per}, there exists a $\gamma\in\p$ for which $x=\mu\gamma\gamma\ldots$ and $stab_{min}(x)=|\gamma|<\infty$. Since $h$ preserves stabilisers, we have that $stab_{min}(h(x))<\infty$ and, by Proposition~\ref{prop_eve_per} again, $h(x)$ is eventually periodic. Now, since $h$ is homeomorphism, if $x$ is an isolated point then $h(x)$ is also an isolated point. For $y\in X_2$ and $h^{-1}:X_2\longrightarrow X_1$ the proof is analogous. 

Next we show $(4)\Rightarrow(1)$. Our proof is motivated by the ideas in the proof of \cite[Prop. 8.3]{DeaRen}. By \cite[Lemma 8.8]{DeaRen}, since $(h,l,k,l',k')$ is a continuous orbit equivalence between $(X_1,\sigma)$ and $(X_2,\tau)$, there is a continuous homomorphism of groupoids $\theta_{k,l}:G(X_1,\sigma)\rightarrow G(X_2,\tau)$, defined by $(x,m-n,x')\mapsto (h(x),l_m(x)-k_m(x)-l_n(x')+k_n(x'),h(x'))$ whenever $\sigma^m(x)=\sigma^n(x')$ (notice that $\theta_{k,l}((x,m-n,x'))$  belongs to $G(X_2,\tau)$ since $(h,l,k,l',k')$ is a continuous orbit equivalence of $(X, \sigma)$ to $(Y, \tau)$, and by \cite[Lemma 8.7]{DeaRen} there is a continuous cocycle $c_{(h,l,k,l',k')} : \G(X,\sigma) \to \Z$ such that
$c_{(h,l,k,l',k')}(x,m-n,x')=l_m(x)-k_m(x)-l_n(x')+k_n(x')$). Furthermore, for each $x\in X_1$, there is a group homomorphism $\pi_x: stab(x)\rightarrow stab(h(x))$, where $\pi_x(m-n)=l_m(x)-k_m(x)-l_n(x)+k_n(x)$, whenever 
$\sigma^m(x)=\sigma^n(x)$ (notice that $\pi_x$ is well defined by \cite[Lema 8.8]{DeaRen}). In addition, for every $x\in X_1$, $m,n\in\N$, we have that $stab(\sigma^m(x))=stab(\sigma^n(x))$,  $stab(h(\sigma^m(x)))=stab(h(\sigma^n(x)))$ and $\pi_{\sigma^m(x)}=\pi_{\sigma^n(x)}$. 

It remains to prove that $\theta:=\theta_{k,l}$ is bijective and that $\theta^{-1}$ is continuous. First we show that $\theta$ is injective: suppose that $\theta(x_1,n_1,x_1')=\theta(x_2,n_2,x_2')$. Since $h$ is a homeomorphism, $x_1=x_2$ and $x_1'=x_2'$. Then $\theta(x_1,n_1-n_2,x_1)=\theta(x_1,n_1,x_1')\theta(x_1,n_2,x_1')^{-1}=(h(x_1),0,h(x_1))$, i.e $\pi_{x_1}(n_1-n_2)=0$. Note that if $x_1$ is not eventually periodic, by Proposition~\ref{prop_eve_per}, we have that $stab(x_1)=\{0\}$. Since $h$ preserves eventually periodic points, $h(x_1)$ cannot be eventually periodic. Thus, by the same Proposition, $stab(h(x_1))=\{0\}$. So, in this case, $\pi_{x_1}$ is the null homomorphism (which is bijective). If $x_1$ is eventually periodic, since $h$ preserves such points, we have that $h(x_1)$ is eventually periodic. By Proposition~\ref{prop_eve_per}, $stab(x_1)$ and $stab(h(x_1))$ are not null subgroups of $\Z$. Since, by hypothesis, $(h,l,k,l',k')$ satisfies the equations (\ref{equa1}) and (\ref{equa2}), by \cite[Lemma 8.8]{DeaRen} we have that $stab(.)$, $stab(h(.))$ and $x\mapsto \pi_x$ are constant on orbits and, with this,  $\pi_{x_1}:stab(x_1)\rightarrow stab(h(x_1))$ is a group bijection. Since $\pi_{x_1}(n_1-n_2)=0$ and $\pi_{x_1}$ is bijective (being $x_1$ eventually periodic or not), we have that $n_1-n_2=0$, i.e $n_1=n_2$, which shows the injectivity of $\theta$.

To show the surjetivity of $\theta$, let $(y,n,y')\in G(X_2,\tau)$. Then $\theta(\theta'(y,n,y'))\overset{(*)}{=}(y,m,y')$ for some $m\in \Z$. Hence, $n-m\in stab(y)$, because if $a\in stab(x)$ and $b\in stab(x)$, then $a-b\in stab(x)$. Since $\pi_{h^{-1}(y)}:stab(h^{-1}(y))\rightarrow stab(y)$ is bijective, we have that $n-m\overset{(**)}{=} \pi_{h^{-1}(y)}(p)$ for some $p\in stab(h^{-1}(y))$. Thus, 
\begin{align*}
\theta(h^{-1}(y),p+c_{X_1}(\theta'(y,n,y')),h^{-1}(y'))=&\theta(h^{-1}(y),p,h^{-1}(y))\theta(h^{-1}(y),c_{X_1}(\theta'(y,n,y')),h^{-1}(y'))\\
\overset{(*) (**)}{=}&(y,\pi_{h^{-1}(y)}(p),y)(y,m,y')=(y,n-m+m,y')
=(y,n,y'),
\end{align*} where $c_{X_1}$ is the continuous cocycle of $G(X_1,\sigma)$ in $\Z$ defined by $c_{X_1}(x,n,y)=n$.

It remains to be shown that $\theta^{-1}$ is continuous. Suppose that $(y_n,m_n,y_n')\overset{(***)}{\longrightarrow} (y,m,y')\in G(X_2,\tau)$. Let's show that $\theta^{-1}(y_n,m_n,y_n')\longrightarrow \theta^{-1}(y,m,y')$. Fix $p,q\in\N$, and open $U,V\subseteq X_1$ such that $h^{-1}(y)\in U$, $h^{-1}(y')\in V$, $\sigma^p|_U=\sigma^q|_V$ are homeomorphism such that $\sigma^p(h^{-1}(y))=\sigma^q(h^{-1}(y'))$ (these exist because $(y,m,y')\overset{\theta^{-1}}{\longmapsto} (h^{-1}(y),``\ldots",h^{-1}(y'))\in G(X_1,\tau)$ and the groupoid is étale), and  $\theta^{-1}(y,m,y')=(h^{-1}(y),p-q,h^{-1}(y'))$. Thus $\theta^{-1}(y,m,y')\in Z(U,p,q,V)$ (these cylinders appeared first in the proof of Theorem~\ref{teo_grupoides}). Now, let $U',V'$ be neighborhoods of $U,V$, with $h^{-1}(y)\in U'$, $h^{-1}(y')\in V'$, such that $\theta(x,p-q,x')=(h(x),m,h(x'))$ whenever $x\in U'$, $y\in V'$ and $\sigma^p(x)=\sigma^q(x')$. Fix $N$ such that $y_n\in h(U')$, $y_n'\in h(V')$ and $m_n=m$ for all $n\geq N$ ($N$ exists by convergence $(***)$ and since $h$ is a homeomorphism). Then $\theta^{-1}(y_n,m_n,y_n')\in Z(U',p,q,V')\subseteq Z(U,p,q,V)$ for all $n>N$, i.e $\theta^{-1}(y_n,m_n,y_n')\longrightarrow \theta^{-1}(y,m,y')$.
\end{proof}

We have the following immediate consequence of the theorem above, see \cite[Theorem 8.2]{DeaRen}.

\begin{corolario}\label{corolario3.36}
Let $\G_1$ and $\G_2$ be ultragraphs that satisfy Condition~(RFUM2). Let $(X_1,\sigma)$ and $(X_2,\tau)$ be the corresponding shift spaces. Consider $h:X_1\rightarrow X_2$ a homeomorphism. Then the following are equivalent:
\begin{enumerate}[(1)]
	
	\item There is a stabiliser-preserving continuous orbit equivalence from $(X_1,\sigma)$ to $(X_2,\tau)$ with underlying homeomorphism $h$;
	
	\item There is a groupoid isomorphism $\theta :G(X_1,\sigma)\rightarrow G(X_2,\tau)$ such that $\theta|_{X_1}=h$; 
	
	\item There is an isomorphism $\phi:C^*(G(X_1,\sigma))\longrightarrow C^*(G(X_2,\tau))$ such that $\phi(C_0(X_1))=C_0(X_2)$, with $\phi(f)=f\circ h^{-1}$ for all $f\in C_0(X_1)$;
	
	\item There is an isomorphism $\varphi:C_0(X_1)\rtimes_\alpha\F\longrightarrow C_0(X_2)\rtimes_\alpha\F$ such that $\varphi(C_0(X_1))=C_0(X_2)$, with $\varphi(f)=f\circ h^{-1}$ for all $f\in C_0(X_1)$;
	
	\item There is a continuous orbit equivalence from $(X_1,\sigma)$ to $(X_2,\tau)$ that preserves isolated eventually periodic points and that satisfies the equations (\ref{equa1}) and (\ref{equa2}), where $h$ is underlying homeomorphism.
\end{enumerate}
\end{corolario}
\begin{proof}
The equivalence between (1), (2), and (3) follows from \cite[Theorem 8.2]{DeaRen}. Item $(5)$ is equivalent to item $(1)$ by Theorem~\ref{teo_3}. Finally, item $(4)$ is equivalent to item $(3)$ by \cite{Abadie} and Theorem~\ref{teo_grupoides}.
\end{proof}

\subsection*{Eventual Conjugacy}
We finish this section describing eventual conjugacy between ultragraph shift spaces. As with continuous orbit equivalence, we will rely on general results of \cite{DeaRen}.

\begin{defi}
Let $(X,\sigma)$ and $(Y,\tau)$ be shift spaces associated to ultragraphs that satisfies Condition~(RFUM2). We say that $(X,\sigma)$ and $(Y,\tau)$ are eventually conjugate if there is a stabiliser-preserving continuous orbit equivalence $(h,k,l,k',l')$ from $(X,\sigma)$ to $(Y,\tau)$ such that $l(x)=k(x)+1$, for all $x\in X$.
\end{defi} 

In addition notice that for any Deaconu-Renault system $(X,\sigma)$ (and in particular for ultragraph shift spaces) there is an action $\gamma^X:\T\rightarrow Aut\left(C^*(G(X,\sigma))\right)$ that satisfies $\gamma^X_z(f)(x,n,x')=z^nf(x,n,x')$, for all $z\in\T$, $(x,n,x')\in G(X,\sigma)$, and $f\in C_c(G(X,\sigma))$. We can now state the last result of this section, which is a  generalization of \cite[Theorem 4.1] {Car_Jam}.

\begin{teorema}
Let $(X,\sigma)$ and $(Y,\tau)$ be ultragraph shift spaces associated to ultragraphs that satisfy Condition~(RFUM2). The following are equivalent:
\begin{enumerate}
\item There is an eventual conjugacy from $(X,\sigma)$ to $(Y,\tau)$ with underlying homeomorphism $h$;
\item There is an isomorphism $\theta:G(X,\sigma)\longrightarrow G(Y,\tau)$ such that $\theta|_X=h$ and $c_X=c_Y\circ\theta$; 
\item There is an isomorphism $\phi:C^*(G(X,\sigma))\longrightarrow C^*(G(Y,\tau))$ such that $\phi(C_0(X))=C_0(Y)$, with $\phi(f)=f\circ h^{-1}$ for all $f\in C_0(X)$ and $\phi\circ \gamma_z^X=\gamma_z^X\circ \phi$ for all $z\in\T$.
\item There is an isomorphism $\varphi:C_0(X)\rtimes_\alpha\F\longrightarrow C_0(Y)\rtimes_\alpha\F$ such that $\varphi(C_0(X))=C_0(Y)$, with $\varphi(f)=f\circ h^{-1}$ for all $f\in C_0(X)$ and $\varphi\circ \gamma_z^X=\gamma_z^X\circ \varphi$ for all $z\in\T$.
\end{enumerate}
\end{teorema}
\begin{proof}
The proof follows from \cite[Theorem 8.10]{DeaRen} and Corollary~\ref{corolario3.36}.
\end{proof}


\section{KMS states in C*(\texorpdfstring{$\G$}))}\label{KMS}

KMS states are the focus of intense research in both Mathematics and Physics. Among other things, there is significant evidence that the KMS data is a useful invariant of a dynamical system (although KMS states do not behave well under Morita equivalence). Furthermore, it is known that the addition of sinks in graphs results in new KMS states, see for example \cite{GM}. 

Given the exposed above, it is interesting to characterize KMS states of ultragraphs that satisfy Condition~(RFUM2). In particular, this characterization provides for an invariant of ultragraph shift spaces, and for a way to construct KMS states even when dealing with ultragraphs with sinks. The results we present generalize those in \cite{danielgilles} and in \cite{CarLar}.

We start recalling the definition of KMS (and ground) states associated with a certain one-parameter group of automorphisms of $C^*(\G)$.

\begin{defi}
	Given a C*-algebra $A$ and a homomorphism (a dynamics) $\sigma: \R\longrightarrow Aut(A)$, an element $a\in A$ is called \emph{analytic} provided that $t\longmapsto \sigma_t(a)$ extends to an entire function $z\longmapsto \sigma_z(a)$ on $\C$.
\end{defi}
\begin{obs}
	The analytical elements of a C*-algebra $A$ form a dense subset in $A$ \cite[Cap. 8.12]{Padersen}.
\end{obs}
\begin{defi}\
	For $\beta\in (0,\infty)$, a \emph{KMS$_\beta$-state} of $(A,\sigma)$ is a state $\psi$ of $A$ which satisfies the KMS$_\beta$ condition:
	$$\psi(ab)=\psi(b\sigma_{i\beta}(a))$$ for all $a,b\in A$ analytic.
\end{defi}

It is well known from \cite[Section 5.3.7]{bratteli2012operator} that it is sufficient to show that the KMS condition is satisfied for a subset of analytical elements whose span is dense in $A$.

\begin{defi}
	A state $\phi$ of $A$ is a \emph{ground state} of $(A,\sigma)$ if for every $a,b\in A$ analytic, the entire function $z\longmapsto \phi(a\sigma_z(b))$ is bounded in the upper-half plane, i.e $$\underset{Im z>0}{sup}|\phi(b\sigma_z(a))|<\infty$$ for every $a,b\in A$ analytic. 
\end{defi}

Again, for $\phi$ to be a ground state it is sufficient to have boundedness for a set of elements that spans a dense subalgebra of $A$.

In \cite[Lemma 3.1]{danielgilles} a characterization of KMS states in terms of states of a core subalgebra of $C^*(\G)$ was described. We quickly recall this characterization below. First we need to construct a one parameter group of automorphisms of $C^*(\G)$. For this, let $N: \G^1\longrightarrow \R_+^*$ be a positive function for which there is a constant $k$  such that $N(e)>k$, for every $e\in \G^1$. We extend $N$ to the set of finite paths, say $N:\G^*\longrightarrow \R_+^*$, by making $N(A)=1$ for every $A\in\G^0$, and $N(e_1e_2\ldots e_m)=N(e_1)N(e_2)\ldots N(e_m)$.

\begin{prop}\label{prop_87}
	Let $\G$ be an ultragraph and $N:\G^*\longrightarrow \R_+^*$ as above. There is a strongly continuous action $\rho^C:\R\longrightarrow Aut(C^*(\G))$ such that $\rho_t^C(p_A)=p_A$, for every $A\in\G^0$, and $\rho_t^C(s_e)=N(e)^{it}s_e$, for every $e\in\G^1$. 	
\end{prop}
\begin{proof}
	The demonstration is exactly the same as the one in \cite[Lemma 3.1]{danielgilles}.
\end{proof}

Thus, from a positive function $N$, we get a strongly continuous action that is a one-parameter group of automorphisms of $C^*(\G)$. Now recall from \cite[Remark 2.10]{tom2003} that  $$C^*(\G)=\overline{span}\{s_\mu p_A s_\nu^* :\mu, \nu\in \G^*, A\in\G^0\}.$$ Define the \emph{core subalgebra} $C(\G)^\gamma$ as the fixed point subalgebra, i.e $$C^*(\G)^\gamma=\{a\in C^*(\G);\gamma_z(a)=a,\forall z\in\mathbb{T}\},$$
where $\gamma$ is the gauge action, that is, \label{gaugeaction} $\gamma:\mathbb{T}\longrightarrow Aut(C^*(\G))$ is defined by $\gamma_z(s_\mu)=zs_\mu$ and $\gamma_z(p_A)=p_A$, for every $z\in\mathbb{T},\ \mu\in \G^*\backslash \G^0, \ A\in G^0$. Finally, recall (see \cite{danielgilles}) that there is a conditional expectation $\Psi:C^*(\G)\longrightarrow C^*(\G)^\gamma$, where $$ \Psi(s_\mu p_A s_\nu^*)=\delta_{|\mu|,|\nu|}s_\mu p_A s_\nu^*,$$for all $\mu,\nu\in\G^*,A\in\G^0$. With this set up we have the following result.

\begin{prop}\label{prop_4.9}\cite[Prop. 3.5]{danielgilles}
	Let $N:\G^1\longrightarrow \R_+^*$ be such that $N(\mu)\neq 1$ for all $\mu\in\G^*\backslash\G^0$ and let $\rho$ be the action of Proposition~\ref{prop_87}. Suppose that $\beta\in \R$ and that $\phi$ is a KMS$_\beta$-state on $C^*(\G)$. Then the restriction $\psi:=\phi|_{C^*(\G)^\gamma}$ satisfies:
	$$\psi(s_\mu p_A s_\nu^*)=\delta_{\mu,\nu}N(\mu)^{-\beta}\psi(p_{A\cap r(\mu)}).$$
	On the other hand, given a KMS state $\psi$ on $C^*(\G)^\gamma$ that satisfies the above equality, we have that $\phi=\psi\circ\Psi$ is a KMS$_\beta$-state on $C^*(\G)$, where $\Psi$ is the conditional expectation defined above. Furthermore, the correspondence obtained is an affine bijection.
\end{prop} 


\subsection*{KMS states in $C^*(\G)$ realized as $C_0(X)\rtimes_\alpha \F$}

In this subsection we describe KMS states associated to a one-parameter group of automorphisms, using the isomorphism between $C^*(\G)$ and $C_0(X)\rtimes_\alpha \F$ described in Section~\ref{ultracross}. For this we will build on ideas presented in \cite{CarLar}, \cite{danielgilles}, and \cite{exel_laca}.

From \cite[Theo. 4.3]{exel_laca}, given any function $N:\G^1\longrightarrow (1,\infty)$, there exists a unique strongly continuous one-parameter group $\sigma$ of automorphisms of $C_0(X)\rtimes \F$ such that:
\begin{equation}\label{eq_4.1}
\sigma_t(b)=(N(e))^{it}b\text{ and }\sigma_t(c)=c,
\end{equation}
for every $e\in\G^1$, every $b\in C(X_e)\delta_e$, and every $c\in C_0(X)\delta_0$.

Note that if $N(e)=exp(1)$ for every $e\in \G^1$, then $\sigma_t$ is $2\pi$-periodic and thus induce a strongly continuous action $\beta: \mathbb{T}\longrightarrow Aut(C_0(X)\rtimes \F)$ such that $\beta_z(1_e\delta_e)=z1_e\delta_e\text{ and }\beta_z(f\delta_0)=f\delta_0,$
for every $z\in \mathbb{T},\ e\in \G^1$, and $f\in C_0(X)$. Given this, we have the following result.

\begin{prop}
	Let $\G$ be an ultragraph that satisfies Condition~(RFUM2). Let $\gamma$ be the Gauge action on $C^*(\G)$ defined on the page \pageref{gaugeaction}, let $\beta$ be the action considered above and $\Phi:C^*(\G)\longrightarrow C_0(X)\rtimes_\alpha \F$ be the isomorphism of Theorem \ref{teo_1}. Then, for all $z\in\mathbb{T}$, we have
	$$\Phi\circ\gamma_z=\beta_z\circ \Phi.$$
  
\end{prop}
\begin{proof}
	Let $z\in \mathbb{\T}$. Since $C^*(\G)$ is generated by $\{s_e:e\in\G^1\}\cup\{p_A:A\in\G^0\}$, it is enough to show that equality occurs for these elements. For $A\in \G^0$, we have that $\Phi\circ\gamma_z(p_A)=\Phi(p_A)=1_A\delta_0=\beta_z(1_A\delta_0)=\beta_z\circ \Phi(p_A)$. For $e\in\G^1$, we have that $\Phi\circ\gamma_z(s_e)=\Phi(zs_e)=z1_e\delta_e=\beta_z(1_e\delta_e)=\beta_z\circ \Phi(s_e)$. 
\end{proof}

Our next step is to identify certain states on the fixed algebra core subalgebra $C^*(\G)^\gamma$ in terms of certain states in $C_0(X)$. For this it is important to identify $C_0(X)$ inside $C^*(\G)$.

\begin{obs}
	\label{obs81}
	A useful tool to work with $C_0(X)$ is to work with the dense subalgebra $D$ generated by all characteristic functions $1_c$, $1_A$ and $\alpha_c(1_{c^{-1}} 1_A)$ (see Lemma~\ref{lema3}). Notice that for $c\in \displaystyle\bigcup_{n=1}^{\infty} (\G^1)^n$ and $A\in \G^0$, we have that $\alpha_c(1_{c^{-1}}.1_{r(c)})=1_c$ and $\alpha_c(1_{c^{-1}}.1_A)=1_{r(c)\cap A}\circ \theta_{c^{-1}}$.  Thus, $D$ is generated by characteristic functions $1_A$ and $1_{r(c)\cap A}\circ \theta_{c^{-1}}$, where $A\in \G^0$ and $c\in \G^*\backslash\G^0$.
\end{obs}

\begin{lema}\label{lema_4.12}
	Let $\G$ be an ultragraph that satisfy Condition~(RFUM2). Then $C_0(X)$ is $*$-isomorphic to the C*-subalgebra of $C^*(\G)$ generated by $$\{s_cp_As_c^*|A\in\G^0,c\in \G^*\backslash \G^0\}\cup\{p_A|A\in\G^0\}=\{s_cp_As_c^*|A\in\G^0,c\in \G^*\},$$
	where in the last set we are assuming that if $c=B\in\G^0$ then $s_c:=p_B$.
\end{lema}
\begin{proof}
	Clearly $C_0(X)$ is identified with $C_0(X)\delta_0$ in $C_0(X)\rtimes_\alpha \F$. Let D be the subalgebra of Remark~\ref{obs81}, and $\phi$ be the isomorphism of Theorem~\ref{teo_1}. Notice that $\Phi^{-1}$ takes $1_A\delta_0$ to $p_A$, $1_c\delta_c$ to $s_c$, and $1_{c^{-1}}\delta_{c^{-1}}$ to $s_c^*$. Restricting $\Phi^{-1}$ to the elements $1_A\delta_0$ and $\alpha_c(1_{c^{-1}1_A})\delta_0$ we have the $*$-isomorphism desired, since $1_A\delta_0\longmapsto p_A$, and $\alpha_c(1_{c^{-1}1_A})\delta_0=1_c\delta_c1_A\delta_01_{c^{-1}}\delta_{c^{-1}}\longmapsto s_cp_As_c^*$, for all $A\in\G^0$ and $c\in\G^*\backslash\G^0$.
\end{proof}

\begin{prop}\label{prop_4.13}
Let $\G$ be an ultragraph that satisfy Condition~(RFUM2). Then there is an affine bijection between the set of states $\omega$ on $C_0(X)$ that satisfy the condition $\omega(f\circ \theta_e^{-1})=N(e)^{-\beta}\omega(f)$, for all $e\in \G^1$ and all $f\in C_0(X_{e^{-1}})$, and the set of states $\psi$ on $C^*(\G)^\gamma$ that satisfy $\psi(s_cp_As_d^*)=\delta_{c,d}N(c)^{-\beta}\psi(p_{A\cap r(c)})$, for all $c,d\in\G^*\backslash\G^0$ and $A\in\G^0$.  
\end{prop}
\begin{proof}

Let $\psi$ be a state on $C^*(\G)^\gamma$ that satisfies the equality $\psi(s_cp_As_d^*)=\delta_{c,d}N(c)^{-\beta}\psi(p_{A\cap r(c)})$, for all $c,d\in\G^*\backslash\G^0$ and $A\in\G^0$. Define the state $\omega$ on $C_0(X)$ by $\omega =\psi\circ \Phi^{-1}$, where $\Phi:C^*(\G)\rightarrow C_0(X)\rtimes \F$ is the isomorphism of Theorem~\ref{teo_1}. A sketch of what we will do is the following:
\vspace*{-20pt}
$$\xymatrix{ C^*(\G) \ar@/^2pc/[rrrr]^\Phi & \supseteq C^*(\G)^\gamma \ar@<2pt>[rr]^{\Phi} \ar[d]^\psi & & \ar@<2pt>[ll]^{\Phi^{-1}} C_0(X)\delta_0 \ar[d]^\omega  \subseteq & C_0(X)\rtimes_\alpha \F \\
  &	\C  & & \C }$$
To prove the required condition of $\omega$, it is enough to check it on a dense subset (we will use the one in Remark~\ref{obs81}). Before we show the condition, note that
$$\Phi(s_ep_As_e^*)=\Phi(s_e)\Phi(p_A)\Phi(s_e^*)=1_e\delta_e1_A\delta_01_{e^{-1}}\delta_{e^{-1}}=\alpha_e(1_{e^{-1}}1_A)\delta_0={1_{r(e)\cap A}\circ \theta_{e^{-1}}}\delta_0.$$ 
Hence, since $C_0(X)\cong C_0(X)\delta_0$, we have 
\begin{equation}\label{eq_4.2}
\Phi^{-1}(1_{r(e)\cap A}\circ \theta_{e^{-1}})=s_ep_As_e^*.
\end{equation}
Let $e\in\G^1$. Notice that if $1_A\in C_0(X_{e^{-1}})$ then $A\subseteq X_{e^{-1}}=X_{r(e)}$. Hence,
\begin{align*}
&\omega(1_A\circ\theta_{e^{-1}})= \omega(1_{A\cap r(e)}\circ\theta_{e^{-1}})=
\psi\circ\Phi^{-1}(1_{A\cap r(e)}\circ\theta_{e^{-1}}) \overset{(\ref{eq_4.2})}{=}\psi(s_ep_As_e^*)=\\
&\overset{(\text{condition of } \psi)}{=}N(e)^{-\beta}\psi(p_{A\cap r(e)})=N(e)^{-\beta}\psi(\Phi^{-1}(1_{A\cap r(e)}))=N(e)^{-\beta}\omega(1_{A\cap r(e)})=N(e)^{-\beta}\omega(1_A).
\end{align*}
Now, to prove the condition the other type of generator, let $c\in\G^*\backslash \G^0$ and $A\in \G^0$.
Suppose that $1_{A\cap r(c)}\circ\theta_{c^{-1}}\in C_0(X_{e^{-1}})$. Since $1_{A\cap r(c)}\circ\theta_{c^{-1}} : X_{e^{-1}}\rightarrow \C$, if $x\in Dom(1_{A\cap r(c)}\circ\theta_{c^{-1}})$ then $x\in X_{e^{-1}}\cap X_c=X_{r(e)}\cap X_c$. Therefore $s(c)\in r(e)$, and hence $r(c)=r(ec)$. Thus,
\begin{align*}
&\omega(1_{A\cap r(c)}\circ\theta_{c^{-1}} \circ\theta_{e^{-1}})= \omega(1_{A\cap r(ec)}\circ\theta_{{(ec)}^{-1}})
=\psi\circ\Phi^{-1}(1_{A\cap r(ec)}\circ\theta_{{(ec)}^{-1}})
=\psi(s_{ec}p_As_{ec}^*)
=&\\
&=N(ec)^{-\beta}\psi(p_{r(ec)\cap A})
=N(e)^{-\beta}N(c)^{-\beta}\psi(p_{r(c)\cap A})
=N(e)^{-\beta}\psi(s_cp_As_c^*)\\
&=N(e)^{-\beta}\psi(\Phi^{-1}(1_{A\cap r(c)}\circ \theta_{c^{-1}}))
=N(e)^{-\beta}\omega(1_{A\cap r(c)}\circ \theta_{c^{-1}}),
\end{align*}
and the desired condition is satisfied.

On the other hand, given a state $\omega$ on $C_0(X)$ that satisfies the condition $\omega(f\circ \theta_e^{-1})=N(e)^{-\beta}\omega(f)$, for all $e\in \G^1$ and all $f\in C_0(X_{e^{-1}})$, consider $\psi=\delta_{c,d} \omega\circ \Phi:C^*(\G)\rightarrow \C$. Let $c=c_1\ldots c_n\in G^*\backslash \G^0$ and $A\in \G^0$. Then if $c\neq d$, we have that $\psi(s_cp_As_d^*)=0=N(c)^{-\beta}.0=N(c)^{-\beta}.\psi(p_{A\cap r(e)})$, and if $c=d$ we have that
\begin{align*}
\psi(s_cp_As_d^*)=\omega(\Phi(s_cp_As_d^*))&=\omega(1_{A\cap r(c)}\circ \theta_{c^{-1}})
=\omega(1_{A\cap r(c)}\circ \theta_{c_n^{-1}}\circ \theta_{c_{n-1}^{-1}}\circ \ldots \circ\theta_{c_1^{-1}})\\
&=\omega([1_{A\cap r(c)}\circ \theta_{c_n^{-1}}\circ \theta_{c_{n-1}^{-1}}\circ \ldots \circ\theta_{c_2^{-1}}] \circ\theta_{c_1^{-1}})\\
&=N(c_1)^{-\beta}\omega([1_{A\cap r(c)}\circ \theta_{c_n^{-1}}\circ \theta_{c_{n-1}^{-1}}\circ \ldots \circ\theta_{c_3^{-1}}] \circ\theta_{c_2^{-1}})=\ldots\\
=N(c_1)^{-\beta}\ldots N(c_n)^{-\beta}\omega(1_{A\cap r(c)})&=N(c)^{-\beta}\omega(1_{A\cap r(c)})=N(c)^{-\beta}\omega(\phi(p_{A\cap r(c)})
=N(c)^{-\beta}\psi(p_{A\cap r(c)}).
\end{align*} 
\end{proof}
\vspace*{-1cm}
\subsection*{A measure on the shift space $X$}

By the Riesz representation theorem we can include a description of KMS states of $C^*(\G)$ in terms of measures, which in turn can be described in terms of a certain set of functions on the generalized vertices of $\G$. We start the section developing the measure theoretical concepts, which require a more concrete description of the generalized vertices in $\G$. 

\begin{obs} We will follow the same outline done in \cite{danielgilles} and \cite{edgar}, and so to define the desired measure we are first going to define a semi-ring generated by the cylinders sets.
\end{obs}

\begin{lema}\label{Lema_77}
Let $\G$ be an ultragraph that satisfy Condition~(RFUM2). Then each $A\in\G^0$ can be written uniquely as $A=\displaystyle\bigcup_{i=1}^n A_i$, where there is an unique $k$ such that $|A_k|<\infty$ and $|\varepsilon (A_k)|<\infty$, $A_j$ is a minimal infinite emitter or a minimal sink, and $A_j\cap A_k=\emptyset$ $\forall j\neq k$. In addition, $A_k$ is made up of regular vertices and (non minimal) sinks .
\end{lema}	
\begin{proof}
The existence follows directly from Condition~(RFUM2) and from the characterization of the generalized vertices given in Proposition~\ref{prop_vert_gen}.

To provide uniqueness, suppose that $A=\displaystyle\bigcup_{i=1}^n A_i$, and $A=\displaystyle\bigcup_{i=1}^m A_i'$, with $|A_k|,|A_k'|<\infty$, $|\varepsilon(A_k)|, |\varepsilon(A_k')|<\infty$, and $A_j\cap A_k=\emptyset$, $A_j'\cap A_k'=\emptyset$ whenever $j\neq k$.

If $A_i\in r(e)\cap A_\infty$, i.e, if $A_i$ is a minimal infinite emitter ($i\neq k$), by Lemma~\ref{Lema_1}, we have that $|A_i|=1$ or $|A_i|=\infty$:
\begin{itemize}
\item If $|A_i|=1$, then $\exists $ $A_j'$ such that $A_i\subseteq A_j'$ for some $j\neq k$. Hence $A_j'$ is also a minimal infinite emitter. Again by Lemma~\ref{Lema_1}, we have that $|A_j'|=1$ or $|A_j'|=\infty$. By the minimality of $A_j'$, and since $A_i$ (which is contained in $A_j'$) has cardinality equal to 1, we cannot have $|A_j'|=\infty$. Therefore $|A_j'|=1$ and hence $A_i=A_j'$.
 
\item If $|A_i|=\infty$, then $\exists$ $j\neq k$ such that $|A_i\cap A_j'|=\infty$. Since minimal infinite emitters have no infinite sinks (in $\G^0$), we have that $|\varepsilon(A_i\cap A_j')|=\infty$. Hence,  $A_j'\cap A_i$ is a minimal infinite emitter and so $A_j'$ is a minimal infinite emitter. Since $A_i$ and $A_j'$ are minimal, by Proposition~\ref{obs3} we have that $A_i=A_j'$. 
\end{itemize}

If $A_i\in r(e)\cap A_s$, i.e if $A_i$ is a minimal sink ($i\neq k$), then $|A_i|=\infty$. Hence $\exists$ $j\neq k$ such that $|A_i\cap A_j'|=\infty$. Since minimal sinks have no infinite emitters (in $\G^0$), we also have that $|\varepsilon(A_i\cap A_j')|<\infty$. Therefore, $A_j'\cap A_i$ is a minimal sink and so $A_j'$ is a minimal sink. Since $A_i$ and $A_j'$ are minimal sinks and $|A_i\cap A_j'|=\infty$, by Proposition~\ref{prop1.16} we have that $A_i=A_j'$. Therefore $\displaystyle\bigcup_{i\neq k} A_i=\displaystyle\bigcup_{i\neq k} A_i'$ and $A_k=A_k'$.

The last statement of the lemma is clear, since $A_j$ is a minimal infinite emitter or minimal sink for every $j\neq k$, and so any sink that (eventually) is in $A$, and is not contained in some minimal set, must belong to $A_k$. In addition, as the ultragraph satisfies (RFUM2), $|A_k|<\infty$ and $|\varepsilon (A_k)|<\infty$, so that $A_k$ contains all the other finite regular vertices (which are not contained in some minimal set).
\end{proof}

\begin{obs}
	Note that, for each $(\beta,B)\in\p$, we can identify $D_{(\beta,B)}$ with $D_{(\beta,B),F,S}$, where $F=\emptyset$ and $S=\emptyset$. Thus, every cylinder can be written in the form  $D_{(\beta,B),F,S}$. 
\end{obs}
\begin{corolario}
\label{cor_78}
Let $(\beta,B)\in \p$, $F\subseteq \varepsilon(B)$ finite and $S\subseteq B$ finite ($F$ and $S$ can be the empty set). Then $D_{(\beta,B),F,S}$ can be written as a finite and disjoint union of cylinders from the following collection: $$\{D_{(\beta,A),F',S'}:(\beta,A)\in X_{fin}, A\text{ minimal}, F'\subseteq\varepsilon(A), |F'|<\infty, S\subseteq B,|S'|<\infty,\}$$ union with a finite, disjoint union of elements of the collection of cylinders sets $$\{D_{(\beta,A)}:(\beta,A)\in \p:|A|<\infty,|\varepsilon (A)|<\infty\}.$$
\end{corolario} 
\begin{proof}
By Lemma~\ref{Lema_77}, if $(\beta,B)\in \p$ then $B=\displaystyle\bigcup_{i=1}^n A_i$ and there exists an unique $k$ such that $|A_k|<\infty$ and $|\varepsilon(A_k)|<\infty$. Furthermore, for each $j\neq k$, we have that $A_j$ is a minimal infinite emitter or a minimal sink, and $A_j\cap A_k=\emptyset$. For each $i\neq k$, consider $F_i=\{e\in \G^1:s(e)\in \cup_{j\neq k}(A_i\cap A_j)\}$,  $S_i=\{v\in G^0_s:v\in \cup_{j\neq k}(A_i\cap A_j)\}$, $V=\{s(e):e\in \cup_{i\neq k}F_i,e\notin F\}$ and $W=\{v:v\in \cup_{i\neq k}S_i, s\notin S\}$. Of course, each of these sets is finite. Thus, 
$$D_{(\beta,B),F,S}=\bigsqcup_{i\neq k} D_{(\beta,A_i),F_i\cup F,S_i\cup S}\bigsqcup_{s(e)\in V} D_{(\beta,s(e))} \bigsqcup_{v\in W} D_{(\beta,\{v\})} \bigsqcup_{v\in A_k} D_{(\beta,\{v\})}. $$ 

Note that each $s(e)\in V$ is a regular vertex, otherwise, we get a contradiction with the minimality of $A_i$. Hence,
we have that
$$D_{(\beta,B),F,S}=\bigsqcup_{i\neq k} D_{(\beta,A_i),F_i\cup F,S_i\cup S}\bigsqcup_{v\in V\sqcup W\sqcup A_k} D_{(\beta,\{v\})}.$$ 
Since $\G$ satisfies Condition~(RFUM2), all previous joins are finite.
\end{proof}

In view of the Corollary \ref{cor_78} we  infer the following:

\begin{prop}
The collection of cylinders  $\{D_{(\beta,A),F,S}:(\beta,A)\in X_{fin}, A\text{ minimal}, F\subseteq\varepsilon(A),|F|<\infty,  S\subseteq A,|S|<\infty\}$ union with the collection of cylinders $\{D_{(\beta,A)}:(\beta,A)\in \p:|A|<\infty, \varepsilon (A)<\infty\}$, form a basis for the topology in $X$ defined in Proposition~\ref{prop_base_top}.
\end{prop}
\begin{proof}
This follows directly from Corollary~\ref{cor_78} and the fact that if $(\alpha,A)\in \p$, $|A|<\infty$ and $|\varepsilon(A)|<\infty$, then $D_{(\alpha,A)}=\displaystyle\bigsqcup_{e\in\varepsilon(A)} D_{(\alpha e,r(e))}\bigsqcup_{v\in A\cap G_s^0} D_{(\alpha,\{v\})}$. 
\end{proof}

In light of the above, following the general ideas of \cite[Sec. 5.5]{edgar}, we build a semi-ring.

\begin{prop}\label{prop_4.15}
Let $\G$ be an ultragraph that satisfies Condition~(RFUM2). Then the collection $\widehat{S}$ of all cylinders of the form  $\{D_{(\beta,A),F,S}:(\beta,A)\in X_{fin}, A\text{ minimal}, F\subseteq\varepsilon(A),|F|<\infty,S\subseteq A,|S|<\infty\}$, union with the collection of cylinders of the form  $\{D_{(\beta,A)}:(\beta,A)\in \p:|A|<\infty,|\varepsilon (A)|<\infty\}$, union with the empty set, is a semi-ring. 
\end{prop}
\begin{proof}

Let $C,C_0\subseteq \widehat{S}$ be such that $C_0\subseteq C$. We have to show that there are finite disjoint sets $C_1,C_2,\ldots,C_n\in \widehat{S}$ such that $C\backslash C_0=\displaystyle\bigsqcup_{i=1}^n C_i$.
Suppose first that $C$ is of the form $D_{(\beta,B),F,S}$, for some $(\beta,B)\in X_{fin}$ with $B$ minimal, $F\subseteq \varepsilon(B)$ finite, and $S\subseteq B$ finite. If $C=C_0$ then the result is clear. Suppose then that $C\backslash C_0\neq \emptyset$. Since $C_0\subsetneq C$, we have that $C_0=D_{(\beta\beta',A),F',S'}$ where it can happen that $\beta\beta'=\beta$ with $A\subsetneq B$, $F'=\emptyset$ and $S'=\emptyset$. We divide the proof in some cases:
\begin{itemize}
\item If $C_0=D_{(\beta,A),F',S'}$, where $(\beta,A)\in X_{fin}$ with $A$ minimal, then, as $(\beta,B)\in X_{fin}$ with $B$ minimal, we have $A=B$ and thus $C_0=D_{(\beta,B),F',S'}$. It follows that $F'\supseteq F$ and that $S'\supseteq S$, and hence $C\backslash C_0=\displaystyle\bigsqcup_{e\in F'\backslash F} D_{(\beta e,r(e))}\bigsqcup_{v\in S'\backslash S} D_{(\beta,\{v\})}$. Now note that, since $r(e)\in \G^0$, by Corollary~\ref{cor_78} each $D_{(\beta e,r(e))}$ is a finite and disjoint union of elements of $\widehat{S}$, and the same holds for elements $D_{(\beta,\{v\})}$. Thus, $C\backslash C_0$ is a finite and disjoint union of elements of $\widehat{S}$. 
\item If $C_0=D_{(\beta,A)}$, with $|A|<\infty$ and $|\varepsilon (A)|<\infty$, then $C\backslash C_0= D_{(\beta,B),F\cup\varepsilon(A),S\cup A}\in \widehat{S}$. 
\item If $C_0=D_{(\beta\beta',A),F',S'}$, with $|\beta'|\geq 1$ and $(\beta\beta',A)\in X_{fin}$, suppose that $\beta'=\beta_1'\ldots\beta_n'$. Since $r(\beta')\in \G^0$, by Lemma~\ref{Lema_77}, we have that $r(\beta')=\cup_{i=1}^m A_i$ (unique $A_i$'s), and that there exists an unique $k$ for which $|A_k|<\infty$, $|\varepsilon(A_k)|<\infty$, and also $A_k\cap A_j=\emptyset$ for every $j\neq k$, and $A_j$ is a minimal infinite emitter or a minimal sink for all $j\neq k$. Thus, since $A\subseteq r(\beta')$ and $A$ is minimal, we have that $A=A_i$ for some $i\neq k$. Hence, for all $j\neq i,k$, we have that $|\varepsilon (A_j\cap A_i)|=|\varepsilon(A_j\cap A)|<\infty$ and that $|A_j\cap A_i|=|A_j\cap A|<\infty$. Therefore, we can write: 
\begin{align*}
C\backslash C_0=&D_{(\beta,B),F\cup\{\beta_1'\},S}\sqcup D_{(\beta\beta_1',r(\beta_1')),\{\beta_2'\}}\sqcup\ldots\sqcup D_{(\beta\beta_1'\ldots\beta_{n-1}',r(\beta_{n-1}')),\{\beta_n'\}}\bigsqcup\\
&\bigsqcup_{j\neq i,k} D_{(\beta\beta',A_j),\varepsilon(A_j\cap A_i),(A_j\cap A_i)} \bigsqcup_{e\in F'\cup \varepsilon(A_k)} D_{(\beta\beta'e,r(e))}  \bigsqcup_{v\in S'\cup A_k}  D_{(\beta\beta',\{v\})}.
\end{align*}
By Corollary~\ref{cor_78}, each $D_{(\beta\beta_1'\ldots\beta_p',r(\beta_p'),\{\beta_{p+1}'\})}$, with $1\leq p\leq n-1$, is a finite and disjoint union of elements of $\widehat{S}$. Also, by Corollary~\ref{cor_78}, each $D_{(\beta \beta' e,r(e))}$ is a finite and disjoint union of elements of $\widehat{S}$, just like the elements $D_{(\beta\beta',\{v\})}$. Thus we obtain the desired description of $C\backslash C_0$.
\item If $C_0=D_{(\beta\beta',A)}$, with $|A|<\infty$ and $|\varepsilon(A)|<\infty$, we use the description of $r(\beta')$ of the previous item and then we have that $A\subseteq A_k$. So, we have the following description of $C\backslash C_0$:
\begin{align*}
C\backslash C_0=&D_{(\beta,B),F\cup\{\beta_1'\},S}\sqcup D_{(\beta\beta_1',r(\beta_1')),\{\beta_2'\}}\sqcup\ldots\sqcup D_{(\beta\beta_1'\ldots\beta_{n-1}',r(\beta_{n-1}')),\{\beta_n'\}}\bigsqcup\\
&\bigsqcup_{j\neq k} D_{(\beta\beta',A_j),\varepsilon(A),A} \bigsqcup_{e\in \varepsilon(A_k\backslash A)} D_{(\beta\beta'e,r(e))} \bigsqcup_{v\in A_k\backslash A}  D_{(\beta\beta',\{v\})}.
\end{align*}
that again, by Corollary~\ref{cor_78}, is a finite and disjoint union of elements of $\widehat{S}$. 

For the case where $C$ is of the form $D_{(\beta,A)}$, with $\varepsilon(A)<\infty$, the proof is analogous.
\end{itemize}
\end{proof} 

Before we define the desired measure on the semi-ring $\hat{S}$ we need to introduce a class of functions from $\G^0$ as follows.

Let $M:\G^1\rightarrow[0,1]$ and $m:\G^0\rightarrow [0,1]$ be functions that together satisfy the following conditions:

\begin{enumerate}[\textbf m1.]\label{emes}
\item $\displaystyle\lim_{A\in\G^0} m(A)=1$;

\item $m(A)=\displaystyle\sum_{e:s(e)\in A}M(e)m(r(e))+\sum_{v\in A\cap G_s^0}m(v)$, if $|A|<\infty$ and $|\epsilon(A)|<\infty$,
\\ where $m(v)=m(\{v\})$;

\item $m(A)\geq \displaystyle\sum_{e\in F}M(e)m(r(e))+ \displaystyle\sum_{v\in S}m(v)$, for every $F\subseteq \varepsilon(A)$, $F$ finite and every $S\subseteq A\cap G_s^0$, $S$ finite;

\item $m(A\cup B)=m(A)+m(B)-m(A\cap B)$.
\end{enumerate} 

We extend $M$ to $\G^*$ defining $M(A)=m(A)$, for all $A\in \G^0$, and $M(\beta)=M(e_1)M(e_2)\ldots M(e_n)$, whenever $\beta=e_1e_2\ldots e_n\in \G^*\backslash \G^0$. Thus we define, for each $(\beta,B)\in\p$, $F\subseteq \varepsilon(B)$, $S\subseteq B$, $F$ and $S$ finite, the function $\kappa:\{\text{Cylinders }D_{(\beta,B),F,S}\}\rightarrow \R_+$ by: 
\begin{equation}\label{kappa}
\kappa(D_{(\beta,B),F,S})=M(\beta)m(B)-\displaystyle\sum_{e\in F}M(\beta e)m(r(e))-\displaystyle\sum_{v\in S}M(\beta)m(v)
\end{equation}
if $|\beta|\geq 1$, by $\kappa(D_{(A,A),F,S})=m(A)-\displaystyle\sum_{e\in \varepsilon(F)}M(e)m(r(e))-\displaystyle\sum_{v\in S}m(v)$ if $A\in \G^0$, $F\subseteq \varepsilon(A)$ and $S\subseteq A$, and finally, $\kappa(\emptyset)=0$.

Clearly by m3., $Im(\kappa)\subseteq \R_+$. We can then restrict $\kappa$ to the set $\widehat{S}$ of Proposition~\ref{prop_4.15} to build a measure. So we have:

\begin{lema}\label{lema_medida}
Let $\G$ be an ultragraph that satisfy Condition~(RFUM2) and $\kappa$ the function defined by equation (\ref{kappa}). Then the restriction of $\kappa$ to the semi-ring $\widehat{S}$ of Proposition~\ref{prop_4.15} is a measure such that $\kappa(\theta_e(V))=M(e)\kappa(V)$, for every $e\in \G^1$ and every subset $V\subseteq X_{e^{-1}}\cap S\subseteq X$. 
\end{lema}
\begin{proof}
Since $\kappa(\emptyset)=0$, we need to show that $\kappa$ is countable additive on $\widehat{S}$. But since the elements that generates $\widehat{S}$ are compact and open, it is enough to show that $\kappa$ is additive on $\widehat{S}$, i.e, if $D_{(\beta,B),F,S}=\displaystyle\bigsqcup_{i=1}^n D_{(\beta_i,B_i),F_i,S_i}$, then $\kappa\left(\displaystyle\sqcup_{i=1}^n D_{(\beta_i,B_i),F_i,S_i}\right)=\sum_{i=1}^n\kappa\left(D_{(\beta_i,B_i),F_i,S_i}\right).$
So, suppose that 
\begin{equation}\label{4.3}
D_{(\beta,B),F,S}=\displaystyle\bigsqcup_{i=1}^n D_{(\beta_i,B_i),F_i,S_i}
\end{equation}

We can assume, without loss of generality, that $\beta_i =\beta\beta_i'$, for every $i$. Note that we can have $|\beta_i'|=0$, but this is not necessarily the case. So,we use induction in $m=\max_{1\leq i\leq n}\{{|\beta_i|-|\beta|}\}=\max_{1\leq i\leq n}{|\beta_i'|}$.

Suppose that $m=0$. Then $\beta_i=\beta$, for all $i$, and  $B_i\subseteq B$, for all $i$. Now we analyze each case:
\begin{itemize}
\item If $B$ is such that $|B|<\infty$ and $|\varepsilon(B)|<\infty$, then $B$ is formed only by a finite number of sinks and regular vertices, and so we have $B=\{v_1,\ldots, v_p\}$, $F=S=\emptyset$, and thus $F_i=\emptyset$ and $S_i=\emptyset$ for all $i$. So we can write $D_{(\beta,B),F,S}=\displaystyle\bigsqcup_{i=1}^n D_{(\beta_i,B_i)}$ (note that we could have, for each $i=1,\ldots,n$, $B_i=\{v_{i_1},\ldots,v_{i_k}\}$, where each $B_i$ contains a finite number of regular vertices and sinks $v_i$'s). Thus,
$$\kappa(D_{(\beta.B)})=M(\beta)m(B)=M(\beta)m(\sqcup_{i=1}^n B_i)\overset{\text{by }m4.}{=}M(\beta)\sum_{i=i}^n m(B_i)=\sum_{i=1}^n \kappa(\beta_i,B_i).$$
\item If $B\in X_{fin}$ with $B$ is minimal, and (\ref{4.3}) holds, by the minimality of $B$ there exists $i_0\in\{1,\ldots,n\}$, such that $B_{i_0}=B$. So, the others $B_i$'s cannot be infinite emitters nor have infinite cardinality, since the union in (\ref{4.3}) is disjoint. This implies that $F_i=\emptyset$ and $S_i=\emptyset$ for all $i\neq i_0$. Furthermore, by (\ref{4.3}), $F_{i_0}\supseteq F$, $S_{i_0}\supseteq S$. Thus $F_{i_0}\backslash F=\displaystyle\bigsqcup_{i\neq i_0} \varepsilon(B_i)$ and  $S_{i_0}\backslash S=\displaystyle\bigsqcup_{i\neq i_0} B_i\cap G_s^0$. Therefore, we have:
\begin{align*}
\kappa\left(D_{(\beta,B),F,S}\right)=M(\beta)m(B)-&\displaystyle\sum_{e\in F} M(\beta e)m(r(e))-\displaystyle\sum_{v\in S} M(\beta )m(v)\\
=M(\beta)m(B)-\displaystyle\sum_{e\in F_{i_0}} M(\beta e)m(r(e))+&\displaystyle\sum_{e\in F_{i_0}\backslash F} M(\beta e)m(r(e))-\\
-\displaystyle\sum_{v\in S_{i_0}} M(\beta )m(v)+&\displaystyle\sum_{v\in S_{i_0}\backslash S} M(\beta )m(v)\\
=M(\beta)m(B)-\displaystyle\sum_{e\in F_{i_0}} M(\beta e)m(r(e))+&\displaystyle\sum_{i=1,i\neq i_0}^n\displaystyle\sum_{e\in \varepsilon(B_i)} M(\beta e)m(r(e))-\\
-\displaystyle\sum_{v\in S_{i_0}} M(\beta )m(v)+&\displaystyle\sum_{i=1,i\neq i_0}^n\displaystyle\sum_{v\in B_i\cap G_s^0} M(\beta)m(v)\\
=M(\beta)m(B)-\displaystyle\sum_{e\in F_{i_0}} M(\beta e)m(r(e))-&\displaystyle\sum_{v\in S_{i_0}} M(\beta )m(v)+\\ +\displaystyle\sum_{i=1,i\neq i_0}^n&\left(\displaystyle\sum_{e\in \varepsilon(B_i)} M(\beta e)m(r(e))+\displaystyle\sum_{v\in B_i\cap G_s^0} M(\beta)m(v)\right)\\
=\kappa (D_{(\beta,B_{i_0},F_{i_0},S_{i_0})}) +\displaystyle\sum_{i=1,i\neq i_0}^n M(\beta)&\left(\displaystyle\sum_{e\in \varepsilon(B_i)} M( e)m(r(e))+\displaystyle\sum_{v\in B_i\cap G_s^0} m(v)\right)\\
=\kappa (D_{(\beta,B_{i_0},F_{i_0},S_{i_0})}) +&\displaystyle\sum_{i=1,i\neq i_0}^n M(\beta)M(B_i)\\=\kappa (D_{(\beta,B_{i_0},F_{i_0},S_{i_0})})+&\displaystyle\sum_{i=1,i\neq i_0}^n\kappa(D_{(\beta,B_i)})=\displaystyle\sum_{i=1}^n\kappa(D_{(\beta,B_i)},F_i,S_i).
\end{align*}
Note that for all $i\neq i_0$, we use that $\kappa(D{(\beta,B_i)})=\kappa(D{(\beta,B_i),F_i,S_i})=M(\beta)m(B_i)$, since $F_i=\emptyset$ and $S_i=\emptyset$ for all $i\neq i_0$. 
\end{itemize}

Now, suppose that the previous equality holds for all $0\leq k< m$. 
Again we will do the proof for the case where $(\beta,B)\in X_{fin}$ with $B$ minimal, and the other case is analogous to what was done when $m=0$. Using the minimality of $B$ once again, there exist $i_0$ such that $\beta_{i_0}=\beta$, $B_{i_0}=B$, $F_{i_0}\supseteq F$ and $S_{i_0}\supseteq S$. Suppose, without loss of generality, that $i_1$ is such that $\beta_i=\beta$ whenever $1\leq i \leq i_1$, and $\beta_i\neq \beta$ whenever $i_1 < i \leq n$. By the case $m=0$ seen previously, $F_i=\emptyset$ and $S_i=\emptyset$ for all $1\leq i \leq i_1$, $i\neq i_0$. Furthermore, $A_i$ is neither an infinite emitter nor has infinite cardinality for any $1\leq i \leq i_1$, $i\neq i_0$. Also $\displaystyle\bigsqcup_{i\neq i_0, i=1}^{i_1}B_i\cap G_s^0= S_{i_0}\setminus S$, however $\displaystyle\bigsqcup_{i\neq i_0, i=1}^{i_1}\varepsilon(B_i) \subseteq F_{i_0}\setminus F$. Define $C=s\left((F_{i_0}\setminus F)\Big\backslash \displaystyle\bigsqcup_{i\neq i_0,i=1}^{i_1}\varepsilon(B_i)\right)$. Of course $|\varepsilon(C)|<\infty$. Note that
\begin{center}
$D_{(\beta,B),F,S}=D_{(\beta,C)}\displaystyle\sqcup\bigsqcup_{i=1}^{i_1} D_{(\beta_i,B_i),F_i,S_i}$,
and that
$D_{(\beta,C)}=\displaystyle\bigsqcup_{i=i_1+1}^{n} D_{(\beta_i,B_i),F_i,S_i}.$
\end{center}
By the case $m=0$, we have that 
\begin{equation}\label{eq:sum0}
\kappa(D_{(\beta,B),F,S})=\kappa (D_{(\beta,C)})+\sum_{i=1}^{i_1} \kappa (D_{(\beta_i,B_i),F_i,S_i}).
\end{equation}
Since $|\varepsilon(C)|<\infty$ we have that  $D_{(\beta,C)}=\displaystyle\bigsqcup_{s(e)\in C} D_{(\beta e,r(e))}$ and the union is finite. On the other hand, by Corollary~\ref{cor_78}, m2. and m4., we can write the finite union 
$D_{(\beta,C)}=\bigsqcup_j D_{(\beta e_j, A_j)},$
where each $D_{(\beta e_j, A_j)}\in\widehat{S}$, in such a way that
\begin{equation} \label{eq:sum1}\kappa(D_{(\beta,C)})=\sum_j\kappa( D_{(\beta e_j, A_j)}).\end{equation}
Thus, $D_{(\beta e_j, A_j)}=D_{(\beta e_j, A_j)}\cap D_{(\beta,C)}=\displaystyle\bigsqcup_{i=i_1+1}^{n}D_{(\beta e_j, A_j)}\cap D_{(\beta_i,B_i),F_i,S_i},$
and note that $D_{(\beta e_j, A_j)}\cap D_{(\beta_i,B_i),F_i,S_i} = D_{(\gamma_j,D_j)}$ for some element of $\widehat{S}$ where  $|\gamma_j|=|\beta_i|$ (because $|\beta_i|\geq |\beta|+1=|\beta e_j|$). This implies that $\max_j\{|\gamma_j|-|\beta e_i|\}<m$. By the induction hypothesis,
\begin{equation} \label{eq:sum2}\kappa(D_{(\beta e_j, A_j)})=\sum_{i=i_1+1}^{n}\kappa(D_{(\beta e_j, A_j)}\cap D_{(\beta_i,B_i),F_i,S_i}).\end{equation}

On the other hand, for $i_1+1\leq i\leq n$, 
$D_{(\beta_i,B_i),F_i,S_i}=D_{(\beta_i,B_i),F_i,S_i}\cap D_{(\beta,C)}=\displaystyle\bigsqcup_j D_{(\beta e_j, A_j)}\cap D_{(\beta_i,B_i),F_i,S_i},$
again using the case $m=0$ (even if any intersections are empty) we have \begin{equation} \label{eq:sum3}\kappa(D_{(\beta_i,B_i),F_i,S_i})=\sum_j \kappa(D_{(\beta e_j, A_j)}\cap D_{(\beta_i,B_i),F_i,S_i}). \end{equation}

Putting the equations (\ref{eq:sum1}), (\ref{eq:sum2}) and (\ref{eq:sum3}) together we conclude that
\\ $\kappa(D_{(\beta,C)})=\displaystyle\sum_{i=i_1+1}^{n} \kappa(D_{(\beta_i,B_i),F_i,S_i}),$
and so, by (\ref{eq:sum0}), we have that $\kappa$ is in fact a measure.

It remains to be shown that $\kappa (\theta_e(V))= M(e)\kappa(V)$. 
But note that if $\widehat{S}\ni V=D_{(\beta,B),F,S}$ and $V\subseteq X_{e^{-1}}$, then  $\theta_e(V)=D_{(e\beta,B),F,S}$. Thus, 
\begin{align*}
\kappa(\theta_e(V))&=\kappa(D_{(e\beta,B),F,S})
=M(e\beta)m(B)-\sum_{e\in F} M(e\beta f)m(r(f))-\sum_{v\in S\cap G_s^0} M(e\beta )m(v)\\
&=M(e)M(\beta)m(B)-\sum_{e\in F} M(e)M(\beta f)m(r(f))-\sum_{v\in S\cap G_s^0} M(e)M(\beta )m(v)\\
&=M(e)\left(M(\beta)m(B)-\sum_{e\in F} M(\beta f)m(r(f))-\sum_{v\in S\cap G_s^0} M(\beta )m(v)\right)\\
&=M(e)\kappa(D_{(\beta,B),F,S})=M(e)\kappa(V).
\end{align*}
and the result follows.
\end{proof}

Using Carathéodory extension theorem we can extend the measure above to a unique measure $\mu$ defined on the $\sigma$-algebra generated by $\widehat{S}$. Of course this measure is related to the functions $M$ and $m$ defined above. We make this relation precise below. 

\begin{prop}\label{prop_4.20}
Let $\G$ be an ultragraph that satisfies Condition~(RFUM2) and $M:\G^1\rightarrow [0,1]$. Then there is a convex map between the set of functions $m:\G^0\rightarrow [0,1]$ that satisfies m1., ..., m4. of page \pageref{emes} and the set of regular, Borel, probability measures $\mu$ on $X$ satisfying $\mu(\theta_e(V))=M(e)\mu(V)$, for all $e\in \G^1$ and for all Borel measurable subsets $V\subseteq X_{e^{-1}}$.
\end{prop}
\begin{proof}
By Lemma~\ref{lema_medida}, we have a measure $\kappa$ in the semi-ring $\widehat{S}$. Hence, by the Carathéodory extension theorem, there is a unique measure $\mu$ defined on the $\sigma$-algebra generated by $\widehat{S}$. Since $\widehat{S}$ forms a countable base, $\mu$ is defined in the Borel $\sigma$-algebra.

By definition of $\kappa$, we have that $\mu(D_{(A,A)})=\kappa(D_{(A,A)})=m(A)$, for all $D_{(A,A)}\in S$. Hence, considering that $\mu$ is measure, $\G^0 $ is a directed set via inclusion, and using the principle of inclusion exclusion, together with m4., we have
$$\mu(X)=\mu\left(\displaystyle\bigcup_{A\in\G^0} D_{(A,A)}\right)=\sum_{A_\in\G^0} \mu(D_{(A,A)})=\sum_{A_\in\G^0}m(A)=\lim_{A\in \G^0} m(A)\overset{m1.}{=}1.$$ 
Thus, $\mu$ is probability measure, which also implies that it is regular, since every finite Borel measure is regular. Also, a convex combination of functions $m$ is preserved when we move to measures $\kappa$ in $\widehat{S}$ (and therefore to measures $\mu$ in $X$), due to Equation (\ref{4.3}).

Finally, to prove that $\mu(\theta_e(V))=M(e)\mu(V)$ for all Borel measurable $V\subseteq X_{e^{-1}}$ and $e\in \G^1$ we have, it is enough to check the equality for the elements of the semi-ring $\widehat{S}$, since $\theta_e$ is a bijection and preserves unions and intersections. But this was done in the proof of Lemma~\ref{lema_medida}.
\end{proof}

As we mentioned before, Riesz representation theory connects states with measures. We make a precise statement below. 

\begin{prop}\label{prop_4.21}
Let $\G$ be an ultragraph that satisfies Condition~(RFUM2) and $M:\G^1\rightarrow [0,1]$. Then there is an affine bijection between the set of states $\omega$ in $C_0(X)$ that satisfy the condition $\omega(f\circ\theta_e^{-1})=M(e) \omega(f)$, for every $e\in\G^1$ and every $f\in C_0(X_{e^{-1}})$, and the set of  regular, Borel, probability measures $\mu$ on $X$ that satisfy $\mu(\theta_e(V))=M(e)\mu(V)$, $\forall e\in\G^1$ and for every Borel measurable subset $V\subseteq X_{e^{-1}}$.
\end{prop}
\begin{proof}
The proof follows from the Riesz representation theorem, and is analogous to the proof of \cite[Lema 4.7]{CarLar} and \cite[Prop. 4.8]{danielgilles}.
\end{proof}

\subsection*{Description of the KMS states of $C^*(\G)$}

We now have all the ingredients to characterize KMS states in five different ways. Before we do this we prove one last result, which is a generalization of \cite[Proposition 4.4]{danielgilles} and  \cite[Lemma 4.9]{CarLar}.
\begin{prop}\label{prop_4.22}
	Let $\G$ be an ultragraph that satisfies Condition~(RFUM2) and $M:\G^1\rightarrow [0,1]$. Then there is a convex, injective map between the set of functions $m:\G^0\rightarrow [0,1]$ satisfying m1., ..., m4. of page \pageref{emes}, and the set of states $\omega$ on $C_0(X)$ that satisfies the condition $\omega(f\circ \theta_e^{-1})=N(e)^{-\beta}\omega(f)$, for all $e\in \G^1$ and all $f\in C_0(X_{e^{-1}})$. Furthermore, such correspondence is of the form $\omega\longmapsto(A\overset{m}{\mapsto} \omega(1_A))$, $\forall A\in\G^0$.
\end{prop}
\begin{proof}
Let $\omega$ be a state on $C_0(X)$ such that $\omega(f\circ\theta_{e^{-1}})=M(e)\omega(f)$, for all $e\in\G^1$ and $f\in C_0(X_{e^{-1}})$. Let $m:\G^0\longmapsto [0,1]$ be given by $m(A)=\omega(1_A)$. We check that $m$ satisfies the conditions $mi.$'s. 
\begin{enumerate}[\textbf m1.]
\item The set $\{A:A\in \G^0\}$ is directed via inclusion, hence $\{1_A:A\in\G^0\}$ is an increasing approximate unity for $C_0(X)$. Thus, since $\omega$ is state, $1=\displaystyle\lim_{A\in\G^0} \omega(1_A)=\lim_{A\in\G^0}m(A)$.
\item First note that if $e\in\G^1$ then $$\omega(1_e)=\omega(1_{e^{-1}}\circ\theta_{e^{-1}})=M(e)\omega(1_{e^{-1}})=M(e)\omega(1_{r(e)})=M(e)m(r(e)).$$
Furthermore, note that if $|A|<\infty$ and $|\varepsilon(A)|<\infty$ then, by Corollary~\ref{cor_78}, we have that $$ 1_A=\displaystyle\sum_{e:s(e)\in A}^{finite}1_e+\sum_{v\in A\cap G^0_{s}}^{finite} 1_v.$$ 
Thus
\begin{align*}
m(A)=\omega(1_A)&=\omega\left(\displaystyle\sum_{e:s(e)\in A}1_e+\sum_{v\in A\cap G^0_{s}} 1_v\right)=\displaystyle\sum_{e:s(e)\in A}\omega(1_e)+\sum_{v\in A\cap G^0_{s}} \omega(1_v)\\
&=\displaystyle\sum_{e:s(e)\in A}M(e)m(r(e))+\sum_{v\in A\cap G^0_{s}} m(v).
\end{align*}
\item Let $A\in\G^0$, $F$ be a finite subset of $\varepsilon(A)$ and $S$ be a finite subset of $A\cap G_s^0$. Then, again by Corollary~\ref{cor_78}, we have that $1_A\geq \displaystyle\sum_{e\in F} 1_e
+\sum_{v\in A\cap G^0_{s}} 1_v$. Hence, $$m(A)=\omega(1_A)\geq\omega\left(\displaystyle\sum_{e\in F} 1_e
+\sum_{v\in A\cap\G^0_{sink}} 1_v\right)=\sum_{e\in F}M(e)m(r(e))+\sum_{v\in G_s^0}m(v).$$
\item Follows from the linearity of $\omega$.
\end{enumerate}

Therefore, for a state $\omega$ on $C_0(X)$ we assign a function $m:\G^0\rightarrow [0,1]$ that satisfies $m1.,\ldots,m4.$. To prove the injectivity of such correspondence, suppose that there are  $\omega_1$ and $\omega_2$ states such that $\omega_i(f\circ\theta_{e^{-1}})=M(e)\omega_i(f)$, $i=1,2$, $\forall f\in C_0(X_{e^{-1}})$ and $e\in \G^1$. Suppose further that  $\omega_1(1_A)=\omega_2(1_A)$, $\forall A\in \G^0$. We must show that $\omega_1=\omega_2$. By Remark~\ref{obs81}, we just need to show that $\omega_1(1_{r(c)\cap A}\circ \theta_{c^{-1}})=\omega_2(1_{r(c)\cap A}\circ \theta_{c^{-1}})$, for every $c\in \G^*\backslash \G^0$ and every $A\in\G^0$. To do this, first note that $$1_A\circ\theta_{ef}^{-1}=1_A\circ\theta_{f^{-1}e^{-1}}=1_A\circ\theta_{f^{-1}}\circ\theta_{e^{-1}}=1_{\theta_f(A)}\circ\theta_{e^{-1}}.$$
Therefore, if $c=c_1\ldots c_n\in \G^*$ and $A\in\G^0$, we have:
\begin{align*}
\omega_1\left(1_{r(c)\cap A}\circ \theta_{c^{-1}}\right)&=\omega_1\left(1_{r(c)\cap A}\circ \theta_{c_n^{-1}\ldots c_1^{-1}}\right)=\omega_1\left(1_{\theta_{c_n}(r(c)\cap A)}\circ \theta_{c_{n-1}^{-1}\ldots c_1^{-1}}\right)=\ldots\\
&=\omega_1\left(1_{\theta_{c_2\ldots c_n}(r(c)\cap A)}\circ \theta_{c_1^{-1}}\right)=M(c_1) \omega_1\left(1_{\theta_{c_2\ldots c_n}(r(c)\cap A)}\right)=\\
&=M(c_1) \omega_1\left(1_{\theta_{c_3\ldots c_n}(r(c)\cap A)}\circ \theta_{c_2}^{-1}\right)=M(c_1) M(c_2)\omega_1\left(1_{\theta_{c_4\ldots c_n}(r(c)\cap A)}\circ \theta_{c_3}^{-1}\right)=\ldots\\
&=M(c_1) \ldots M(c_n)\omega_1\left(1_{r(c)\cap A}\right)=M(c_1) \ldots M(c_n)\omega_2\left(1_{r(c)\cap A}\right)=\ldots=\\
&=\omega_2\left(1_{r(c)\cap A}\circ \theta_{c^{-1}}\right).
\end{align*}
\end{proof}
We can now state and prove our main result regarding KMS states of C*-algebras associated with ultragraphs that satisfy Condition~(RFUM2). This generalizes \cite[Theorem 4.3]{danielgilles} and  \cite[Theorem 4.1]{CarLar}.

\begin{teorema}\label{teo_do_cap}
Let $\G$ be an ultragraph that satisfies Condition~(RFUM2). Consider the function $N:\G^1\rightarrow (1,\infty)$ and $0\leq \beta<\infty$, and consider the following sets:  
\begin{itemize}
\item[$A^\beta$:] the set of KMS$_\beta$ states for $C^*(\G)$;
\item[$B^\beta$:] the set of states $\omega$ of $C_0(X)$ that satisfy the condition $\omega(f\circ \theta_e^{-1})=N(e)^{-\beta}\omega(f)$, for all $e\in\G^1$ and all $f\in C_0(X_{e^{-1}})$;
\item[$C^\beta$:] The set of regular, Bore probability measure $\mu$ on $X$ that satisfy $\mu(\theta_e(V))=N(e)^{-\beta}\mu(V)$, for all $e\in \G^1$ and for all Borel measurable subset $V\subseteq X_{e^{-1}}$;
\item[$D^\beta$:] The set of functions $m:\G^0\rightarrow [0,1]$ that satisfy
\begin{enumerate}[\textbf m1.]
	\item $\displaystyle\lim_{A\in\G^0} m(A)=1$;
	\item $m(A)=\displaystyle\sum_{e:s(e)\in A}N(e)^{-\beta}m(r(e))+\sum_{v\in A\cap G_{s}^0}m(v)$, if $|A|<\infty$ and $|\epsilon(A)|<\infty$;
	\item $m(A)\geq \displaystyle\sum_{e\in F}N(e)^{-\beta}m(r(e)) +\displaystyle\sum_{v\in S}m(v)$, for every  finite subset $F\subseteq \varepsilon(A)$ and every finite subset $S\subseteq A\cap G_s^0$, $S$;
	\item $m(A\cup B)=m(A)+m(B)-m(A\cap B)$.
\end{enumerate}
\item[$E^\beta$:] The set of states $\psi$ on $C^*(\G)^\gamma$ that satisfies $\psi(s_cp_As_d^*)=\delta_{c,d}N(c)^{-\beta}\psi(p_{A\cap r(c)})$, for all $c,d\in\G^*\backslash\G^0$ and $A\in\G^0$.
\end{itemize}
Then there exists a convex isomorphism between $A^{\beta},B^{\beta},C^{\beta},D^{\beta}$ and $E^{\beta}$.
\end{teorema} 
\begin{proof}
By Proposition \ref{prop_4.9} there is an affine isomorphism between $A^{\beta}$ and $E^\beta$. From proposition~\ref{prop_4.13} we get an affine isomorphism between $B^\beta$ and $E^\beta$. 

We use the Proposition~\ref{prop_4.21} to show the affine isomorphism between $B^\beta$ and $C^\beta$. For this, note that there is a correspondence between the functions $N:\G^1\rightarrow (1,\infty)$ and the functions $M:\G^1\rightarrow [0,1]$ (just let $M(e)=N(e)^{-\beta}$).  

We already proved that $A^\beta\iff E^\beta\iff B^\beta \iff C^\beta$. It remains to show the equivalence with $D^\beta$. Let $M(e)=N(e)^{-\beta}$. From Proposition~\ref{prop_4.22} we obtain an affine map from $B^\beta$ to $D^\beta$ and, taking the same $M$ in Proposition~\ref{prop_4.20}, we obtain an affine map from $D^\beta$ to $C^\beta$. Finally, the desired result follows from the fact that the maps from $B^\beta$ to $D^\beta$, from $D^\beta$ to $C^\beta$, and from $C^\beta$ to $B^\beta$, compose to the identity. 
\end{proof}

\subsection*{Ground States}\label{ground}

In this last section we characterize ground states in ultragraph C*-algebras that satisfy Condition~(RFUM2). Our results extend \cite[Theo. 5.1]{danielgilles}.

Let $\G$ be an ultragraph that satisfies Condition~(RFUM2). Define
\begin{itemize}
	\item[$A^{gr}$:] the set of ground states in $C^*(\G)$,
	\item[$B^{gr}$:] the set of states $\omega$ of $C_0(X)$ such that $\omega(1_e)=0 $, for every $e\in \G^1$,
	\item[$C^{gr}$:] the set of all regular, Borel probability measure $\mu$ on $X$ such that $\mu(A) = 0$, for every $e\in \G^1$ and every Borel measurable subset $A$ of $X_{e}$,
	\item[$D^{gr}$:] the set of functions $m:\G^0\to [0,1]$ satisfying 
 \begin{enumerate}[\textbf m1.]
	\item $\displaystyle\lim_{A\in\G^0} m(A)=1$;
	\item $m(A)=0$ if $|A|<\infty$ and $|\epsilon(A)|<\infty$;
	\item $m(A\cup B)=m(A)+m(B)-m(A\cap B)$.
\end{enumerate}
\end{itemize}

We now have the following.

\begin{teorema}\label{teoremaground} Let $\G$ be an ultragraph that satisfies Condition~(RFUM2). Then there is an affine isomorphism between $A^{gr}$, $B^{gr}$, $C^{gr}$, and $D^{gr}$.
	
\end{teorema}

\begin{proof}The proof follows the same lines of \cite[Theo. 5.1]{danielgilles}. The existence of an affine isomorphism between $A^{gr}$ and the set of the states $\phi$ on $C_0(X)$ such that $\phi(f)=0 $, for all $e\in \G^1$ and $f\in C(X_e)$, follows from \cite[Theorem 4.3]{exel_laca}. Since $\phi$ is a state, and $1_e$ is unity for $C(X_e)$, it follows that if $\phi(1_e)=0$ then $\phi(f)=0$ for all $f\in C(X_e)$. Hence we have that $A^{gr}$ is isomorphic to $B^{gr}$, via an affine isomorphism.

As with KMS states, an isomorphism between $B^{gr}$ and $C^{gr}$ is obtained analogously to what was done in \cite[Proposition 4.8]{CarLar}.

Finally, an affine isomorphism between $B^{gr}$ and $D^{gr}$ is obtained by applying Propositions~\ref{prop_4.20} and \ref{prop_4.22} with $M(e)=0$ for all $e\in \G^1$, and proceeding as in the proof of Theorem \ref{teo_do_cap}.
\end{proof}

\begin{obs}
In \cite[Obs 4.4]{CarLar} the authors related the description of KMS states of graph C*-algebra via partial crossed products and via Deaconu-Renault groupoids. Since we realized C*-algebras of ultragraphs that satisfy Condition~(RFUM2) both as a partial crossed product and as a groupoid C*-algebra an analogous observation holds in this case.
\end{obs}

We finish the paper describing the KMS and ground states of the C*-algebra associated to an ultragraph with sinks. Furthermore, for the example we consider, it is known that the ultragraph Leavitt path algebra associated to it (over the field $\Z_2)$ is neither a Leavitt path algebra nor a algebraic Exel-Laca algebra (see \cite[Example 6.12]{leavittultragraph} and \cite{ultrapartialleavitt}.
\begin{exemplo}
Consider the following\label{ultragrafo_exemplo_geral} ultragraph $\G$, which has only one edge:

\hspace{3cm}
\setlength{\unitlength}{1mm}
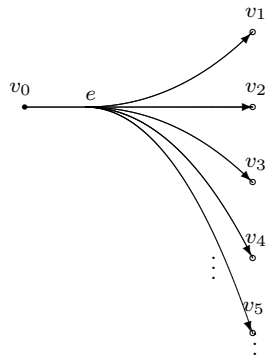
\begin{picture}(80,45)
	\put(0,32){\scriptsize$v_0$}
	\put(2,30){\circle*{0.7}}
	\put(31,42){\scriptsize$v_1$}
	\put(32,40){\circle{0.8}}
	\put(31,32){\scriptsize$v_2$}
	\put(32,30){\circle{0.8}}
	\put(31,22){\scriptsize$v_3$}
	\put(32,20){\circle{0.8}}
	\put(31,12){\scriptsize$v_4$}
	\put(32,10){\circle{0.8}}
	\put(30.7,3){\scriptsize$v_5$}
	\put(32,0){\circle{0.8}}
	\put(10,31){\scriptsize$e$}

	\linethickness{0.2mm}
	\put(2,30){\line(1,0){10}}
	\thinlines
	\qbezier(10,30)(22,30)(32,40)
	\qbezier(10,30)(22,30)(32,30)
	\qbezier(10,30)(22,30)(32,20)
	\qbezier(10,30)(22,30)(32,10)
	\qbezier(10,30)(22,30)(32,0)
	\put(31.8,39.8){\vector(1,1){0}}
	\put(32,30){\vector(1,0){0}}
	\put(31.8,20.2){\vector(1,-1){0}}
	\put(31.8,10.2){\vector(1,-2){0}}
	\put(31.8,0.2){\vector(1,-3){0}}
	\begin{rotate}{90}
	\put(7,-27){\scriptsize$\ldots$}
	\end{rotate}
	
	\begin{rotate}{90}
	\put(-3.2,-31){\scriptsize$\ldots$}
	\end{rotate}

\end{picture}
	
\vspace{10pt}

Note that $\G^0$ is formed by finite subsets of $G^0$, in addition to finite unions with $r(e)=\{v_1,v_2,v_3,\ldots\}$. In particular, $G^0\in\G^0$. The shift space associated is
$$X=\{(e,v_i)\}_{i\geq 1}\cup \{(v_i,v_i)\}_{i\geq 1}\cup \{(e,r(e)),(r(e),r(e))\}.$$

Our goal is to find $m$ that satisfies the conditions $\textbf{m}1$ to $\textbf{m}4$ that define the set $D^\beta$. For such a function to exist we must have, for $N$ and $\beta$ given, that:
\begin{equation}\label{eq.example}
m(\{v_0\})+m(r(e))=m(G^0)=1.
\end{equation}
Also, $m$ must satisfy \textbf{m2}, \textbf{m4}, and since $s(e)\in \{v_0\}$, we must have 
\begin{equation}\label{eq.example2}
m(\{v_0\})=N(e)^{-\beta}m(r(e)),
\end{equation}
i.e, $m(\{v_0\})$ depends on $N$, $\beta$ and of $r(e)$.

Now, joining $(\ref{eq.example})$ and $(\ref{eq.example2})$, we have that $m(r(e))=\displaystyle\frac{1}{1+N(e)^{-\beta}}$, i.e $m(r(e))$ depends only on $N$ and $\beta$.

For \textbf{m3} to be satisfied we must have that, for all $k\in \N^*$,
$ m(r(e))\geq \displaystyle\sum_{i=1}^k m(v_i)$, i.e
$$\displaystyle\frac{1}{1+N(e)^{-\beta}}\geq \sum_{i=1}^k m(v_i)$$

By the definition of $N$ we have that $N(e)\in (1,\infty)$, and hence $N(e)^{-\beta}\in (0,1)$ (thus, if $\beta\longrightarrow \infty$ then $N(e)^{-\beta}\longrightarrow 0$). Therefore, for each $i\in \N^*$, consider 
\[m(\{v_i\})=\frac{1-N(e)^{-\beta}}{1+N(e)^{-\beta}}\left(N(e)^{-\beta}\right)^{(i-1)}.\]
 
To make equations more readable, let $r:=N(e)^{-\beta}$ and $a:=\frac{1-N(e)^{-\beta}}{1+N(e)^{-\beta}}=\frac{1-r}{1+r}$. Then

\[ m(r(e))= \sum_{i=1}^\infty m(v_i) = \sum_{i=1}^\infty a.r^{(i-1)}=a\sum_{i=1}^\infty r^{(i-1)}=a\left(\frac{1}{1-r}\right) =\left(\frac{1-r}{1+r}\right).\left(\frac{1}{1-r}\right)=\frac{1}{1+r}.\]
Furthermore, since $1\geq a$, we have that for all $k\in\N^*$  it holds that

\[ m(r(e))=\frac{1}{1+r}\geq \frac{1-r^k}{1+r}\geq \left(\frac{1-r}{1+r}\right)\frac{1-r^k}{1+r}=\sum_{i=1}^k m(v_i),\]
and so \textbf{m3} is valid for any finite subset $S$ of $r(e)$. 

Also notice that, independently on $N$ and  $\beta$, if we choose $m(\{v_0\})=1$ and $m(r(e))=m(\{v_i\})=0$ for $i\in \N^*$, we have that \textbf{m1-m4} is satisfied.

Finally, for ground states, the function $m:\G^0\rightarrow [0,1]$ is in $D^{gr}$ if, and only if, $m(\{v_i\})=0$ for all $i\in \mathbb{N}$, and $m(r(e))=1$. 

\end{exemplo}


\end{document}